\documentclass{article}

\usepackage{amssymb}
\usepackage[utf8]{inputenc}
\usepackage{graphicx}
\usepackage{authblk}
\usepackage{amsfonts} 
\usepackage{mathtools} 
\usepackage{amsmath}
\usepackage{amsthm} 
\usepackage{mathrsfs} 
\usepackage{xcolor}
\usepackage{bbm}
\usepackage{bm}
\usepackage{pifont}
\usepackage[T1]{fontenc}
\usepackage[framemethod=tikz]{mdframed}
\usepackage{longtable}
\usepackage{indentfirst} 
\usepackage{extarrows}
\usepackage{centernot}
\usepackage{hyperref}
\usepackage[hang,flushmargin]{footmisc}
\usepackage[bottom=2.0cm,top=2.0cm,left=1.5cm,right=1.5cm]{geometry}
\usepackage[style=numeric,defernumbers,backend=biber,giveninits=true,sorting=nyt]{biblatex}

\usepackage{tikz}
\usetikzlibrary{arrows.meta, positioning, matrix}
\usetikzlibrary{calc}

\newcommand{\RomanNumeralCaps}[1]{\MakeUppercase{\romannumeral #1}}

\hypersetup{colorlinks=true,linkcolor=ccr,anchorcolor=ccr,citecolor=ccr}

\DeclareUnicodeCharacter{02B9}{'}
\DeclareUnicodeCharacter{FE20}{'}
\DeclareUnicodeCharacter{FE21}{'}

\definecolor{ccr}{RGB}{10,110,150}
            
\title{From CKLS Process to CIR-type Process: Using a Twice-differentiable Mapping and Generalized Girsanov's Theorem}
\author[1]{Boyuan Ning}
\author[2]{Yasutaka Shimizu}
\affil[1]{ningboyuan@akane.waseda.jp}
\affil[2]{shimizu@waseda.jp}
\affil[1,2]{Department of Pure and Applied Mathematics, Waseda University, \authorcr 3-4-1 Ohkubo, Shinjuku-ku, Tokyo, 169-8555, Japan}
\date{October 6th, 2026}

\begin{document}
\numberwithin{equation}{section}

\maketitle

\begin{abstract}
    We construct a twice-differentiable mapping $T(x): \mathbb{R}_{+}\to\mathbb{R}_{+}$ satisfying $\frac{dT(x)}{dx}x^{k}=L\big[T(x)\big]^{\frac{1}{2}}$ for a given constant $L$ and apply it to the CKLS short-rate process $\lambda_{t}$, which solves the stochastic differential equation (SDE) of the form $d\lambda_{t}=(a-b\lambda_{t})dt+\sigma (\lambda_{t})^{k}dW_{t}$. By It\^{o}'s lemma,the transformed process $X_{t}\overset{\text{def}}{=}T(\lambda_{t})$ obeys an SDE whose diffusion term is proportional to $(X_{t})^{\frac{1}{2}}$ and whose drift is a non-linear function of $\lambda_{t}$. A critical review of an earlier study on the same transformation reveals substantial errors in its model specification, derivations, and proofs. Next, a generalized Girsanov transformation of measure is introduced to shift the drift. Under the equivalent measure $\mathbb{Q}$, the dynamics of $X_{t}$ reduces to the classical Cox-Ingersoll-Ross (CIR) form. Leveraging well-known properties concerning uniqueness, strongness, and positivity of $\lambda_{t}$ induced by the Yamada-Watanabe-Engelbert theorem, we show that the combined twice-differentiable mapping and Girsanov step is valid precisely when $L>0$, $a>0, b>0, \sigma>0$ and, most importantly, $\frac{1}{2}<k<1$ (which is the parameter range of particular relevance in financial applications) or $k=\frac{1}{2}$ with $2a\geq\sigma^{2}$ (which reduces to the CIR process with Feller's condition satisfied). The CIR representation allows us to obtain the transition density, moment formulas, stationary density, and boundary behavior of $X_{t}$, and to derive the transition density of $\lambda_{t}$ under $\mathbb{Q}$. Finally, we explain why the available moment estimates alone do not establish Novikov's or Kazamaki's condition, and then prove directly that the Dol\'{e}ans-Dade exponential associated with our Girsanov transformation is a true martingale on every finite time interval (thus can be called Radon-Nikod\'{y}m derivative), thus we have the soundness of the entire procedure combining $T(x)$ and the Girsanov transformation validated. Our argument adapts a recent result, rather than relying on Novikov's or Kazamaki's conditions, that extends the classical martingale criterion: by applying Feller's explosion test together with his boundary classification, it provides a necessary and sufficient condition under which the Radon-Nikod\'{y}m derivative is a true martingale.\\
   
    
\end{abstract}

\section{Introduction}
\subsection{The Chan-Karolyi-Longstaff-Sanders model and some of its properties}
\noindent In single-factor models, the evolution of the short rate can be given by a stochastic differential equation (SDE) defined on some filtered probability space $(\Omega,\mathcal{F},\{\mathcal{F}_{t}\}_{t\in[0,T]},\mathbb{P})$:
\begin{equation}
    d\lambda_{t}=(a-b\lambda_{t})dt+\sigma (\lambda_{t})^{k}dW_{t},
\end{equation}
where $W_{t}$ is a Wiener process on the given probability space. The parameters include the constant initial value $\lambda_{0}\in\mathbb{R}_{+}$, the drift intercept $a\in\mathbb{R}_{+}$ (standing for long-term mean level times mean-reversion speed level) and the mean-reversion speed level $b\in\mathbb{R}_{+}$ constituting the drift term, the volatility (diffusion) $\sigma\in\mathbb{R}_{+}$ and the elasticity (of volatility) $k\in\mathbb{R}_{+}$ constituting the volatility (diffusion) term. This model is commonly known as the Chan-Karolyi-Longstaff-Sanders (CKLS) model, which was first proposed by \hyperlink{Chan1992}{Chan et al. (1992)} to model the short-term interest rate. The stochastic process $\lambda_{t}$, which is the solution to this SDE, is generally called the CKLS process.\\

\noindent \textbf{Remark 1.1}\\
\noindent (1) The CKLS model, as presented in Equation (1.1), describes a broad range of interest rate processes, encompassing several well-known interest rate models:
\begin{longtable}{|p{.6\textwidth}|p{.09\textwidth}|p{.09\textwidth}|p{.09\textwidth}|}
\caption*{Variants of CKLS model under different parametric specifications}\\
\hline
\textbf{Model/Process} & $a$ & $b$ & $k$ \\
\hline
\endfirsthead
\hline
\endhead
\hline
\endfoot
\hline
\endlastfoot
Merton (\hyperlink{Merton1974}{Merton, 1974})& Any & 0 & 0 \\ 
\hline
Va\v{s}\'{\i}\v{c}ek \hyperlink{Vasicek1977}{(Va\v{s}\'{\i}\v{c}ek, 1977)} & Any & Any & 0 \\
\hline
Cox-Ingersoll-Ross (CIR) (\hyperlink{Cox1985}{Cox et al., 1985}) & Any & Any & 1/2 \\
\hline
Dothan (\hyperlink{Dothan1978}{Dothan, 1978}) & 0 & 0 & 1 \\
\hline
Geometric Brownian motion & 0 & Any & 1 \\
\hline
Brennan and Schwartz (\hyperlink{Brennan1980}{Brennan and Schwartz, 1980}) & Any & Any & 1 \\
\hline
Cox-Ingersoll-Ross Variable-Rate (CIR VR) (\hyperlink{Cox1989}{Cox et al., 1980}) & 0 & 0 & 3/2 \\
\hline
Constant Elasticity of Variance (CEV) (\hyperlink{Cox1996}{Cox, 1996}) & 0 & Any & Any \\
\end{longtable}
\noindent (2) Many scholars, particularly in financial fields, may treat the CKLS model as a generalization of the CEV model and name it "Mean-reverting CEV model" (e.g. by \hyperlink{Tsumurai2020}{Tsumurai (2020)}), "CEV model with linear drift" (e.g. by \hyperlink{AitSahalia}{A\"{i}t-Sahalia (1999)}) or simply "Mean-reverting stochastic volatility model" (in the context of Heston model, e.g. by \hyperlink{Andersen2007}{Andersen and Piterbarg (2007)}). Indeed, if one sets $a=0$, the CKLS model degenerates into a CEV model. For more information about the CEV model, one is recommended to refer to $\textbf{Lemma 3.5}$ in this paper. On the other hand, the terminology "Mean-reversion" refers to the observed phenomenon that the price of an asset, no matter how volatile it can be, will eventually move back towards its average over time, and significant deviations in price are typically unsustainable for long periods. The CKLS model embodies this principle by suggesting that short-term interest rates will revert to their long-term average. Specifically, if $a>0$, and the current short-term interest rate $\lambda_{t}$ exceeds its long-term average $a/b$, then the expected change in the interest rate will be negative-valued, and vice versa. Essentially, the mean $a/b$ serves as a balancing point for the process, earning it the descriptive name $\verb+"+$mean-reversion$\verb+"+$.\\
\noindent (3) Two seminal contributions merit special attention. First, the paper by \hyperlink{Andersen2007}{Andersen and Piterbarg (2007)} represents a watershed in CKLS-related research. The authors demonstrate that the condition $k>\frac{1}{2}$ is sufficient for the pathwise uniqueness and almost-sure strict positivity of solutions to the CKLS model. They also derive the model's stationary density by leveraging its ergodic properties. Prior empirical studies — such as \hyperlink{Chan1992}{Chan et al. (1992)} and a series of investigations in the late 1990s — primarily treated the CKLS model as an econometric tool, with limited focus on its analytical structure. In contrast, the work of \hyperlink{Andersen2007}{Andersen and Piterbarg (2007)} is widely recognized as the first to systematically generalize the CIR model within the broader CKLS framework. Second, \hyperlink{Mao2006}{Mao et al. (2006)}, along with subsequent developments by \hyperlink{Mao2013}{Mao and Szpruch (2013)}, \hyperlink{Wu2008}{Wu et al. (2008)}, and \hyperlink{Yang2020}{Yang et al. (2020)}, investigate the applicability of the Euler-Maruyama method to CKLS processes in more general settings. These studies emphasize that the weak and strong convergence — as well as the almost-sure stability — of Euler approximations to equation (1.1) are nontrivial properties that require careful and rigorous justification, thereby providing a rigorous theoretical basis for validating parametric estimation techniques.$\hfill\blacksquare$\\

\noindent \textbf{Theorem 1.2 (Some key properties of the solution to the CKLS model)}\\ 
The solution to (1.1), denoted by $\lambda_{t}$ for ${t\in[0,T]}$, has the following properties:\\
(1) $+\infty$ is an unattainable boundary for $\lambda_{t}$, $\forall k>0$.\\
(2) For $k>\frac{1}{2}$, $\lambda_{t}$ is a pathwise unique strong solution, being strictly positive-valued almost surely.\\
(3) For $0<k<\frac{1}{2}$, $0$ is always an attainable boundary for the solution $\lambda_{t}$.\\ 
(4) For $k=\frac{1}{2}$, $\lambda_{t}$ is a pathwise unique and strong solution. For case $2a\geq\sigma^{2}$, $\lambda_{t}$ is strictly positive-valued; for case $2a<\sigma^{2}$, $\lambda_{t}$ can reach $0$ with probability one (but will immediately bounce upward to a positive level after reaching $0$).\\ 
(5) For $0<k<\frac{1}{2}$, the following statement concerns the instantaneously reflecting extension of the diffusion at zero. As $T\to+\infty$ and $t\in[0,+\infty)$ (which also means the time filtration should be modified as $\{\mathcal{F}_{t}\}_{t\in[0,+\infty)}$), the CKLS process $\lambda_{t}$ is positive Harris recurrent (provided that $0<C_{k}<+\infty$, for which a simple sufficient condition is $b>0$) with a unique stationary density (invariant probability measure $\pi_{0}$, $\pi_{0}(dx)=p_{\infty}dx$) $p_{\infty}(x)=C_{k}x^{-2k}e^{\Lambda(x;k)}$, where
$$ \Lambda(x;k)=\left\{
\begin{aligned}
    &\frac{2}{\sigma^{2}}\Big(\frac{ax^{1-2k}}{1-2k}-\frac{bx^{2-2k}}{2-2k}\Big),&&k\in(0,\frac{1}{2}) \cup(\frac{1}{2},1) \cup (1,+\infty);\\
    &\frac{2}{\sigma^{2}}\big(a\text{log}x-bx\big),&&k=\frac{1}{2}\text{, (Cox-Ingersoll-Ross model)};\\
    &\frac{2}{\sigma^{2}}\big(-\frac{a}{x}-b\text{log}x\big),&&k=1\text{, (Brennan and Schwartz model)},
\end{aligned}
\right.
$$
with the constant $C_{k}=\Big(\int_{0}^{\infty}u^{-2k}e^{\Lambda(u;k)}du\Big)^{-1}$.\\
(6) For $k\in(\frac{1}{2},1)$, for any $p\geq 0 $, we have $\mathbb{E}[\text{sup}_{t\in[0,T]}(\lambda_{t})^{p}]<+\infty$ and $\mathbb{E}[\text{sup}_{t\in[0,T]}(\lambda_{t})^{-p}]<+\infty$. Therefore, it also holds that $\mathbb{E}[(\lambda_{t})^{p}]<+\infty$ and $\mathbb{E}[(\lambda_{t})^{-p}]<+\infty$.\\
(7) For $k\in(\frac{1}{2},1)$, if $0<C_{k}<\infty$, then for any $q$ with $\mathbb{R}_{+}\ni x\mapsto x^{q}\in \mathcal{L}^{1}(\pi_{0})$, the time average integral of order $q$ of $\lambda_{t}$ has the following $a.s.$ ergodic limit as $T\to+\infty$:
\begin{equation}\nonumber
     \frac{1}{T}\int_{0}^{T}\lambda_{t}^{q}dt\xrightarrow[T\to+\infty]{a.s.}\int_{0}^{\infty}x^{q}p_{\infty}(x)dx;
\end{equation}
Particularly, the integral $\int_{0}^{\infty}x^{q} p_{\infty}(x)dx<+\infty$ for $k\in(\frac{1}{2},1)$ and for any $q\in\mathbb{R}$ (when $q=0$, the above convergence holds trivially).\footnote{Note that when $k=\frac{1}{2}$, the CKLS process degenerates into the CIR process, and its stationary density is of a Gamma type with shape parameter $2a/\sigma^{2}$. Consequently, for $q\in\mathbb{R}$, $\int_{0}^{\infty}x^{q}p_{\infty}(x)dx<\infty
~\Leftrightarrow~q>-\frac{2a}{\sigma^{2}}$. Equivalently, for $p>0$, the stationary negative moment of order $p$ is finite if and only if $p<2a/\sigma^{2}$.}\\
\noindent \textit{Proof}: See the Appendix.\qed\\

\noindent \textbf{Remark 1.3}\\
(1) The drift $=a-bx$ is globally Lipschitz, hence globally H\"{o}lder (of order 1). The diffusion function $\sigma x^{k}$ is globally Lipschitz only for $k=0$ or $k=1$; is globally H\"{o}lder of order $k$ for every $0<k\leq 1$ (in particular H\"{o}lder-1 when $k=1$), and trivially H\"{o}lder-1 for $k=0$; is locally Lipschitz for $k>1$ on every compact subset of $[0,+\infty)$ and for $\frac{1}{2}\leq k<1$ on every compact subset of $(0,+\infty)$, while for $k<0$ and $0<k<\frac{1}{2}$, it is Lipschitz on every compact subset of $(0,+\infty)$.\\
(2) For $0<k<\frac{1}{2}$, to ensure that the process for $\lambda_{t}$ is unique, positively recurrent and achieves a stationary distribution, it is a standard approach to impose a boundary condition: a standard way is to assume that $\lambda_{t}$ is instantaneously reflected at the origin.\\
(3) For the case $k=\frac{1}{2}$ with $2a<\sigma^{2}$, the origin acts as a strong reflector, meaning that the duration $\lambda_{t}$ remaining at zero is $0$ in terms of the Lebesgue measure; therefore, there is no need for a specific boundary condition at $\lambda_{t}=0$.$\hfill\blacksquare$\\

\noindent However, similar to many other interest rate models, the CKLS model generally lacks a closed-form analytical solution. This is partly due to its role as a generalized framework—deliberately designed by econometricians to unify various model variants, as noted in \textbf{Remark 1.1 (3)}. One of our main contributions is to identify conditions under which a twice-differentiable transformation of the CKLS process has CIR dynamics under an equivalent probability measure, allowing its distributional properties to be studied through the CIR representation (see Sections 2 and 3).

\subsection{Some literature Review}
\subsubsection{On applications of model (1.1) in financial engineering via numerical solutions}
\noindent (\romannumeral1) Traditional methods devised for pricing financial derivatives include those based on Taylor expansions (see, e.g., \hyperlink{Stehlikova2013}{Stehl{\'\i}kov{\'a} (2013)}) and analytical approximation formulae (see e.g. \hyperlink{Choi2007}{Choi and Wirjanto (2007)} (In this paper, the authors consider a CKLS-type interest rate model under the physical measure, featuring a nonlinear drift term involving the market price of risk. Under the risk-neutral measure, the price of the zero-coupon bond satisfies a specific partial differential equation. The authors derive an analytical approximation formula for zero-coupon bond prices under the CKLS framework.), \hyperlink{Stehlikova2009}{Stehl{\'\i}kov{\'a} and \v{S}ev\v{c}ovi\v{c} (2009)}). A noteworthy contribution is the seminal work by \hyperlink{BaroneAdesi1999}{Barone-Adesi et al. (1999)}, which introduces a numerical approach known as the Box method. The simulated solution is subsequently employed to price zero-coupon bonds and bond options within the CKLS framework. The paper also presents a comparative analysis of bond and option prices obtained using both the Crank-Nicolson and Box methods. Building on this idea, several empirical studies — including \hyperlink{Byers1998}{Byers and Nowman (1998)}, Nowman and Sorwar (\hyperlink{Nowman1999a}{1999a}, \hyperlink{Nowman1999b}{1999b}), \hyperlink{Nowman2005}{Nowman and Sorwar (2005)}, and \hyperlink{Ma2008}{Ma et al. (2008)} — provide supporting evidence for the use of the CKLS model in modeling interest rates in financial markets. Furthermore, \hyperlink{Tangman2011}{Tangman et al. (2011)} propose an innovative computational technique for approximating the prices of zero-coupon bonds and bond options under the CKLS framework. This method employs a second-order finite difference approximation to discretize the pricing partial differential equations. In addition, it utilizes an exponential time integration scheme enhanced by optimal rational approximations derived via the Carath\'{e}odory-Fej\'{e}r method to solve the resulting semi-discrete system. In a similar vein, \hyperlink{Khor2012}{Khor et al. (2012)} and \hyperlink{Khor2014}{Khor and Pooi (2014)} adopt polynomial approximations of the first four moments of $\lambda_{t}$ to complete the discretization.\\

\noindent (\romannumeral2) As already noted in \textbf{Remark 1.1 (3)}, a rigorous proof that the Euler-Maruyama scheme converges to the exact solution was obtained only after the method had already been widely adopted in practice. For the case $k\in[\frac{1}{2},1)$ (when $k=\frac{1}{2}$, $2a\geq \sigma^{2}$), \hyperlink{Mao2006}{Mao et al. (2006)} investigate and verify the applicability of the Euler-Maruyama method for the CKLS process with a more general setting (In this paper, the authors use the term "the mean-reverting $k$-process" to denote the CKLS process. The hybrid variant introduces state-dependent parameters - specifically $a_{X_{t}}$, $b_{X_{t}}$, $\sigma_{X_{t}}$ - determined by a Markov chain $X_{t}$. In this way, the authors essentially generalize the CKLS model to a more flexible framework, also known as the regime-switching CKLS model. The analysis is then carried out by fixing $X_{t}=i$ to examine the model under a given regime.). To overcome difficulties arising from regime switching and non-Lipschitz coefficients, the authors develop several novel numerical techniques to establish the convergence of the Euler-Maruyama scheme. Building on this work, \hyperlink{Mao2013}{Mao and Szpruch (2013)} analyze the strong convergence and almost sure stability of Euler-Maruyama-type methods for SDEs with nonlinear, non-Lipschitz coefficients and further prove the global almost sure asymptotic stability of these schemes in such settings. For the case $k>1$, \hyperlink{Wu2008}{Wu et al. (2008)} characterize the analytical properties of the CKLS model and establish weak convergence in probability of the Euler-Maruyama approximation, drawing on results from \hyperlink{Mao2006}{Mao et al. (2006)}. Subsequently, \hyperlink{Yang2020}{Yang et al. (2020)} examine the moment convergence of the truncated Euler-Maruyama method at any fixed time $T$, and prove its strong convergence. More recently, for the CKLS and CEV processes, \hyperlink{Lileika2020}{Lileika and Mackevi{\v{c}}ius (2020)} propose a first-order split-step weak approximation method, which generates two-valued random variables at each discretization step and avoids regime switching near the origin. Last but not least, in the work by \hyperlink{Tsumurai2020}{Tsumurai (2020)}, the solution of the CKLS model (which is called the CEV-type process in the paper) is subjected to a non-linear transformation, leading to a new SDE for the transformed process. A numerical approximation of this SDE is then constructed by defining a piecewise continuous function based on a given threshold $\varepsilon$; this function can be shown to satisfy a global Lipschitz condition. As a result, the approximated SDE is globally Lipschitz and thus admits Malliavin differentiability of the CKLS process. The author proves the convergence of this approximation both in $\mathcal{L}^{2}$ and almost surely. Leveraging the Malliavin differentiability of the approximated solution, a Malliavin calculus-based analysis for the CKLS process is performed for a Heston-type model in which the CKLS process serves as the volatility component, and the arbitrage analysis as well as the computation of Greeks are conducted accordingly.

\subsubsection{On parameter estimation for model (1.1) and its econometric applications}
\noindent (\romannumeral1) The seminal work by \hyperlink{Chan1992}{Chan et al. (1992)} employs the Euler-Maruyama method in conjunction with the generalized method of moments (GMM), as developed by \hyperlink{Hansen1982}{Hansen (1982)}, to estimate parameters, conduct inference, and compute test statistics for model evaluation. Applying GMM to U.S. Treasury bill data, the authors obtain an estimate of $k=1.449$. In contrast to traditional approaches such as Bayesian and quasi-maximum likelihood (QML) methods — which impose strict distributional assumptions on the transition density of the process — GMM relies primarily on the asymptotic properties of sample means, as guaranteed by the central limit theorem. This makes GMM a preferred method for estimating the CKLS model due to its flexibility with respect to distributional assumptions. However, several alternative studies — such as \hyperlink{Brenner1996}{Brenner (1996)}, \hyperlink{Nowman1998}{Nowman (1998)}, \hyperlink{Beuermann2005}{Beuermann et al. (2005)}, among others — advocate for QML estimation for the following reasons: \ding{172} \hyperlink{Broze1995}{Broze et al. (1995)} report that GMM performs poorly when $k>1$, reflecting heightened sensitivity of volatility to the current interest rate level. This observation is also supported by empirical findings in \hyperlink{Byers1998}{Byers and Nowman (1998)}, which document instances of $k>1$, thereby highlighting the potential shortcomings of GMM in such regimes and motivating the use of QML. \ding{173} \hyperlink{Dahlquist1996}{Dahlquist (1996)} argues that GMM estimators yield less powerful statistical tests compared to their QML counterparts. \ding{174} \hyperlink{Broze1995}{Broze et al. (1995)} further emphasize that QML estimation is generally more efficient than GMM. \ding{175} A notable advantage of QML over GMM, as discussed by \hyperlink{Nowman1997}{Nowman (1997)}, lies in its ability to incorporate more precise estimators, thereby enhancing the overall accuracy of model estimation. In most studies utilizing QML, the discretization scheme of \hyperlink{Bergstrom1984}{Bergstrom (1984)} is adopted for Gaussian cases. Of particular note is the study by \hyperlink{Nowman1997}{Nowman (1997)}, which uses the same treasury bill dataset as \hyperlink{Chan1992}{Chan et al. (1992)} and yields an estimate of $k=1.361$, closely aligning with the earlier result. For non-Gaussian settings — especially where data exhibit leptokurtosis — the scheme proposed by \hyperlink{Newey1997}{Newey and Steigerwald (1997)} is employed, often in combination with Student's t-distributed innovations. Lastly, for a Bayesian approach to inference within the CKLS framework that incorporates MCMC techniques, Li et al. (2010, unpublished working paper; available at \url{https://www.academia.edu/207922/Bayesian_Analysis_of_CKLS_models_for_US_{s}hort_term_Interest_Rate}; accessed 17 July 2025) present a representative study in which the model integrates an ARMA-GARCH error structure based on the asymmetric exponential power distribution.\\

\noindent (\romannumeral2) Considering parametric estimation methods besides using the Euler-Maruyama scheme, two classic methods are worthy of highlighting. \hyperlink{AitSahalia1999}{A\"{i}t-Sahalia (1999)} (see also \hyperlink{AitSahalia2002}{A\"{i}t-Sahalia (2002)}) proposes a likelihood-based estimation method for diffusion models observed at discrete intervals, using a Hermite polynomial expansion of the transition density. A key step is the Lamperti transform, which standardizes the diffusion coefficient by converting the original process $X_{t}$ into a new process $U_{t}$ with unit diffusion coefficient. This facilitates the Hermite expansion of the transition density $p_{U}$ around a Gaussian density. The approximate transition density of the original process $p_{X}$ is then recovered through the inverse transform and the Jacobian formula. The resulting approximated log-likelihood function is maximized to obtain an estimator $\hat\theta_{n}^{(J)}$, which the author proves to be asymptotically normal under suitable conditions. Though the method requires the analytical invertibility of the Lamperti transform and Hermite coefficients, it offers high estimation accuracy for stationary diffusion processes. \hyperlink{Shoji1998}{Shoji and Ozaki (1998)} (see also \hyperlink{Ozaki1992}{Ozaki (1992)}) propose a method called local linearization (LL), which replaces the drift by a linear function on each sampling interval while treating the diffusion coefficient as constant. (If the diffusion coefficient is not constant, we may first apply a Lamperti transform to obtain a model with unit diffusion coefficient and then use the same LL machinery.) This step-wise linear SDE has a Gaussian transition law whose mean and variance can be written in closed form, allowing the log-likelihood to be evaluated and maximized directly. By capturing local curvature in the drift, LL is markedly more accurate than the Euler scheme, most notably for nonlinear dynamics, yet remains computationally light.  Overall, LL provides an efficient and precise route to maximum-likelihood estimation for discretely observed diffusion processes.\\

\noindent (\romannumeral3) In several recent studies, innovative approaches have been developed for constructing parametric estimators, including maximum and quasi-maximum likelihood estimators for the drift parameters in both continuous- and discrete-time settings, with their performance validated using high-frequency data over infinite time horizons. The asymptotic normality of such estimators has also been a topic of active investigation. Additionally, these studies also address the estimation of the diffusion coefficient (see, e.g., \hyperlink{Mazzonetto2024}{Mazzonetto and Nieto (2024)}, \hyperlink{Lyu2025}{Lyu and Nkurunziza (2025)}, \hyperlink{Wei2020}{Wei (2020)}). The computation of realized volatility naturally leads to nonparametric estimators of the volatility parameter $\sigma$ and the elasticity coefficient $k$. Building on this idea, a novel nonparametric method for jointly estimating $k$ and $\sigma$—involving complex number techniques—has been proposed by \hyperlink{Dokuchaev2017}{Dokuchaev (2017)}.\\

\noindent (\romannumeral4) Recently, for the case $\frac{1}{2}<k<1$, \hyperlink{Mishura2022}{Mishura et al. (2022)} treat the Dol\'{e}ans-Dade exponential (Radon-Nikod\'{y}m derivative) as a likelihood function and derived the expression for the MLE of the unknown drift parameters through continuous observations of a sample path. The strong consistency and asymptotic normality of this maximum likelihood estimator are also derived. This approach draws inspiration from the method developed by the series papers by \hyperlink{Alaya2012}{Alaya and Kebaier (2012)} and \hyperlink{Alaya2013}{Alaya and Kebaier (2013)}, for the case when the diffusion parameter is assumed known, where the drift parameters of the CIR model are estimated, and the strong consistency and asymptotic normality of the maximum likelihood estimator are proven based on Laplace transform techniques, in both ergodic and non-ergodic settings. \hyperlink{Mishura2022}{Mishura et al. (2022)} also introduce a new strongly consistent estimator for drift parameters by extending the estimation techniques previously suggested by \hyperlink{Dehtiar2022}{Dehtiar et al. (2022)} for the CIR model. See also \hyperlink{DeRossi2010}{De Rossi (2010)} and \hyperlink{Overbeck1997}{Overbeck and Ryd\'{e}n (1997)}, \hyperlink{Overbeck1998}{Overbeck (1998)} for more about parametric estimation of the CIR model.

\subsubsection{On fixed-\texorpdfstring{$k$}{k} specializations and extensions of model (1.1)}
\noindent (\romannumeral1) We will no longer dwell on the best-known and most extensively studied CIR model—many of its properties (when $k$ is assumed to be $\frac{1}{2}$), including its connections to the CEV model and to the Bessel process, are discussed at several other points in this paper, with full references provided. Instead, we now turn to another frequently overlooked variant within the CKLS family: the $3/2$-model. The choice of the $3/2$-model is in fact supported by empirical evidence provided by studies such as \hyperlink{Chan1992}{Chan et al. (1992)} and \hyperlink{Nowman1997}{Nowman (1997)}. In an early phase, the monograph \hyperlink{Lewis2000}{Lewis (2000)} offers a thorough survey of option-pricing techniques under a range of stochastic-volatility specifications. Besides discussing the GARCH diffusion and risk-adjusted processes, the paper treats several CKLS special cases that arise for particular values of~$k$—notably the CIR model, the so-called "$3/2$" model, and the OU process. \hyperlink{Ahn1999}{Ahn and Gao (1999)} derive a closed-form bond pricing formula under the $3/2$-model for interest rates using the Girsanov theorem, where the drift term takes a distinctive quadratic form in the interest rate. The quadratic drift structure is adopted so as to make the process exhibit a substantial nonlinear mean-reverting behavior when the interest rate exceeds its long-run mean. In addition, the authors document that this unique type of SDE admits a concave relationship between interest rates and yields. \hyperlink{Carr2007}{Carr and Sun (2007)} propose and analyze the $3/2$ model with a quadratic drift term to describe the normal volatility of instantaneous variance, showing that it is theoretically sound and empirically well supported. Moreover, the authors further demonstrate its analytical tractability by deriving a closed-form expression for the joint Fourier-Laplace transform, highlighting its applicability in pricing volatility derivatives.\\

\noindent (\romannumeral2) Over the past decade, numerous studies have explored generalizations of the CKLS model. The first type of variant involves modifying the structure or parameters of the CKLS model. Most recently, \hyperlink{Mazzonetto2024}{Mazzonetto and Nieto (2024)} introduce a variant of the CKLS model, which is a continuous-time, self-exciting and ergodic process, called the threshold CKLS process, which incorporates the presence of multiple thresholds governing shifts in dynamics. \hyperlink{Lyu2025}{Lyu and Nkurunziza (2025)} extend the CKLS model by letting the mean-reversion level be a deterministic periodic function, improving its fit to realistic rate dynamics. Using transition semigroup theory, they prove the ergodicity and positive Harris recurrence of the discrete chain. They derive unrestricted and restricted MLEs with joint asymptotic normality under local alternatives, propose a class of shrinkage estimators, and show via simulation that these (and in particular the positive-part SE) outperform the standard UMLE. \hyperlink{Cai2015}{Cai and Wang (2015)} investigate the asymptotic behavior of the CKLS model with small random perturbation $\sqrt{\epsilon}$ and obtain the central limit theorem and the moderate deviation principle for the solution of this model when $\epsilon\to 0$. \hyperlink{Baldi2011}{Baldi and Caramellino (2011)} establish Freidlin-Wentzell large deviation estimates for the same model under minimal assumptions for diffusion processes on the positive half-line, applicable to the CKLS model with non-Lipschitz but H\"{o}lder continuous coefficients. \hyperlink{Feng2012}{Feng et al. (2012)} drive the stock‐price SDE with a fast mean-reverting CKLS-type volatility process and use viscosity-solution techniques to prove a rare-event large-deviation principle governing short-time, out-of-the-money option prices. Building on this LDP, they derive asymptotic formulas for both option prices and their implied volatility in two multiscale regimes ($\delta=\varepsilon^{2}$ and $\delta=\varepsilon^{4}$). The paper's main financial contribution is to furnish rigorous, model-agnostic short-maturity approximations for option valuations and implied volatility under nonlinear volatility dynamics beyond the classical Heston case.\\

\noindent (\romannumeral3) The second type of variant replaces the Brownian motion in the CKLS model with other stochastic processes. \hyperlink{Wei2020}{Wei (2020)} proposes a least squares estimator for the CKLS model driven by small L\'{e}vy noises using discrete observations. The estimator is constructed from a contrast function that captures the weighted squared deviation between the observed increments and their Euler-Maruyama approximation. The paper derives the explicit form of the estimator, analyzes the estimation errors, and proves the consistency as the diffusion coefficient $\sigma$ approaches $0$ and the sample size approaches $+\infty$. This expression of the estimator closely resembles the approach in the work by \hyperlink{Mishura2022}{Mishura et al. (2022)} for the CKLS model with Wiener noise, reflecting the well-known asymptotic equivalence between least squares estimation and maximum likelihood estimation (see, e.g., \hyperlink{Skouras2000}{Skouras (2000)}, \hyperlink{Mendy2013}{Mendy (2013)}) in the context of SDE parameter estimation. For the CKLS model driven by fractional Brownian motion, researchers have examined the process from various perspectives. \hyperlink{Kubilius2021}{Kubilius and Med\v{z}i{\=u}nas (2021)} study the CKLS model driven by fractional Brownian motion with non-Lipschitz diffusion functions and without linear growth conditions. By applying the Lamperti transform, the authors derive conditions that ensure the positivity of the solutions and show that for the fractional CKLS model with $k>1$, the trajectories are not necessarily positive-valued. The authors further establish the almost sure convergence rate of the backward Euler approximation scheme and provide a strongly consistent and asymptotically normal estimator of the Hurst index $H>\frac{1}{2}$ for positive-valued solutions (see also \hyperlink{Gyongy1996}{Gy{\"{o}}ngy and Krylov (1996)}). Schl{\"u}chtermann and Yang (2016, unpublished working paper; available at \url{https://www.researchgate.net/publication/299670926_{n}ote_on_fractional_CLKS-type_{s}tochastic_differential_equation_-path-wise_and_in_the_Wick_{s}ense}; accessed 17 July 2025) show that a generalized fractional CKLS model with time-varying drift (including the CIR model with positive-valued thresholds, see \hyperlink{Zahle1998}{Z\"{a}hle (1998)}), has a positive-valued solution, both for the pathwise integral and in the Wick sense (see \hyperlink{Holden1996}{Holden et al. (1996)} or \hyperlink{Hu2003}{Hu and {\O}ksendal (2003)} for this concept). Moreover, \hyperlink{Zhao2022}{Zhao and Xu (2022)} address the inverse problem of estimating the time-varying diffusion $\sigma_{t}$ and the elasticity parameter $k$ in the fractional CKLS model for European options from a limited number of market observations. Tikhonov regularization and the ADMM algorithm are applied to ensure the stability of the solution and efficient optimization. In the framework of rough path analysis, \hyperlink{Marie2014}{Marie (2014)} considers the CKLS‐type mean‐reverting SDE driven by a general centered Gaussian rough path, thus treating the classical CKLS model as a "rough" variant. By leveraging rough-path techniques, the author proves global existence and uniqueness of the solution, establishes continuity and differentiability of the associated It\^{o} map, derives $\mathcal{L}^{p}$-convergent Euler approximations with explicit rates, and obtains a large-deviation principle and density for the underlying process. Finally, they showcase the model's applicability by formulating and analyzing a pharmacokinetic mean‐reversion model, illustrating how this variant of the CKLS model can capture dynamics beyond finance.

\subsection{Revisiting Hu et al. (2015): Errors in Model Formulation, Derivations, Propositions and Proofs}
\noindent After completing our paper, we became aware of the paper by \hyperlink{Hu2015}{Hu et al. (2015)}, which had already investigated the same problem using a similar approach and arrived at comparable conclusions. However, there appear to be several significant problems with the results presented in the paper.\\

\noindent First of all, upon thorough examination, we identified a major issue in the initial derivation of their paper, which appears in Equation (3) on page 70 of \hyperlink{Hu2015}{Hu et al. (2015)}. Concisely speaking, this fundamental error originates from the incorrect calculation of the square root of $(1-\gamma)^{2}$ as $1-\gamma$ rather than $\vert 1-\gamma \vert$, leading to a sign error in the stated transformed diffusion coefficient when $C>0$ and $\gamma>1$ (to be specific, the root of $\frac{C^{2}}{4(1-\gamma)^{2}}$ is $\lvert \frac{C}{2(1-\gamma)}\lvert$, rather than $\frac{C}{2(1-\gamma)}$). As a result, for $C>0$ and $\gamma>1$, the transformed diffusion coefficient has the opposite sign to the prescribed positive square-root coefficient when the driving Wiener process is kept unchanged.\\

\noindent Hu et al. (2015) explicitly assume $C>0$ at the beginning of Section 2. Under this assumption, the distinction between the signed expression $Cr_{t}^{1-\gamma}/[2(1-\gamma)]$ and the positive square root $\sqrt{f(r_{t})}$ is essential. To be specific, for the deduction $\sigma C \frac{C}{2(1-\gamma)}(r_{t})^{1-\gamma}=\sigma C\sqrt{f(r_{t})}=\sigma C\sqrt{Y_{t}}$ to hold, one must ensure that the diffusion coefficient satisfies $\sigma C>0$ and that the term under the square root, $f(r_{t})$, always remains positive-valued. This requires both $C>0$ and $\frac{(r_{t})^{1-\gamma}}{1-\gamma}>0$. When $\gamma>1$, by the strict positivity of the CKLS process (recall \textbf{Theorem 1.2}), $r_{t}$ remains strictly positive-valued, making $(r_{t})^{1-\gamma}$ remain positive-valued for sure (because $r_{t}>0$ together with $1-\gamma<0$ makes $(r_{t})^{1-\gamma}>0$). Note the fact that $A^{-B}=\frac{1}{A^B}>0$ for $A,B>0$), the signed expression $\frac{C}{2(1-\gamma)}(r_{t})^{1-\gamma}$ becomes negative-valued and therefore cannot equal the positive square root $\sqrt{f(r_{t})}$. Thus, with $C>0$ and the driving Wiener process kept unchanged, the case $\gamma>1$ does not yield the prescribed positive square-root diffusion term. Accordingly, in the present construction we restrict attention to $C>0$ and $\frac{1}{2}\leq\gamma<1$, with $2a\geq\sigma^{2}$ when $\gamma=\frac{1}{2}$.\\

\noindent Secondly, in equation (7) on page 72 of \hyperlink{Hu2015}{Hu et al. (2015)}, the authors mistakenly wrote $\sigma$ instead of $\sigma^{2}$ in expressions where the latter should appear. This error propagates into equations (8) and (9), producing the factor $y^{-\gamma/\sigma}$ instead of the correct factor $y^{-\gamma}$ in the scale density. Consequently, the restriction $\gamma/\sigma\geq1$ used in their boundary classification is an artifact of the incorrect coefficient and does not provide the claimed non-exit condition for the auxiliary diffusion. If the authors had derived the expression of $p(x)$ correctly, they would have gotten the following expression: $\lim_{x\downarrow0}p(x)=\lim_{x\downarrow0}\text{exp}\{-\frac{b}{\sigma^{2}(1-\gamma)}\}\int_{1}^{x}y^{-\gamma}\text{exp}\{\frac{b}{\sigma^{2}(1-\gamma)}y^{2(1-\gamma)}\}dy$ when $\frac{1}{2}\leq \gamma<1$, $b>0$, $\frac{\gamma}{\sigma}\geq 1$. Since $\frac{1}{2}\leq\gamma<1$, $0<1-\gamma\leq\frac{1}{2}$ and $0<2(1-\gamma)\leq 1$, if one assumes that $b>0$ (thus $K\overset{\text{def}}{=}\frac{b}{\sigma^{2}(1-\gamma)}>0$, $e^{K}>0$, $e^{-K}>0$), it would be $\lim_{x\downarrow0}p(x)=\lim_{x\downarrow0}e^{-K}\int_{1}^{x}y^{-\gamma}\text{exp}\{Ky^{2(1-\gamma)}\}dy$. As $y\downarrow0$, one may notice that $Ky^{2-2\gamma}\to0$ because $0<2(1-\gamma)\leq 1$, so $\text{exp}\{Ky^{2(1-\gamma)}\}\to1$, and thus consequently $y^{-\gamma}\text{exp}\{Ky^{2(1-\gamma)}\}\sim y^{-\gamma}$. Since $\gamma<1$, this asymptotic equivalence implies integrability near zero. As a result $\lim_{x\downarrow0}p(x)=-e^{-K}\int_{0}^{1}y^{-\gamma}\text{exp}\{Ky^{2(1-\gamma)}\}dy\in(-\infty,0)$, i.e. $\lim_{x\downarrow0}p(x)$ is some negative value, not $-\infty$. The implications of this corrected boundary calculation for the claimed martingale property are examined below.\\

\noindent Lastly and most importantly, the choice of the constant drift term $\sigma^{2}C^{2}/4$ in the target square-root diffusion in \hyperlink{Hu2015}{Hu et al. (2015)} leads to a contradiction between its boundary behavior and the claimed equivalent change of measure. To see this, consider the parameter regime $\frac{1}{2}<\gamma<1$, with $a,b,\sigma,C>0$ and $r_{0}>0$. The constant drift and diffusion coefficients of their target equation give the associated squared Bessel dimension
\begin{equation}\nonumber
	d=\frac{4(\sigma^{2}C^{2}/4)}{(\sigma C)^{2}}=1<2.
\end{equation}
Thus, zero is accessible from every positive initial value: the target diffusion has a strictly positive probability of hitting zero before any prescribed time $T>0$. This boundary property holds irrespective of the sign of its linear drift coefficient. Moreover, changing the positive constant $C$ cannot remove this accessibility, since $C$ cancels from the expression for $d$.\\

\noindent However, the original CKLS process remains strictly positive under $\mathbb{P}$ in this parameter regime. Consequently, the transformed process
\begin{equation}\nonumber
	Y_{t}=f(r_{t}), \qquad f(x)=\frac{C^{2}}{4(1-\gamma)^{2}}x^{2(1-\gamma)},
\end{equation}
\noindent also remains strictly positive under $\mathbb{P}$. For any fixed $T>0$, define $A_{T}=\big\{\inf_{0\leq t\leq T}Y_{t}=0\big\}$. Continuity and strict positivity imply $\mathbb{P}(A_{T})=0$. If the proposed stochastic exponential $R$ were a true martingale, it would define a probability measure $\mathbb{Q}$ on $\mathcal{F}_{T}$ by
\begin{equation}\nonumber
	\left.\frac{d\mathbb{Q}}{d\mathbb{P}}\right|_{\mathcal{F}_{T}}=R_{T},\qquad\text{and hence}\qquad \mathbb{Q}(A_{T})=\mathbb{E}_{\mathbb{P}}\big[R_{T}\mathbbm{1}_{A_{T}}\big]=0.
\end{equation}
In fact, this measure would be equivalent to $\mathbb{P}$, since $R_{T}>0$ almost surely. The resulting impossibility of hitting zero contradicts the strictly positive hitting probability required by the target square-root equation. Thus, no such absolutely continuous change of probability measure can produce the asserted dynamics for $Y=f(r)$ on the whole interval $[0,T]$.\\

\noindent This boundary incompatibility is also revealed by correcting the authors' application of Feller's test for explosions. For $1/2<\gamma<1$, their Girsanov kernel is
\begin{equation}\nonumber
	q(x)=\frac{\gamma\sigma}{2}x^{\gamma-1}-\frac{a}{\sigma}x^{-\gamma}.
\end{equation}
Substitution into the auxiliary equation gives
\begin{equation}\nonumber
	d\widetilde{r}_{t}=\Big(\frac{\gamma\sigma^{2}}{2}\widetilde{r}_{t}^{2\gamma-1}-b\widetilde{r}_{t}\Big)dt+\sigma\widetilde{r}_{t}^{\gamma}dB'_{t}, \qquad \widetilde{r}_{0}=r_{0}.
\end{equation}
Using the notation introduced above, the scale density and scale function are
\begin{equation}\nonumber
	\rho(y)=\text{exp}\Big\{-2\int_{1}^{y}\frac{\frac{\gamma\sigma^{2}}{2}z^{2\gamma-1}-bz}{\sigma^{2}z^{2\gamma}}dz\Big\}=e^{-K}y^{-\gamma}e^{Ky^{2(1-\gamma)}},\qquad p(x)=\int_{1}^{x}\rho(y)dy.
\end{equation}
As established above, $p(0+)\in(-\infty,0)$.\\

\noindent To determine accessibility of zero, the remaining integrability condition in Feller's test must also be checked. As $y\downarrow0$,
\begin{equation}\nonumber
	\rho(y)\sim e^{-K}y^{-\gamma}, \qquad p(y)-p(0+)\sim\frac{e^{-K}}{1-\gamma}y^{1-\gamma}.\qquad \text{Consequently,}\quad
	\frac{p(y)-p(0+)}{\sigma^{2}y^{2\gamma}\rho(y)}\sim\frac{1}{\sigma^{2}(1-\gamma)}y^{1-2\gamma}.
\end{equation}
Because $1-2\gamma>-1$, this expression is integrable near zero:
\begin{equation}\nonumber
	\int_{0}^{\varepsilon}\frac{p(y)-p(0+)}{\sigma^{2}y^{2\gamma}\rho(y)}dy<\infty
\end{equation}
for sufficiently small $\varepsilon>0$. Together with $p(0+)>-\infty$, this establishes that zero is accessible and that the auxiliary diffusion exits $(0,\infty)$ with positive probability.\\

\noindent The same boundary conclusion can also be obtained using the scale and test functions discussed in \hyperlink{Karatzas2012}{Karatzas and Shreve (2012)} [\textit{Brownian Motion and Stochastic Calculus}, p.~348, Problem~5.27]. Denoting these functions here by $\psi$ and $\phi$, respectively, and taking the reference point $c=1$, we define
\begin{equation}\nonumber
	\psi(x)=p(x), \qquad \phi(x)=\int_{1}^{x}\rho(y)\Big(\int_{1}^{y}\frac{2dz}{\sigma^{2}z^{2\gamma}\rho(z)}\Big)dy.
\end{equation}
For $0<x<1$, reversing the integration limits and applying Tonelli's theorem gives
\begin{equation}\nonumber
	\phi(x)=\int_{x}^{1}\rho(y)\Big(\int_{y}^{1}\frac{2dz}{\sigma^{2}z^{2\gamma}\rho(z)}\Big)dy=2\int_{x}^{1}\frac{\psi(z)-\psi(x)}{\sigma^{2}z^{2\gamma}\rho(z)}dz.
\end{equation}
Hence, by monotone convergence and the preceding integrability estimate,
\begin{equation}\nonumber
	\phi(0+)=2\int_{0}^{1}\frac{p(z)-p(0+)}{\sigma^{2}z^{2\gamma}\rho(z)}dz<\infty.
\end{equation}
On the other hand, $\psi(+\infty)=+\infty$ implies $\phi(+\infty)=+\infty$. Thus, Feller's test again shows that the auxiliary diffusion can exit $(0,\infty)$ through zero with positive probability. The finiteness of $\phi(0+)$ follows from the explicit calculation above; it cannot be inferred solely from the finiteness of $\psi(0+)$.\\

\noindent Since the original CKLS process does not exit this state space, Corollary 2.2 of \hyperlink{Mijatovic2012}{Mijatovi\'c and Urusov (2012)}, together with the finite-horizon remark following their Theorem 2.1, implies that $R$ is a strict local martingale. In particular,
\begin{equation}\nonumber
	\mathbb{E}_{\mathbb{P}}[R_{T}]<1, \qquad T>0.
\end{equation}
Therefore, correcting the coefficient and completing Feller's test yields a conclusion opposite to the claimed non-exit and true-martingale properties. The additional restrictions $\gamma/\sigma\geq1$ and $b>0$ do not resolve this contradiction. The proposed density fails to define the probability measure required for their Girsanov representation in this parameter regime, and the subsequent explicit-solution and distributional conclusions cannot be asserted for the original CKLS process under the claimed equivalent measure. Although localization away from zero permits stopped changes of measure, it does not justify the unstopped representation stated in their results.

\subsection{Structure of the paper}
\noindent In \textbf{Section 2}, the main result of this paper is narrated as follows: We introduce a certain twice-differentiable mapping that maps the general CKLS model to an intermediate/transitional SDE of a specific expression, whose drift term is yet cumbersome and intractable, where its parameters are constructed by the original parameters of the CKLS model. In this procedure, particular attention should be paid to the domains that the parameters could take values from, i.e. the parameter space is strictly restrained. Next, we apply the Cameron-Martin-Girsanov-Maruyama Measure Transform theorem to the "immature" process and obtain a process of a CIR type with concise and tractable parameters. Yet, the available moment estimates alone do not establish Novikov's or Kazamaki's condition, so it remains to prove that the induced Dol\'{e}ans–Dade exponential is a true martingale and hence defines the Radon-Nikod\'{y}m density for the measure change. \textbf{Section 3} gives several subsidiary results obtained for our model based on established theories about the CIR process. These subsidiary results may potentially be used in real-world financial studies. Finally, in \textbf{Section 4}, we explain why the available moment estimates alone do not establish Novikov's or Kazamaki's condition, and then outline the proof of the martingale property using an auxiliary diffusion. After this, we detail how the key proof of the claim that the induced Dol\'{e}ans-Dade exponential is a true martingale (Radon-Nikod\'{y}m derivative), which is left unproven in Section 2, is achieved. A foundational result is given by \hyperlink{Mijatovic2012}{Mijatovi\'c and Urusov (2012)}  who establish necessary and sufficient conditions under which a generalized Girsanov transformation yields a Radon-Nikod\'{y}m derivative that is a true martingale. Their characterization is formulated in terms of Feller's boundary classification and the associated explosion test for one-dimensional diffusion processes, as developed by \hyperlink{Feller1952}{Feller (1952)}.

\section{Main result: Applying the Cameron-Martin-Girsanov-Maruyama Measure Transform on \texorpdfstring{$\lambda_{t}$}{lambda	extunderscore t} that has been transformed by a twice-differentiable Function Parameterized by \texorpdfstring{$k$}{k} to Derive a Cox-Ingersoll-Ross-Type Model}
\noindent Consider a twice-differentiable function $T: \mathbb{R}_{+}\xrightarrow{}\mathbb{R}_{+}$ and a constant $L\in\mathbb{R}$ such that:
\begin{equation}\nonumber
    \frac{dT(x)}{dx}\cdot x^{k}=L\cdot\big[T(x)\big]^\frac{1}{2}.
\end{equation}
Solving this ordinary differential equation with the help of the technique of separation of variables gives:
\begin{equation}
\big[T(x)\big]^{\frac{1}{2}}=\left\{
\begin{aligned}
    &\frac{L}{2(1-k)}x^{1-k}+\text{constant, when }k\neq 1;\\
    &\frac{L}{2}\text{log}x+\text{constant, when }k=1.
\end{aligned}
\right.
\end{equation}
We choose a zero value to the integral constant, obtaining:
\begin{equation}
T(x)=\left\{
\begin{aligned}
    &\frac{L^{2}}{4(1-k)^{2}}x^{2(1-k)}\text{, when }k\neq 1;\\
    &\frac{L^{2}}{4}\big(\text{log}x\big)^{2}\text{, when }k=1.
\end{aligned}
\right.
\end{equation}
Note that $[T(x)]^{\frac{1}{2}}$ is strictly positive on the prescribed state space, which is an important requirement in the following discussion. For case $k=1$, when $L>0$, $x$ should be taken from $(1,+\infty)$; when $L<0$, $x$ should be taken from $(0,1)$; for case $k\neq 1$, $T(x)$ always takes positive values, as does $[T(x)]^{\frac{1}{2}}$, and so does the product of $x^{1-k}$ and $\frac{L}{1-k}$. Since $L$ and $1-k$ are deterministic after a certain model together with its parameters assigned certain values, $\frac{L}{1-k}$ is also determinant. \textbf{This requires one to assume that $x^{1-k}$ is either strictly positive-valued or strictly negative-valued for all $x$ in the prescribed domain.} When $x^{1-k}>0$, the assumption $\frac{L}{1-k}>0$ is needed; when $x^{1-k}<0$, the assumption $\frac{L}{1-k}<0$ is needed.\\

\noindent In other words, we need to first specify appropriate values for $k$ and $L$, where the value of $k$ determines the range of $x^{1-k}$. If $x^{1-k}$ is not strictly positive-valued or strictly negative-valued over the domain, then the chosen values of $k$ and $L$ are inappropriate.\\

\noindent Due to the original setting $\lambda_{t}\lvert_{t=0}=\lambda_{0}>0$, it is easy to see that $(\lambda_{t})^{1-k}$ cannot always be strictly negative no matter what value $k$ is given. As a result, $x^{1-k}$ must be assumed to be always strictly positive-valued, and $\frac{L}{1-k}$ must be strictly positive-valued. Therefore, either the case $L>0$ and $1-k>0$ or the case $L<0$ and $1-k<0$ must be assumed.\\

\noindent Based on this, we now consider the cases $k=1$ and $0<k<\frac{1}{2}$ separately. To be more specific, when $k=1$, the state space of the CKLS process is $(0,+\infty)$ (recall \textbf{Theorem 1.2 (2)}), over which $\text{log}x$ takes both positive and negative values. Therefore, for any fixed nonzero $L$, the expression for $[T(x)]^{\frac{1}{2}}$ obtained above cannot remain strictly positive-valued throughout this state space. As a result, the case $k=1$ is excluded from the present construction, which requires a globally defined, strictly positive-valued mapping on $(0,+\infty)$.\\

\noindent When $0<k<\frac{1}{2}$, zero is an attainable boundary of the CKLS diffusion (recall \textbf{Theorem 1.2 (3)}). Although $T(x)$ extends continuously to $[0,+\infty)$ with $T(0)=0$, its second derivative is singular at zero. Moreover, the Girsanov kernel to be introduced later contains negative powers of the state variable and is also singular at zero. Thus, the It\^{o} and Girsanov calculations performed in the interior of the state space do not, by themselves, justify extending the construction across this boundary.\\

\noindent There is a further restriction arising from the equivalent change of measure that we seek to construct. As will be seen below, choosing the constant drift term of the target CIR process as $\sigma^{2}L^{2}/2$ makes Feller's condition hold at equality, so that this process cannot reach zero from a positive initial value. Since $T(0)=0$ and $T(x)>0$ for $x>0$, the processes $\lambda_{t}$ and $X_{t}=T(\lambda_{t})$ reach zero at the same time. Consequently, on any finite horizon over which the original CKLS process has a positive probability of hitting zero under $\mathbb{P}$, the desired CIR representation would assign probability zero to the same event under $\mathbb{Q}$. This contradicts the equivalence of the two measures on that horizon.\\

\noindent As a result, the case $0<k<\frac{1}{2}$ must also be excluded from the present equivalent-measure construction. Transformations stopped before the boundary, or constructions requiring only absolute continuity, would need separate analysis. Having acknowledged these restrictions, we will not discuss the cases $0<k<\frac{1}{2}$ and $k=1$ further within the present framework.\\

\noindent We now turn to the remaining cases $k\geq\frac{1}{2}$ and $k\neq1$, with the additional assumption $2a\geq\sigma^{2}$ when $k=\frac{1}{2}$. Under these conditions, $\lambda_{t}$ remains strictly positive almost surely (recall \textbf{Theorem 1.2 (2) and (4)}). Therefore, $(\lambda_{t})^{1-k}$ is strictly positive-valued, regardless of whether $1-k$ is positive or negative. To be more specific, when $0<1-k\leq\frac{1}{2}$, the function $h(x)=x^{1-k}$ is increasing on $(0,+\infty)$; when $1-k<0$, it is decreasing on $(0,+\infty)$. In either case, $h(x)>0$ for all $x>0$. As a result, the requirement that $[T(x)]^{\frac{1}{2}}$ be strictly positive throughout the prescribed domain implies that $\frac{L}{1-k}>0$. Hence, either the case $L>0$ and $1-k>0$, or the case $L<0$ and $1-k<0$, must be assumed.\\

\noindent A summary of this basic setting is to be referred to in (2.4) where soon we may see that $L>0$ is also an indispensable requirement. That is, the case $L<0$ and $1-k<0$ will invalidate a key property/effect of the mapping $T(x)$.\\

\noindent Having obtained (2.2), we have:
\begin{align}\nonumber
   \frac{dT(x)}{dx}&=\frac{L^{2}}{2(1-k)}x^{1-2k},\\
   \frac{d^{2}T(x)}{dx^{2}}&=\frac{L^{2}(1-2k)}{2(1-k)}x^{-2k},\nonumber\\
   T^{-1}(x)&=\Big[\frac{2(1-k)}{L}\Big]^{\frac{1}{1-k}}x^{\frac{1}{2(1-k)}}\nonumber.
\end{align}

\noindent \textbf{Remark 2.1}\\
(1) Under the restriction (2.4), the inverse of $T$, whose existence is guaranteed by the inverse function theorem, increases strictly over $(0,+\infty)$.\\
(2) The idea of introducing this transform is that $\frac{dT(x)}{dx}$ times $x^{k}$ will have the exact same expression as $L$ times $\big[T(x)\big]^{\frac{1}{2}}$, that is, $\frac{dT(x)}{dx}x^{k}=L\big[T(x)\big]^{\frac{1}{2}}$. The reason for introducing such a transform will be seen immediately afterwards.$\hfill\blacksquare$\\

\noindent Of interest now is what happens if the transform is applied to $\lambda_{t}$, that is, $T(\lambda_{t})$. By It\^{o}'s lemma, we have:
\begin{align}\nonumber
    dT(\lambda_{t})&=\Big\{\frac{\partial T}{\partial t}+\frac{\partial T}{\partial \lambda_{t}}(a-b\lambda_{t})+\frac{\sigma^{2}}{2}\frac{\partial^{2} T}{\partial \lambda^{2}_{t}}(\lambda_{t})^{2k}\Big\}dt+\sigma \frac{\partial T}{\partial \lambda_{t}}(\lambda_{t})^{k}dW_{t}\\
    &=\Big\{0+\frac{aL^{2}}{2(1-k)}(\lambda_{t})^{1-2k}-\frac{bL^{2}}{2(1-k)}(\lambda_{t})^{2-2k}+\frac{\sigma^{2}}{2}\frac{L^{2}(1-2k)}{2(1-k)}\Big\}dt+\sigma \frac{L^{2}}{2(1-k)}(\lambda_{t})^{1-2k}(\lambda_{t})^{k}dW_{t}\nonumber\\
    &=\Big\{\frac{aL^{2}}{2(1-k)}(\lambda_{t})^{1-2k}-\frac{bL^{2}}{2(1-k)}(\lambda_{t})^{2-2k}+\frac{\sigma^{2}}{2}\frac{L^{2}(1-2k)}{2(1-k)}\Big\}dt+\sigma L \frac{L}{2(1-k)}(\lambda_{t})^{1-k}dW_{t}.
\end{align}
Based on the fact $\frac{L}{1-k}$ must be strictly positive-valued, we notice that we also need to assume $\sigma L>0$ so that the diffusion coefficient is positive-valued. Thus, we finally realize that we must set $L>0$, and subsequently $k<1$. Let us consider it from a different perspective: If $L<0$, $\sigma L$ in (2.3) will be negative-valued, which is undesirable and falls outside the prescribed form of our transformation, as the diffusion term $\sigma L \frac{L}{2(1-k)}(\lambda_{t})^{1-k}$ will be negative-valued (see later discussion about the CIR model). Also, when $k=\frac{1}{2}$ and $2a<\sigma^{2}$, zero is an attainable boundary of the original CKLS process (recall \textbf{Theorem 1.2 (4)}). However, the target CIR process constructed below, with constant drift term $\sigma^{2}L^{2}/2$, satisfies Feller's condition at equality and cannot reach zero from a positive initial value. Since $T(0)=0$, this is incompatible with equivalence of the two measures on any finite horizon over which the original process has a positive probability of hitting zero. Therefore, the case $k=\frac{1}{2}$ with $2a<\sigma^{2}$ is also excluded from the present equivalent-measure construction. Therefore, in what follows, we will keep assuming that:
\begin{equation}
    L>0 \text{ and } \frac{1}{2}\leq k<1 \text{ (when }k=\frac{1}{2},~2a\geq \sigma^{2}).
\end{equation}

\noindent Now according to (2.1), (2.3) becomes:
\begin{align}\nonumber
    dT(\lambda_{t})=\Big\{\frac{aL^{2}}{2(1-k)}(\lambda_{t})^{1-2k}-\frac{bL^{2}}{2(1-k)}(\lambda_{t})^{2-2k}+\frac{\sigma^{2}L^{2}(1-2k)}{4(1-k)}\Big\}dt+\sigma L \big[T(\lambda_{t})\big]^{\frac{1}{2}}dW_{t}.
\end{align}
Define $X_{t}\overset{\text{def}}{=}T(\lambda_{t})$, we have:
\begin{align}\nonumber
   \lambda_{t}=&T^{-1}(X_{t})=\Big[\frac{2(1-k)}{L}\Big]^{\frac{1}{1-k}}(X_{t})^{\frac{1}{2(1-k)}},\\
    &(\lambda_{t})^{2-2k}=\Big[\frac{2(1-k)}{L}\Big]^{2} X_{t},\nonumber\\
    &X_{t}=\Big[\frac{L}{2(1-k)}\Big]^{2}(\lambda_{t})^{2-2k}.
\end{align}
As a result: 
\begin{align}\nonumber
    dX_{t}&=\Big\{\frac{aL^{2}}{2(1-k)}(\lambda_{t})^{1-2k}-\frac{bL^{2}}{2(1-k)}\Big[\frac{2(1-k)}{L}\Big]^{2}X_{t}+\frac{\sigma^{2}L^{2}(1-2k)}{4(1-k)}\Big\}dt+\sigma L(X_{t})^{\frac{1}{2}}dW_{t}\\
    &=\Big\{\frac{\sigma^{2}L^{2}(1-2k)}{4(1-k)}-2b(1-k) X_{t}+\frac{aL^{2}}{2(1-k)}(\lambda_{t})^{1-2k}\Big\}dt+\sigma L(X_{t})^{\frac{1}{2}}dW_{t}.
\end{align}
The expression of (2.6) is rather tedious and hard to cope with. Indeed, in spite of the fact that the diffusion term of this SDE is already of the form of a CIR model (i.e. the exponent value of the process $X_{t}$ is $\frac{1}{2}$, see later discussions), the drift term still contains both $X_{t}$ and $\lambda_{t}$, that is:
\begin{equation}\nonumber
    \Big\{-2b(1-k)X_{t}+\frac{aL^{2}}{2(1-k)}(\lambda_{t})^{1-2k}\Big\}dt=\Bigg\{2b(k-1)X_{t}+\frac{aL^{2}}{2(1-k)}\Bigg[\Big[\frac{2(1-k)}{L}\Big]^{2} X_{t}\Bigg]^{1-\frac{1}{2}\frac{1}{1-k}}\Bigg\}dt
\end{equation}
is not a linear transform of $X_{t}$ with some constant coefficients when $\frac{1}{2}<k<1$. To this end, making the expression more concise and more easily applicable should be desired. We may refer to using the measure transform technique through the help of the Cameron-Martin-Girsanov-Maruyama theorem.\\

\noindent \textbf{Theorem 2.2} (One-dimensional Cameron-Martin-Girsanov-Maruyama theorem)\\
Let $W_{t}$ be a Wiener process in some filtered probability space $(\Omega,\mathcal{F}_{T},\{\mathcal{F}_{t}\}_{t\in[0,T]},\mathbb{P})$. Let $\theta_{t}=\{\theta_{t}\}_{t\in[0,T]}$ be a progressively measurable process with respect to $\mathcal{F}_{t}$, and assume that $\int_{0}^{T}(\theta_{s})^{2}ds<\infty~\mathbb{P}-a.s.$. Define a stochastic process $\mathcal{E}(\theta)_{t}$ on the same filtered probability space as (called the Dol\'{e}ans-Dade exponential or stochastic exponential of $\theta$ with respect to $W$):
\begin{align}\nonumber
    \mathcal{E}(\theta)_{t}&\overset{\text{def}}{=}\text{exp}\Big\{\int_{0}^{t}\theta_{s}dW_{s}-\frac{1}{2}\int_{0}^{t}(\theta_{s})^{2}ds\Big\},\\
    ~\text{and}~\tilde{W}_{t}&\overset{\text{def}}{=}W_{t}-\int_{0}^{t}\theta_{s}ds\nonumber.
\end{align}
When certain conditions are fulfilled\footnote{Note that this is crucial to our problem, which is also the main result in this paper. There can be various conditions for the justification.}, such as Novikov's condition or Kazamaki's condition:
\begin{equation}\nonumber
    \text{(Novikov) }\mathbb{E}^{\mathbb{P}}\Big[\text{exp}\Big\{\frac{1}{2}\int_{0}^{T}(\theta_{s})^{2}ds\Big\}\Big]<+\infty,~\text{(Kazamaki) }\Big\{\text{exp}\Big\{\frac{1}{2}\int_{0}^{t}\theta_{s}dW_{s}\Big\}\Big\}_{0\leq t\leq T}\text{ is a submartingale},
\end{equation}
for any $T>0$. Then $\mathbb{E}[\mathcal{E}(\theta)_{T}]=1$, and $\mathcal{E}(\theta_{t})$ is a (true) martingale with respect to $\mathbb{P}$. If so, $\mathcal{E}(\theta)_{t}$ is called the Radon-Nikod\'{y}m derivative, and a probability measure $\mathbb{Q}$ can be defined on $(\Omega,\mathcal{F}_{T})$ such that:
\begin{equation}\nonumber
    \frac{d\mathbb{Q}}{d\mathbb{P}}\Big\lvert_{\mathcal{F}_{t}}=\mathcal{E}(\theta)_{t},
\end{equation}
and the relationship between the two probability measures $\mathbb{P}$ and $\mathbb{Q}$ is:
\begin{equation}\nonumber
    \mathbb{Q}(B)=\mathbb{E}^{\mathbb{P}}[\mathcal{E}(\theta)_{t}\mathbbm{1}_{B}]=\mathbb{E}^{\mathbb{P}}\Big[\text{exp}\Big\{\int_{0}^{t}\theta_{s}dW_{s}-\frac{1}{2}\int_{0}^{t}(\theta_{s})^{2}ds\Big\}\mathbbm{1}_{B}\Big],~\forall B\in\mathcal{F}_{t}.
\end{equation}
In addition, the process $\tilde{W}_{t}$ is a $\mathbb{Q}$-Wiener process in the filtered probability space $(\Omega,\mathcal{F}_{T},\{\mathcal{F}_{t}\}_{t\in[0,T]},\mathbb{Q})$.\\

\noindent \textit{Proof}: Proofs of this theorem appear in a considerable amount of literature. Here, we cite only two representative references, see \hyperlink{Karatzas2012}{Karatzas and Shreve (2012)} [Chapter 3 \S 5 pages 190-198] or \hyperlink{Baxter1996}{Baxter and Rennie (1996)} [Chapter 3 \S 4 pages 63-76]. The argument hinges on constructing an absolutely continuous measure $\mathbb{Q}\ll\mathbb{P}$ through a Radon-Nikod\'{y}m density given by a Dol\'{e}ans-Dade exponential martingale, validating it with the Novikov or Kazamaki condition (i.e. the canonical formulation), applying It\^{o}'s lemma plus martingale properties under localized stopping times, and, optionally, using the martingale representation theorem. \qed\\

\noindent \textbf{Lemma 2.3}\\
Consider a stochastic process $Z_{t}$ defined on some filtered probability space $(\Omega,\mathcal{F},\{\mathcal{F}_{t}\}_{t\in[0,T]},\mathbb{P})$. Suppose the SDE of interest has the following expression:
\begin{equation}\nonumber
    dZ_{t}=A(Z_{t})dt+B(Z_{t})dW_{t},
\end{equation}
with $B(Z_{t})\neq 0$ for $t\in[0,T]$. Assume that under an equivalent probability measure $\mathbb{Q}$, where $\tilde{W}_{t}$ denotes the Wiener process under $\mathbb{Q}$.\\

\noindent The drift term of $Z_{t}$ can be changed to $\tilde{A}(Z_{t})$ from $A(Z_{t})$ as a direct result of the application of \textbf{Theorem 2.2} in the following way:
\begin{align}\nonumber
    dZ_{t}&=A(Z_{t})dt+B(Z_{t})dW_{t}=\tilde{A}(Z_{t})dt+B(Z_{t})\Big(\frac{A(Z_{t})-\tilde{A}(Z_{t})}{B(Z_{t})}\Big)dt+B(Z_{t})dW_{t}
    \nonumber\\
    &=\tilde{A}(Z_{t})dt+B(Z_{t})d\Big(W_{t}-\int_{0}^{t}-\frac{A(Z_{s})-\tilde{A}(Z_{s})}{B(Z_{s})}ds\Big)=\tilde{A}(Z_{t})dt+B(Z_{t})d\tilde{W}_{t},\nonumber
\end{align}
with $\tilde{W}_{t}\overset{\text{def}}{=}W_{t}-\int_{0}^{t}q_{s}ds$ where $q_{t}\overset{\text{def}}{=}-\frac{A(Z_{t})-\tilde{A}(Z_{t})}{B(Z_{t})}$. If some conditions such as Novikov's or Kazamaki's are fulfilled, then by \textbf{Theorem 2.2}, $\tilde{W}_{t}$ is a $\mathbb{Q}$-Wiener process $(\Omega,\mathcal{F},\{\mathcal{F}_{t}\}_{t\in[0,T]},\mathbb{Q})$ where $\mathbb{Q}$ is defined as:
\begin{align}\nonumber
    \frac{d\mathbb{Q}}{d\mathbb{P}}\Big\lvert_{\mathcal{F}_{t}}&=\text{exp}\Big\{\int_{0}^{t}q_{s}dW_{s}-\frac{1}{2}\int_{0}^{t}(q_{s})^{2}ds\Big\},\\
    \text{ and }\mathbb{Q}(B)&=\mathbb{E}^{\mathbb{P}}\Big[\text{exp}\Big\{\int_{0}^{t}q_{s}dW_{s}-\frac{1}{2}\int_{0}^{t}(q_{s})^{2}ds\Big\}\mathbbm{1}_{B}\Big],~\forall B\in\mathcal{F}_{t}\nonumber.
\end{align}

\noindent \textbf{Remark 2.4}\\
(1) The Cameron-Martin theorem has been progressively expanded into broader contexts by several authors, including \hyperlink{Maruyama1954}{Maruyama (1954)} and \hyperlink{Maruyama1955}{Maruyama (1955)}, \hyperlink{Girsanov1960}{Girsanov (1960)}, and \hyperlink{VanSchuppen1974}{Van Schuppen and Wong (1974)}, etc. In this context, we keep using the nomenclature Cameron-Martin-Girsanov-Maruyama theorem when referring to this theorem, rather than the Girsanov-Van Schuppen-Wong theorem.\\
(2) As already mentioned, the primarily used martingale criteria were developed by \hyperlink{Novikov1972}{Novikov (1972)} and \hyperlink{Kazamaki1977}{Kazamaki (1977)}. However, in practice, neither Novikov's nor Kazamaki's condition is easy to check. Writing $N_{t}=\int_{0}^{t}f(Z_{u})dW_{u}$, $\langle N\rangle_{t}=\int_{0}^{t}[f(Z_{u})]^{2}du$, and assuming that the stochastic integral is well defined, Novikov's condition requires $\mathbb{E}\big[\text{exp}\{\frac12\langle N\rangle_{T}\}\big]<\infty$. Kazamaki's criterion instead requires the process $\text{exp}\{N_{t}/2\}$ to be a submartingale on $[0,T]$. In concrete financial models, one may have only rough moment bounds for the quantities involved. Such bounds do not by themselves establish the exponential integrability or submartingale property required by the respective criteria. Consequently, even when a practitioner strongly suspects that the stochastic exponential is a true martingale, these criteria may still be difficult to verify, motivating the search for alternative criteria expressed directly in the model's drift and volatility functions.\\
(3) Here we mention 3 lesser‑known alternatives to the Novikov and Kazamaki conditions. In some particular situations, they may be more convenient to use. Because they are peripheral to our main line of argument, we will not elaborate on their precise proofs. \textbf{The first one} is called the Novikov-Krylov condition (\hyperlink{Krylov2002}{Krylov, 2002}), which reads: Fix $T>0$. Let $\theta_{t}$ be a progressively measurable process satisfying $\int_{0}^{T}\theta_{s}^{2}ds<\infty$ $\mathbb{P}-a.s.$, and set $N_{t}=\int_{0}^{t}\theta_{s}dW_{s}$, $0\leq t\leq T$. Then $N$ is a continuous local martingale starting at zero. Assume that
\begin{equation}\nonumber
   \lim_{\varepsilon\to 0+}\varepsilon\text{log}\mathbb{E}\Big[\text{exp}\Big\{\frac{1-\varepsilon}{2}\int_{0}^{T}(\theta_{s})^{2}ds\Big\}\Big]=0.
\end{equation}
Then $\mathbb{E}[\mathcal{E}(\theta)_{T}]=\mathbb{E}[\text{exp}\{N_{T}-\frac{1}{2}\langle N\rangle_{T}\}]=1$. In particular, the conclusion holds whenever Novikov's condition is satisfied. \textbf{The second one} is known as another Novikov's type condition. One may refer to Exercise 1.40 in \hyperlink{Revuz2013}{Revuz and Yor (2013)} [Chapter VIII, page 338] for this. Let $W_{t}$, $t\geq 0$ be a standard Wiener process, $H_{t}$, $t\geq 0$ a predictable process and fix $T>0$. Set $N_{t}\overset{\text{def}}{=}\int_{0}^{t} H_{s}dW_{s}$ for $t\in[0,T]$. Assume there exist constants $A,C>0$ such that
\begin{equation}\nonumber
   \mathbb{E}\Big[\text{exp}\Big\{A\lvert H_{t}\lvert^{2}\Big\}\Big]\leq C,~\forall t\in[0,T].
\end{equation}
Then $\mathbb{E}[\text{exp}\{N_{T}-\frac{1}{2}\int_{0}^{T}(H_{s})^{2}ds\}]=1$. Typical examples of such $H_{t}$ include $H_{t}=b(W_{t})$ when $b(\cdot)$ has at most linear growth, but also any predictable Gaussian process whose mean and variance are bounded on $[0,T]$ (e.g. $H_{t}=\widetilde{W}_{t}$ where $\widetilde{W}$ is an independent standard Wiener process). \textbf{The third one} is called Bene\v{s}'s condition (see \hyperlink{Liptser2013}{Liptser (2013)}): Let $\mathcal{E}_{t}$ denote the solution to the Dol\'{e}ans-Dade equation: $\mathcal{E}_{t}=1+\int_{0}^{t}\mathcal{E}_{s}[\theta_{s}]dW_{s}$ for $t\in[0,T]$, where $W_{t}$ is a standard Wiener process and $\theta_{t}$ is a progressively Borel measurable process with $\int_{0}^{t}(\theta_{s})^{2}ds<+\infty$ almost surely. The process $\mathcal{E}_{t}$ is a martingale provided there exists some constant $K$ such that
\begin{equation}\nonumber
    \lvert\theta_{t}\lvert^{2}\leq K\big[1+\sup_{s\in[0,t]}(W_{s})^{2}\Big],~\forall t\in[0,T].
\end{equation}
(4) A proof of Novikov's condition, Kazamaki's condition, and that Kazamaki's condition is a necessary but not sufficient condition for Novikov's condition is given at the second-to-last part of the Appendix.$\hfill\blacksquare$\\

\noindent Now since 
\begin{equation}\nonumber
    dX_{t}=\Big\{\frac{\sigma^{2}L^{2}(1-2k)}{4(1-k)}-2b(1-k) X_{t}+\frac{aL^{2}}{2(1-k)}(\lambda_{t})^{1-2k}\Big\}dt+\sigma L(X_{t})^{\frac{1}{2}}dW_{t},
\end{equation}
we may want to get rid of $\lambda_{t}$ in the drift term and want the new drift term to be $\frac{\sigma^{2}L^{2}}{2}-2b(1-k)X_{t}$, so that the expression will become much more natural and analytically friendly: the sum of a constant and $X_{t}$ multiplied by a $k$-dependent coefficient.\footnote{An extravagant hope is that we may even get rid of $X_{t}$ in the drift term, yet since the expression $q_{t}=-\frac{A(Z_{t})-\tilde{A}(Z_{t})}{B(Z_{t})}$ in the previous lemma would contain $X_{t}$ as well if we insist doing so, making the expression of $q_{t}$ even more complicated and tricky. Therefore, it would be more wise if we can just keep $-2b(1-k)X_{t}$ in the expression of the drift term.} As a result, based on the fact that the diffusion term is $\sigma$ times $L$ times the square root of $X_{t}$, modifying the drift term in a way like this would possibly lead to the so-called Cox-Ingersoll-Ross (CIR) model describing the evolution of $r_{t}$ (called the Feller square-root process), which is a mean-reverting process as well, defined on some filtered probability space $(\Omega,\mathcal{F},\{\mathcal{F}_{t}\}_{t\in[0,T]},\mathbb{P})$ with the mean-reversion speed $a^{*}\in\mathbb{R}_{+}$, the long-term mean level $b^{*}\in\mathbb{R}_{+}$, the diffusion coefficient $\sigma^{*}\in\mathbb{R}_{+}$, and the initial value $r_{0}\in\mathbb{R}_{+}$\footnote{Note that the structure of the drift term is $a^{*}(b^{*}-r_{t})dt$, which is different from CKLS one where $(a-b\lambda_{t})dt$ is the structure of the drift term.}:
\begin{equation}
    dr_{t}=a^{*}(b^{*}-r_{t})dt+\sigma^{*}(r_{t})^{\frac{1}{2}}dW_{t}.
\end{equation}
More details of the CIR model will be explained in Section 3.\\

\noindent Consequently, we obtain $q_{t}$ and $\tilde{W}_{t}$ in our case:
\begin{align}\nonumber
    q_{t}&=q(\lambda_{t})=-\frac{\big[\frac{\sigma^{2}L^{2}(1-2k)}{4(1-k)}-2b(1-k) X_{t}+\frac{aL^{2}}{2(1-k)}(\lambda_{t})^{1-2k}\big]-\big[\frac{\sigma^{2}L^{2}}{2}-2b(1-k)X_{t}\big]}{\frac{\sigma L^{2}}{2(1-k)}(\lambda_{t})^{1-k}}\\
    &=-\frac{\big(\frac{\sigma^{2}L^{2}(1-2k)}{4(1-k)}-\frac{\sigma^{2}L^{2}}{2}\big)+\frac{aL^{2}}{2(1-k)}(\lambda_{t})^{1-2k}}{\frac{\sigma L^{2}}{2(1-k)}(\lambda_{t})^{1-k}}=-\frac{-\frac{\sigma}{2}+\frac{a}{\sigma}(\lambda_{t})^{1-2k}}{(\lambda_{t})^{1-k}}\nonumber\\
    &=\frac{\sigma}{2}(\lambda_{t})^{k-1}-\frac{a}{\sigma}(\lambda_{t})^{-k},\nonumber\\
    \text{and~}\tilde{W}_{t}&=W_{t}-\int_{0}^{t}q_{s}ds.\nonumber
\end{align}
We define the Dol\'{e}ans-Dade exponential (not yet being a Radon-Nikod\'{y}m derivative until its martingality is proven):
\begin{equation}
    M_{t}\overset{\text{def}}{=}\text{exp}\Big\{\int_{0}^{t}q_{s}dW_{s}-\frac{1}{2}\int_{0}^{t}(q_{s})^{2}ds\Big\}.
\end{equation}

\noindent If we could manage to prove that the Dol\'{e}ans-Dade exponential $M_{t}$ is a martingale with respect to the original probability measure $\mathbb{P}$ (e.g. successfully verifying that Novikov's or Kazamaki's condition is satisfied hence $\mathbb{E}^{\mathbb{P}}[M_{T}]=1$), then the one-dimensional Cameron-Martin-Girsanov theorem implies that $M_{t}$ serves as the Radon-Nikod\'{y}m derivative process $\frac{d\mathbb{Q}}{d\mathbb{P}}\lvert_{\mathcal{F}_{t}}=M_{t}$ for $t\in[0,T]$ defining an equivalent measure $\mathbb{Q}$ on the same filtered space $(\Omega,\mathcal{F},\{\mathcal{F}_{t}\}_{t\in[0,T]})$ with $\tilde{W}_{t}$ being a $\mathbb{Q}$-Wiener process. And most importantly, under the equivalent probability measure $\mathbb{Q}$, $X_{t}$ will admit the following CIR dynamics as we desire:
\begin{equation}
    dX_{t}=\Big(\frac{\sigma^{2} L^{2}}{2}-2b(1-k)X_{t}\Big)dt+\sigma L(X_{t})^{\frac{1}{2}}d\tilde{W}_{t}.
\end{equation}

\noindent Recall after obtaining (2.3), we have emphasized that $L>0$ and $\frac{1}{2}\leq k<1$ (when $k=\frac{1}{2}$, $2a\geq \sigma^{2}$) should be met, so that $\lambda_{t}$ is always strictly positive almost surely, together with $\sigma>0$, $L>0$, $X_{t}>0$, making the diffusion term $\sigma L(X_{t})^{\frac{1}{2}}$ positive-valued. Indeed, in the context of the CIR dynamics (2.7), the parameter $\sigma^{*}=\sigma L$ should be strictly positive-valued, which means that $L$ should be positive in our case, forcing $k<1$; otherwise, (2.9) will have a negative-valued diffusion coefficient. In addition, $a^{*}$ in the model is assumed to be positive-valued as well. In our case, we may let $a^{*}=2b(1-k)$ and check if $a^{*}>0$ is satisfied. Here $b>0$ is the general assumption of the CKLS model (1.1); therefore, when $\frac{1}{2}\leq k<1$, the positivity of $a^{*}$ is easily checked. $\frac{\sigma^{2} L^{2}}{2}$ is always positive-valued, $b^{*}=\frac{\sigma^{2} L^{2}}{2a^{*}}=\frac{\sigma^{2}L^{2}}{4b(1-k)}$ is always positive-valued, which means: When the following parameter setting (2.10) is assumed, (2.9) indeed corresponds to a CIR model.\\

\noindent \textbf{Remark 2.5}\\
\noindent The choice of the constant drift term is motivated by boundary compatibility and analytical tractability. Since the original CKLS process is strictly positive and $T$ preserves positivity, an equivalent change of measure must preserve this property on every finite time interval. For the target CIR dynamics $dX_{t}=(\bar{a}-2b(1-k)X_{t})dt +\sigma L\sqrt{X_{t}}d\tilde{W}_{t}$, Feller's condition for strict positivity from a positive initial value is $2\bar{a}\geq\sigma^{2}L^{2}$. We therefore choose $\bar{a}=\sigma^{2}L^{2}/2$, the smallest constant satisfying this condition. Although larger values are also permitted by Feller's condition, this critical choice offers several simplifications. First, the Girsanov kernel becomes $q(x)=\frac{\sigma}{2}x^{k-1}-\frac{a}{\sigma}x^{-k}$, so the coefficient of $x^{k-1}$ is independent of $k$. Second, the associated auxiliary diffusion has scale density proportional to $x^{-1}\text{exp}\{\frac{b}{\sigma^{2}(1-k)}x^{2(1-k)}\}$, whose integral diverges logarithmically at zero, facilitating the boundary analysis underlying the martingale proof. Third, the stationary Gamma distribution of the resulting CIR model has shape parameter $2\bar{a}/(\sigma^{2}L^{2})=1$, and hence reduces to an exponential distribution. Through $T^{-1}$, this yields a Weibull stationary distribution for the corresponding $\lambda$-dynamics under the transformed measure, simplifying the subsequent distributional analysis. As proved in Section 4, the associated stochastic exponential is a true martingale on every finite time interval, thereby validating the equivalent change of measure.$\hfill\blacksquare$\\

\noindent In the following sections, we always assume that\footnote{The scenario where $\frac{1}{2}<k<1$ is often referred to as the predominant case in the context of the CKLS model. Recall the contents in the literature review part of this paper: A significant number of academic studies have indicated that this particular case is frequently cited in empirical finance research.}:
\begin{align}
    \frac{1}{2}\leq k<1 \text{ (when }k=\frac{1}{2},~&2a\geq \sigma^{2}),~L>0,~a~>0,~b>0,~\sigma>0;\\
    a^{*}=2b(1-k)>0,~&b^{*}=\frac{\sigma^{2}L^{2}}{4b(1-k)}>0,~\sigma^{*}=\sigma L>0.\nonumber
\end{align}

\noindent For a clearer explanation, we present a flowchart (Figure \ref{fig:ckls-firststyle}) illustrating the key steps of the whole procedure  :
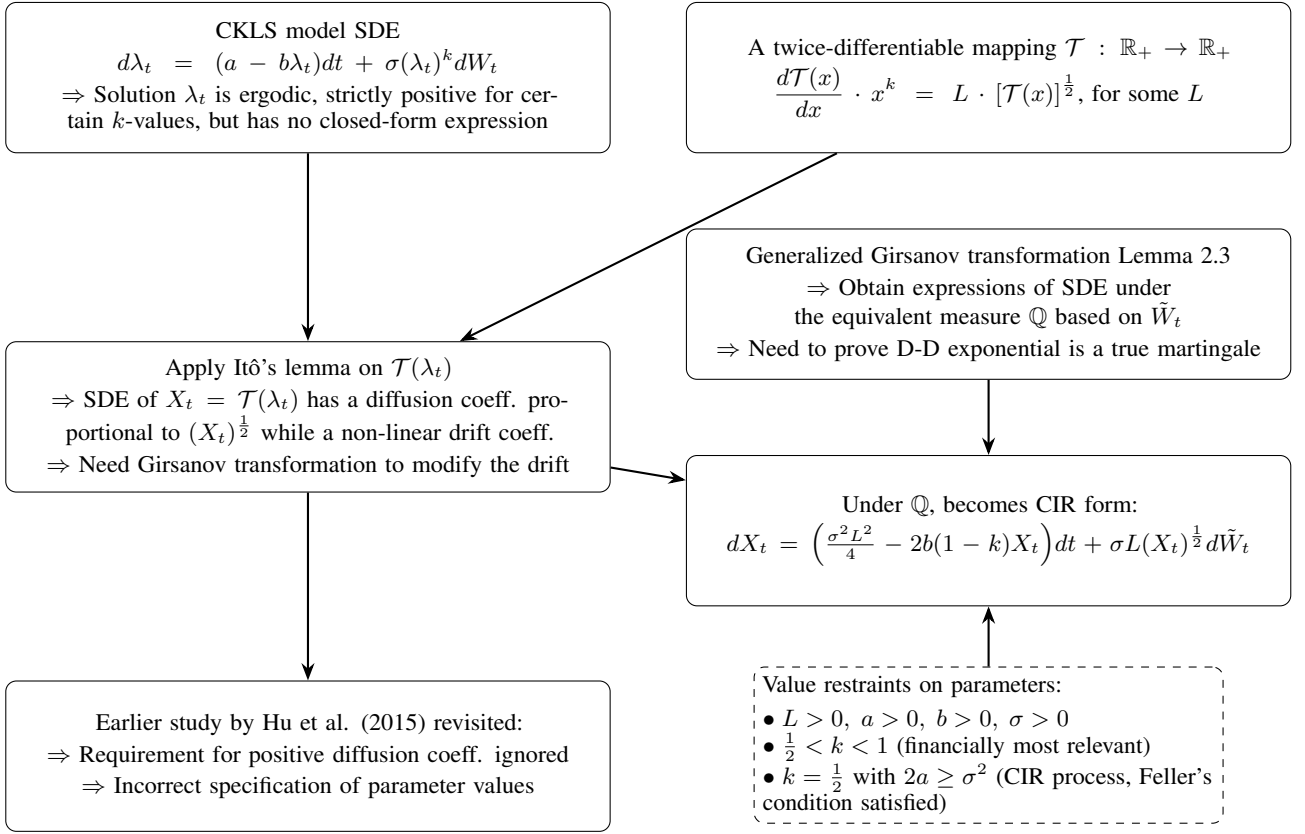
\begin{figure}[htbp]
\centering
\begin{tikzpicture}[
  node distance=2.0cm and 1.0cm,
  proc/.style={draw, rounded corners, align=center,
                  font=\small, 
                  text width=7.8cm,
                  minimum height=2cm},
  note/.style={proc, draw, dashed, align=left, text width=6.0cm},
  arrow/.style={-Stealth, thick}
]

\node[proc] (ckls) {CKLS model SDE\\[2pt]
  $d\lambda_{t}=(a-b\lambda_{t})dt+\sigma(\lambda_{t})^{k}dW_{t}$\\[2pt]
  $\Rightarrow$ Solution $\lambda_{t}$ is ergodic, strictly positive for certain $k$-values, but has no closed-form expression};
  
\node[proc, right=of ckls] (map) {A twice-differentiable mapping $T:\mathbb{R}_{+}\!\to\!\mathbb{R}_{+}$\\[2pt]
  $\displaystyle \frac{dT(x)}{dx}\cdot x^{k}=L\cdot[T(x)]^{\frac{1}{2}}$, for some $L$};

\node[proc, below=2.5cm of ckls] (ito) {Apply It\^{o}'s lemma on $T(\lambda_{t})$\\[2pt]
  $\Rightarrow$ SDE of $X_{t}=T(\lambda_{t})$ has a diffusion coeff. proportional to $(X_{t})^{\frac{1}{2}}$ while a non-linear drift coeff.\\[2pt]
  $\Rightarrow$ Need Girsanov transformation to modify the drift};
  
\node[proc, below=1.0cm of map] (girs) {Generalized Girsanov transformation \textbf{Lemma 2.3}\\[2pt]
  $\Rightarrow$ Obtain expressions of SDE under the equivalent measure $\mathbb{Q}$ based on $\tilde{W}_{t}$\\[2pt]
  $\Rightarrow$ Need to prove D-D exponential is a true martingale};

\node[proc, below=2.5cm of ito] (err) {Earlier study by Hu et al. (2015) revisited:\\[2pt]
  $\Rightarrow$ Requirement for positive diffusion coeff. ignored\\[2pt]
  $\Rightarrow$ Incorrect specification of parameter values};
  
\node[proc, below=1.0cm of girs] (cir) {Under $\mathbb{Q}$, becomes CIR form:\\[2pt]
  $dX_{t}=\Big(\frac{\sigma^{2}L^{2}}{2}-2b(1-k)X_{t}\Big)dt+\sigma L(X_{t})^{\frac{1}{2}}d\tilde{W}_{t}$};

\node[note, below=0.8cm of cir] (cond) {
  Value restraints on parameters:\\[3pt]
  $\bullet\;L>0,\;a>0,\;b>0,\;\sigma>0$\\
  $\bullet\;\frac{1}{2}<k<1$ \text{(financially most relevant)}\\
  $\bullet\;k=\frac{1}{2}$ with $2a\ge\sigma^{2}$ (CIR process, Feller's condition satisfied)
};

\draw[arrow] (ckls) -- (ito);
\draw[arrow] (map)  -- (ito);
\draw[arrow] (ito)  -- (cir);
\draw[arrow] (ito)  -- (err);
\draw[arrow] (girs) -- (cir);
\draw[arrow] (cond) -- (cir);

\end{tikzpicture}
\caption{From CKLS to CIR: Using a Twice-differentiable Mapping and Generalized Girsanov's Theorem.}
\label{fig:ckls-firststyle}
\end{figure}

\section{Subsidiary result: Properties of \texorpdfstring{$X_{t}$}{Xt} and \texorpdfstring{$\lambda_{t}$}{lambdat} under the equivalent probability measure based on general theories of Cox-Ingersoll-Ross model}
\noindent Let the CIR process (the Feller square-root process) $r_{t}$ (the expression of which will be discussed soon later) be the solution to the SDE (2.7); we have:\\

\noindent \textbf{Lemma 3.1} (Feller's condition)\\
(1) If $2a^{*}b^{*}\geq\sigma^{*2}$, the process $r_{t}$ will be strictly positive-valued with probability one. That is, it will never hit 0 in a finite time: $\mathbb{P}(\tau^{r}_{0}=+\infty)=1$ where $\tau^{r}_{0}\overset{\text{def}}{=}\text{inf}\{t\geq 0\lvert r_{t}=0\}$.\\
(2) If $2a^{*}b^{*}<\sigma^{*2}$, then, starting from $r_{0}>0$, the process $r_{t}$ reaches zero in finite time almost surely: $\mathbb{P}(\tau_{0}^{r}<+\infty)=1$. The boundary at zero is instantaneously reflecting.

\noindent \textit{Proof}: See Appendix. \qed\\

\noindent \textbf{Remark 3.2}\\
(1) The diffusion term $\sigma^{*}(r_{t})^{\frac{1}{2}}$ prevents interest rates from becoming negative-valued for all positive values of $a^{*}$ and $b^{*}$.\\
(2) The reason why $2a^{*}b^{*}\geq\sigma^{*2}$ or not matters is that: If so, as the rate $r_{t}$ approaches zero, the level-dependent diffusion term $\sigma^{*}(r_{t})^{\frac{1}{2}}$ diminishes significantly, reducing the impact of random shocks on the rate. As a result, when the rate nears zero, its movement is dominantly determined by the drift, driving the rate upward to a state of equilibrium.\\
(3) The CIR process is also used to model the instantaneous variance in the Heston stochastic volatility model.$\hfill\blacksquare$\\

\noindent In our case, according to (2.5), $X_{t}$ is the square of $\frac{L}{2(1-k)}(\lambda_{t})^{1-k}$, and since we have already assumed that the initial value of $\lambda_{0}>0$ in (1.1), it is clear that $X_{t}\lvert_{t=0}>0$ and $X_{t}\lvert_{t>0}>0$. Under the new probability measure $\mathbb{Q}$, the updated CIR parameters satisfy
\begin{equation}\nonumber
    2a^{*}b^{*}=2\big[2b(1-k)\big]\frac{\sigma^{2}L^{2}}{4b(1-k)}=\sigma^{2}L^{2}=\sigma^{*2}.
\end{equation}
Therefore, by \textbf{Lemma 3.1 (1)}, the process $X_{t}$, starting from $X_{0}>0$, remains strictly positive almost surely under $\mathbb{Q}$. This is consistent with the equivalence of $\mathbb{P}$ and $\mathbb{Q}$ on every finite time interval.\\

\noindent We also record the dynamics of $\lambda_{t}$ and its power transform under $\mathbb{Q}$. Fix any finite time horizon $T>0$. The validity of the equivalent change of measure on $[0,T]$ is established in Section~4. Since
$dW_{t}=d\widetilde{W}_{t}+q(\lambda_{t})dt$, substitution into the original CKLS equation yields
\begin{equation}\nonumber
    d\lambda_{t}=\big[a-b\lambda_{t}+\sigma\lambda_{t}^{k}q(\lambda_{t})\big]dt+\sigma\lambda_{t}^{k}d\widetilde{W}_{t}=\Big[\frac{\sigma^{2}}{2}\lambda_{t}^{2k-1}-b\lambda_{t}\Big]dt+\sigma\lambda_{t}^{k}d\widetilde{W}_{t}, \qquad 0\leq t\leq T.
\end{equation}
Furthermore, define $V_{t}:=(\lambda_{t})^{2-2k}=\frac{4(1-k)^{2}}{L^{2}}X_{t}$ according to (2.5). Multiplying both sides of (2.9) by $\frac{4(1-k)^{2}}{L^{2}}$ gives
\begin{equation}\nonumber
    dV_{t}=\big[2\sigma^{2}(1-k)^{2}-2b(1-k)V_{t}\big]dt+2\sigma(1-k)\sqrt{V_{t}}d\widetilde{W}_{t}, \qquad V_{0}=\lambda_{0}^{2-2k}.
\end{equation}
Thus, $V$ itself is a CIR process under $\mathbb{Q}$.

\noindent \textbf{Lemma 3.3}\\
The CIR model (2.7) has the exact solution (the Feller square-root process):
\begin{align}\nonumber
    r_{t}&=e^{-a^{*}t}r_{0}+b^{*}(1-e^{-a^{*}t})+\sigma^{*}e^{-a^{*}t}\int_{0}^{t}e^{a^{*}s}\big(r_{s}\big)^{\frac{1}{2}}dW_{s}.\nonumber
\end{align}
\noindent \textit{Proof}: See the Appendix.\qed\\
Note that, for any $0\leq t'<t\leq T$, the corresponding representation starting at time $t'$ is
\begin{align}\nonumber
    r_{t}=e^{-a^{*}(t-t')}r_{t'}+b^{*}(1-e^{-a^{*}(t-t')})+\sigma^{*}\int_{t'}^{t}e^{-a^{*}(t-s)}\big(r_{s}\big)^{\frac{1}{2}}dW_{s}.
\end{align}
We shall use this representation when needed.\\

\noindent In our case:
\begin{align}\nonumber
    X_{t}=&e^{2b(k-1)t}X_{0}+\frac{\sigma^{2}L^{2}}{4b(1-k)}[1-e^{2b(k-1)t}]+\sigma Le^{2b(k-1)t}\int_{0}^{t}e^{2b(1-k)s}\big(X_{s}\big)^{\frac{1}{2}}d\tilde{W}_{s}\\
    =&e^{2b(k-1)t}\Big[\frac{L^{2}}{4(1-k)^{2}}\big(\lambda_{0}\big)^{2-2k}\Big]+\frac{\sigma^{2}L^{2}}{4b(1-k)}[1-e^{2b(k-1)t}]+\sigma Le^{2b(k-1)t}\int_{0}^{t}e^{2b(1-k)s}\Big[\frac{L}{2(1-k)}\big(\lambda_{s}\big)^{1-k}\Big]d\tilde{W}_{s}\nonumber\\
    =&\frac{L^{2}\big(\lambda_{0}\big)^{2-2k}}{4(1-k)^{2}}e^{2b(k-1)t}+\frac{\sigma^{2}L^{2}}{4b(1-k)}-\frac{\sigma^{2}L^{2}}{4b(1-k)}e^{2b(k-1)t}+\frac{\sigma L^{2}}{2(1-k)}e^{2b(k-1)t}\int_{0}^{t}e^{2b(1-k)s}\big(\lambda_{s}\big)^{1-k}d\tilde{W}_{s}\nonumber\\
    =&\frac{\sigma^{2}L^{2}}{4b(1-k)}+\frac{L^{2}\big[b\big(\lambda_{0}\big)^{2-2k}-\sigma^{2}(1-k)\big]}{4b(1-k)^{2}}e^{2b(k-1)t}+\frac{\sigma L^{2}}{2(1-k)}e^{2b(k-1)t}\int_{0}^{t}e^{2b(1-k)s}\big(\lambda_{s}\big)^{1-k}d\tilde{W}_{s}.\nonumber
\end{align}
According to (2.5), $\lambda_{t}=\Big[\frac{2(1-k)}{L}\Big]^{\frac{1}{1-k}}(X_{t})^{\frac{1}{2(1-k)}}$. Using the definition of $V_{t}$ above, we conclude that
\begin{align}
    \lambda_{t}&=\Big[\frac{2(1-k)}{L}\Big]^{\frac{1}{1-k}}\Big\{\frac{\sigma^{2}L^{2}}{4b(1-k)}+\frac{L^{2}\big[b(\lambda_{0})^{2-2k}-\sigma^{2}(1-k)\big]}{4b(1-k)^{2}}e^{2b(k-1)t}+\frac{\sigma L^{2}}{2(1-k)}e^{2b(k-1)t}\int_{0}^{t}e^{2b(1-k)s}(\lambda_{s})^{1-k}d\tilde{W}_{s}\Big\}^{\frac{1}{2(1-k)}}\\
    V_{t}&=V_{0}e^{2b(k-1)t}+\frac{\sigma^{2}(1-k)}{b}\Big[1-e^{2b(k-1)t}\Big]+2\sigma(1-k)e^{2b(k-1)t}\int_{0}^{t}e^{2b(1-k)s}(\lambda_{s})^{1-k}d\tilde{W}_{s}.\nonumber
\end{align}

\noindent \textbf{Lemma 3.4}\\
Let $r_{t}$ be the solution to the CIR model (2.7). For $t>0$, its conditional distribution given $r_{0}$ can be computed in closed form. Without loss of generality, given the current value $r_{t}$, we always have the distribution of the future value $r_{t+t^{*}}$ with $t^{*}>0$. For simplicity, we let $t=0$ and $t^{*}=t$. To be specific: For $\gamma\overset{\text{def}}{=}\omega r_{t}$, define $R\overset{\text{def}}{=}2\omega r_{t}=2\gamma$, where $\omega\overset{\text{def}}{=}\frac{2a^{*}}{(1-e^{-a^{*}t})\sigma^{*2}}$. Then, $R$ is a non-central chi-squared distributed random variable with $2(\kappa+1)$ degrees of freedom where $\kappa\overset{\text{def}}{=}\frac{2a^{*}b^{*}}{\sigma^{*2}}-1$ and non-centrality parameter $2\theta$ and $\theta\overset{\text{def}}{=}\omega e^{-a^{*}t}r_{0}$. The transition density function of the process $r_{t}$ (the future value), given the value of $r_{0}$, is:
\begin{equation}\nonumber
    f(r_{t}\lvert r_{0},a^{*},b^{*},\sigma^{*})=\omega e^{-\theta-\gamma}\Big(\frac{\gamma}{\theta}\Big)^{\frac{\kappa}{2}}I_{\kappa}(2\sqrt{\theta\gamma}),\qquad r_{t}>0.
\end{equation}
where $I_{\kappa}(\cdot)$ is a modified Bessel function of the first kind of order $\kappa$: $I_{\kappa}(x)=(\frac{x}{2})^{\kappa}\sum_{n=0}^{+\infty}\frac{(x/2)^{2n}}{n!\Gamma(\kappa+n+1)}$ and the Gamma function $\Gamma(z)=\int_{0}^{\infty}t^{z-1}e^{-t}dt$, $Re(z)>0$, $z\in\mathbb{C}$.\\

\noindent \textit{Proof}: See the Appendix.\qed\\

\noindent In our case: $\omega=\frac{4b(1-k)}{[1-e^{2b(k-1)t}]\sigma^{2}L^{2}}$. $\kappa=\frac{\sigma^{2}L^{2}}{\sigma^{2}L^{2}}-1=0$, $2(\kappa+1)=2$, $\theta=\frac{4b(1-k)}{[1-e^{2b(k-1)t}]\sigma^{2}L^{2}}e^{2b(k-1)t}X_{0}=\frac{b(\lambda_{0})^{2-2k}}{(e^{2b(1-k)t}-1)\sigma^{2}(1-k)}$ and $\gamma=\frac{4b(1-k)X_{t}}{[1-e^{2b(k-1)t}]\sigma^{2}L^{2}}$.\\

\noindent By relationships $X_{t}=T(\lambda_{t})=\frac{L^{2}}{4(1-k)^{2}}\lambda_{t}^{2(1-k)},~T'(\lambda_{t})=\frac{L^{2}}{2(1-k)}\lambda_{t}^{1-2k},~X_{t}=\chi V_{t},~\chi:=\frac{L^{2}}{4(1-k)^{2}}$. The conditional density of $V_{t}=v$ given $V_{0}=v_{0}$ is
\begin{align}\nonumber
	f^{V}(v\lvert v_{0})&=\chi f^{X}\big(\chi v \lvert \chi v_{0}\big)=\chi \omega e^{-\gamma_{1}-\theta_{1}}I_{0}\big(2\sqrt{\gamma_{1}\theta_{1}}\big),\qquad v>0.
\end{align}
where $\gamma_{1}=\chi \omega v_{0}e^{-2b(1-k)t},~\theta_{1}=\chi \omega v$.\\

\noindent The conditional density of $\lambda_{t}=\ell$ given $\lambda_{0}=\ell_{0}$ is
\begin{align}\nonumber
	f^{\lambda}(\ell \lvert \ell_{0})=f^{X}\big(T(\ell)\lvert T(\ell_{0})\big)T'(\ell)=T'(\ell) \omega e^{-(\gamma_{2}+\theta_{2})}I_{0}\big(2\sqrt{\gamma_{2}\theta_{2}}\big),\quad \ell>0.
\end{align}
where $\gamma_{2}=\omega T(\ell_{0})e^{-2b(1-k)t},\quad \theta_{2}=\omega T(\ell),\quad T'(\ell)=\frac{L^{2}}{2(1-k)}\ell^{1-2k}$.

\noindent \textbf{Lemma 3.5}\\
\noindent (a) A CIR process $r_{t}$ can be represented in the following form:
\begin{equation}\nonumber
	r_{t}=e^{-a^{*}t}\mathrm{BESQ}_{(d,R_{0})}\Big(\frac{\sigma^{*2}}{4a^{*}} (e^{a^{*}t}-1)\Big),
\end{equation}
where $\mathrm{BESQ}_{(d,R_{0})}$ denotes a squared Bessel process starting from the initial point $R_{0}=r_{0}$ of dimension $d=\frac{4a^{*}b^{*}}{\sigma^{*2}}$. For the process $X$ in (2.9), this dimension equals $2$.\\

\noindent For parts (b) and (c), assume that $\mu\in\mathbb{R}$, $\gamma>0$, $\eta_{0}>0$, and $K>0$ with $K\neq 1$. Let $\tau_{0}=\inf\{t\geq 0:\eta_{t}=0\}$, with $\inf\varnothing=\infty$. The transformations below are understood for $0\leq t<\tau_{0}$.

\noindent (b) A CEV process $\eta_{t}$ solving the equation
\begin{equation}\nonumber
    d\eta_{t}=\mu\eta_{t}dt+\gamma(\eta_{t})^{K}dW_{t},
\end{equation}
can be represented as a power of a square-root diffusion. Indeed, setting $\delta=2(K-1)$, the process $(\eta_{t})^{-\delta}$ satisfies
\begin{equation}\nonumber
	d\Big(\frac{1}{(\eta_{t})^{\delta}}\Big)=\Big(\mathbf{A}-\mathbf{B}\frac{1}{(\eta_{t})^{\delta}}\Big)dt+\mathbf{\Gamma} \Big(\Big\lvert \frac{1}{(\eta_{t})^{\delta}}\Big\lvert \Big)^{\frac{1}{2}}dW_{t},
\end{equation}
where $\mathbf{A}=\frac{\delta(\delta+1)\gamma^{2}}{2}$, $\mathbf{B}=\delta\mu$, $\mathbf{\Gamma}=-\delta\gamma$. The result follows directly by applying It\^{o}'s lemma to $(\eta_{t})^{-\delta}$.\\

\noindent (c) A CEV process $\eta_{t}$ can be represented in the form
\begin{equation}\nonumber
	\eta_{t}=e^{\mu t}\mathrm{BESQ}_{\big(\frac{2K-1}{K-1}, \eta_{0}^{-2(K-1)}\big)}^{\frac{1}{2(1-K)}}\Big(\frac{(K-1)\gamma^{2}}{2\mu}\big(e^{2(K-1)\mu t}-1\big)\Big).
\end{equation}
{Here $\operatorname{BESQ}_{(d,z)}$ denotes a squared Bessel process of dimension $d$, starting from $z>0$, considered up to its first hitting time of zero, with $d=\frac{2K-1}{K-1}$. For $\frac{1}{2}<K<1$, this dimension is negative, and the squared Bessel dynamics are used only before the process reaches zero. No continuation beyond the zero boundary is asserted here. When $\mu=0$, the time-change argument is interpreted by continuity as $(K-1)^{2}\gamma^{2}t$.\\

\noindent \textit{Proof} See Appendix. See also e.g. \hyperlink{Delbaen2002}{Delbaen and Shirakawa (2002)}.\\

\noindent \textbf{Lemma 3.6}\\
The moments of the CIR $r_{t}$ are:
\begin{align}\nonumber
    \mathbb{E}[r_{t}]&=r_{0}e^{-a^{*}t}+b^{*}(1-e^{-a^{*}t}),\nonumber\\
    \text{Var}(r_{t})&=\frac{r_{0}\sigma^{*2}}{a^{*}}\Big(e^{-a^{*}t}-e^{-2a^{*}t}\Big)+\frac{b^{*}\sigma^{*2}}{2a^{*}}\Big(1-e^{-a^{*}t}\Big)^{2},\nonumber\\
    \text{Cov}(r_{t},r_{t'})&=\frac{r_{0}\sigma^{*2}}{a^{*}}\Big[e^{-a^{*}t'}-e^{-a^{*}(t+t')}\Big]+\frac{b^{*}\sigma^{*2}}{2a^{*}}\Big[e^{a^{*}(t-t')}+e^{-a^{*}(t+t')}-2e^{-a^{*}t'}\Big],~\text{for }0\leq t\leq t'.\nonumber
\end{align}
More generally, for $n\in\mathbb{N}$ and $t>0$, define $\iota=\frac{2a^{*}b^{*}}{\sigma^{*2}}=\kappa+1$, $s_{t}=\frac{\sigma^{*2}}{2a^{*}}\big(1-e^{-a^{*}t}\big)$, $h_{t}=r_{0}e^{-a^{*}t}$. Then
\begin{equation}\nonumber
    \mathbb{E}[(r_{t})^{n}]=\sum_{j=0}^{n}\binom{n}{j}\frac{\Gamma(\iota+n)}{\Gamma(\iota+j)}s_{t}^{n-j}h_{t}^{j}.
\end{equation}
At $t=0$, the formula is understood by continuity and gives $\mathbb{E}\big[(r_{0})^{n}\big]=(r_{0})^{n}$.\\

\noindent \textit{Proof}: See the Appendix.\qed\\

\noindent In our case:
\begin{align}\nonumber
    \mathbb{E}[X_{t}]&=\frac{L^{2}\big(\lambda_{0}\big)^{2-2k}}{4(1-k)^{2}}e^{2b(k-1)t}+\frac{\sigma^{2}L^{2}}{4b(1-k)}[1-e^{2b(k-1)t}];\\ 
    \text{Var}(X_{t})&=\frac{\sigma^{2}L^{4}\big(\lambda_{0}\big)^{2-2k}}{8b(1-k)^3}\big[e^{2b(k-1)t}-e^{4b(k-1)t}\big]+\frac{\sigma^{4}L^{4}}{16b^{2}(1-k)^{2}}\big[1-e^{2b(k-1)t}\big]^{2};\nonumber\\
    \text{Cov}(X_{t},X_{t'})&=\frac{\sigma^{2}L^{4}\big(\lambda_{0}\big)^{2-2k}}{8b(1-k)^3}\big[e^{2b(k-1)t'}-e^{2b(k-1)(t+t')}\big]+\frac{\sigma^{4}L^{4}}{16b^{2}(1-k)^{2}}\big[e^{2b(1-k)(t-t')}+e^{2b(k-1)(t+t')}-2e^{2b(k-1)t'}\big];\nonumber
\end{align}
\begin{equation}\nonumber\hspace{-0.6cm}
    \mathbb{E}[(X_{t})^{n}]=\sum_{j=0}^{n}\binom{n}{j}\frac{\Gamma(n+1)}{\Gamma(j+1)}s_{t}^{n-j}h_{t}^{j},\quad\text{where }s_{t}=\frac{\sigma^{2}L^{2}}{4b(1-k)}\big[1-e^{2b(k-1)t}\big],~h_{t}=\frac{L^{2}\lambda_{0}^{2-2k}}{4(1-k)^{2}}
  e^{2b(k-1)t}.
\end{equation}
Using the relation $V_{t}=\frac{4(1-k)^{2}}{L^{2}}X_{t}$ and the above moment formulas for $X_{t}$, we obtain:
\begin{align}\nonumber
	\mathbb{E}\big[V_{t}\big]&=(\lambda_{0})^{2-2k}e^{2b(k-1)t}+\frac{\sigma^{2}(1-k)}{b}\big[1-e^{2b(k-1)t}\big];\\
	\text{Var}\big(V_{t}\big)&=\frac{2\sigma^{2}{(}1-k)}{b}(\lambda_{0})^{2-2k}\big[e^{2b(k-1)t}-e^{4b(k-1)t}\big]+\frac{\sigma^{4}(1-k)^{2}}{b^{2}}\big[1-e^{2b(k-1)t}\big]^{2};\nonumber\\
	\text{Cov}\big(V_{t},V_{t'}\big)&=\frac{2\sigma^{2}(1-k)}{b}(\lambda_{0})^{2-2k}\big[e^{2b(k-1)t'}-e^{2b(k-1)(t+t')}\big]+\frac{\sigma^{4}(1-k)^{2}}{b^{2}} \Big[e^{2b(1-k)(t-t')}+e^{2b(k-1)(t+t')}-2e^{2b(k-1)t'}\Big];\nonumber\\
	\mathbb{E}\big[\big(V_{t}\big)^{n}\big]&=\Big[\frac{4(1-k)^{2}}{L^{2}}\Big]^{n}\mathbb{E}[(X_{t})^{n}].\nonumber
\end{align}
Here $\iota$ and $s_{t}$ are evaluated at the CIR parameters
of $X$, and $h_{t}=X_{0}e^{-2b(1-k)t}$.

\noindent \textbf{Lemma 3.7}\\
Given the value of $a^{*}$, $b^{*}$ and $\sigma^{*}$ and thus the value of $\kappa=\frac{2a^{*}b^{*}}{\sigma^{*2}}-1$, the asymptotic stationary probability density function of $r_{t}$ with $t$ going to infinity, ranging over $[0,+\infty)$, is of the gamma type (parameters $\kappa+1$ and $\frac{\kappa+1}{b^{*}}$), which means:
\begin{equation}\nonumber
    p_{\infty}(x)=f(x\lvert a^{*},b^{*},\sigma^{*})=\frac{(\frac{\kappa+1}{b^{*}})^{\kappa+1}}{\Gamma(\kappa+1)}x^{\kappa}\text{exp}\Big\{-\frac{\kappa+1}{b^{*}}x\Big\},~x\in(0,+\infty).
\end{equation}

\noindent \textit{Proof}: See the Appendix.\qed\\

\noindent In our case:
\begin{align}\nonumber
    &p_{\infty}(x)=\frac{\big[\frac{4b(1-k)}{\sigma^{2}L^{2}}\big]^{1}}{\Gamma(1)}x^{0}\text{exp}\Big\{\frac{4b(k-1)}{\sigma^{2}L^{2}}x\Big\}
    =\frac{4b(1-k)}{\sigma^{2}L^{2}}\text{exp}\Big\{\frac{4b(k-1)}{\sigma^{2}L^{2}}x\Big\},\quad x\in(0,+\infty).
\end{align}
With $v=\frac{4(1-k)^{2}}{L^{2}}x$ where $v$ denotes the value taken by $V_{t}$. By a linear change of variables, the stationary density of $V_{t}$ is:
\begin{equation}\nonumber
	p_{\infty}^{V}(v)=\frac{L^{2}}{4(1-k)^{2}}p_{\infty}\Big(\frac{L^{2}}{4(1-k)^{2}}v\Big)=\frac{b}{\sigma^{2}(1-k)}\text{exp}\Big\{-\frac{b}{\sigma^{2}(1-k)}v\Big\},\quad v\in(0,+\infty),
\end{equation}
With $x=T(\ell)=\frac{L^{2}}{4(1-k)^{2}}\ell^{2(1-k)}$ where $\ell$ denotes the value taken by $\lambda_{t}$, we deduce $\frac{dx}{d\ell}=\frac{L^{2}}{2(1-k)}\ell^{1-2k}$. By change of variables, the stationary density of $\lambda_{t}$ is
\begin{equation}\nonumber
	p_{\infty}^{\lambda}(\ell)=p_{\infty}(x)\Big\lvert\frac{dx}{d\ell}\Big\lvert=\frac{2b}{\sigma^{2}}\ell^{1-2k}
  \text{exp}\Big\{-\frac{b}{\sigma^{2}(1-k)}\ell^{2(1-k)}\Big\},\quad \ell\in(0,+\infty).
\end{equation}
Thus, under $\mathbb{Q}$, the stationary distribution of $V$ is an exponential distribution with rate parameter $b/[\sigma^2(1-k)]$, whereas the stationary distribution of $\lambda$ is a Weibull distribution with shape parameter
$2(1-k)$ and scale parameter $[\sigma^2(1-k)/b]^{1/[2(1-k)]}$.

\section{Martingale property of \texorpdfstring{$M_{t}$}{Mt}}
\subsection{Limitations of moment estimates in verifying Novikov's condition}
\noindent Let $\lambda_{t}$ be the CKLS process
\begin{equation}\nonumber
	d\lambda_{t}=(a-b\lambda_{t})dt+\sigma(\lambda_{t})^{k}dW_{t},~\lambda_{0}>0,~k\in\big[\frac{1}{2},1\big)~(\text{When }k=\frac{1}{2}, 2a\geq\sigma^{2}).
\end{equation} 
Recall in our case the Dol\'{e}ans-Dade exponential (2.8) has the following expression
\begin{equation}\nonumber
	M_{t}\overset{\text{def}}{=}\text{exp}\Big\{\int_{0}^{t} q_{s}dW_{s}-\frac{1}{2}\int_{0}^{t}(q_{s})^{2}ds \Big\},
\end{equation}
with its kernel:
\begin{equation}\nonumber
	q_{t}=\frac{\sigma}{2}(\lambda_{t})^{k-1}-\frac{a}{\sigma}(\lambda_{t})^{-k}.
\end{equation}
Our goal is to show that $M_{t}$ is a true martingale on every finite horizon. One sufficient condition would be Novikov's condition:
\begin{equation}\nonumber
	\mathbb{E}^{\mathbb{P}}\Big[\text{exp}\Big\{\frac{1}{2}\int_{0}^{T} (q_{s})^{2}ds\Big\}\Big]<+\infty ~\text{holds }\forall T>0,
\end{equation}
under the parameter constraint $\frac{1}{2}\leq k<1$ (when $k=\frac{1}{2}$, $2a\geq \sigma^{2}$). A direct computation gives the expression of the square of $q_{t}$:
\begin{equation}\nonumber
	(q_{t})^{2}=\frac{\sigma^{2}}{4}(\lambda_{t})^{2k-2}-a(\lambda_{t})^{-1}+\frac{a^{2}}{\sigma^{2}}(\lambda_{t})^{-2k}.
\end{equation}
To illustrate the limitation of moment estimates in verifying Novikov's condition, we consider the case $\frac{1}{2}<k<1$. Recall that, by \textbf{Theorem 1.2 (2) and (6)}, the CKLS process $\lambda_{t}$ remains strictly positive almost surely and satisfies, for every $p\geq 0$ and every finite $T>0$,
\begin{equation}\nonumber
	\mathbb{E}^{\mathbb{P}}\Big[\sup_{0\leq t\leq T}(\lambda_{t})^{-p}\Big]<+\infty.
\end{equation}
In other words, $\lambda_{t}$ never hits $0$ and possesses finite negative-power moments of every non-zero order. Because $\frac{1}{2}<k<1$, we have $-1<2k-2<0$. Setting $-p=2k-2$ yields
\begin{equation}\nonumber
	\mathbb{E}^{\mathbb{P}}\Big[\int_{0}^{T}(\lambda_{t})^{2k-2}dt\Big]<+\infty.
\end{equation}
Analogously, choosing $p=1$ and $p=2k>0$ yields the same conclusions:
\begin{equation}\nonumber
	\mathbb{E}^{\mathbb{P}}\Big[\int_{0}^{T}(\lambda_{t})^{-1}dt\Big]<+\infty,\qquad \mathbb{E}^{\mathbb{P}}\Big[\int_{0}^{T}(\lambda_{t})^{-2k}dt\Big]<+\infty.\nonumber
\end{equation}
Combining all three, together with the linearity of expectation, gives
\begin{equation}\nonumber
	\mathbb{E}^{\mathbb{P}}\Big[\int_{0}^{T}(q_{t})^{2}dt\Big]<+\infty.
\end{equation}
Suppose that, for each finite $T>0$, we can verify that the following:
\begin{equation}\nonumber
      \exists \text{ some constant }C_{T}\geq0\text{ such that }\int_{0}^{T}(q_{s})^{2}\mathrm ds \leq C_{T}\text{ holds }a.s..
\end{equation}
Then by the fact that the exponential function $e^{x}$ is monotone increasing, we have the trivial implication:
\begin{equation}\nonumber
      \int_{0}^{T}(q_{t})^{2}dt \leq C_{T}<+\infty
      \Longrightarrow
      \text{exp}\Big\{\frac{1}{2}\int_{0}^{T}(q_{t})^{2}dt\Big\}
      \leq e^{C_{T}/2}
      \Longrightarrow
      \mathbb{E}^{\mathbb{P}}\Big[\text{exp}\Big\{\frac{1}{2}\int_{0}^{T}(q_{t})^{2}dt\Big\}\Big]
      \leq e^{C_{T}/2}<+\infty,
\end{equation}
which verifies Novikov's condition. Novikov's condition also implies Kazamaki's condition, as shown in the Appendix. However, the moment estimates above \textbf{only} give $\mathbb{E}^{\mathbb{P}}\big[\int_{0}^{T}(q_{s})^{2}ds\big]<+\infty$. Since the integral is nonnegative, this implies $\int_{0}^{T}(q_{s})^{2}ds<+\infty$ $\mathbb{P}-a.s.$ Nevertheless, this estimate \textbf{does not} by itself establish
\[
    \mathbb{E}^{\mathbb{P}}\Big[\text{exp}\Big\{\frac{1}{2}\int_{0}^{T}(q_{s})^{2}ds\Big\}\Big]<+\infty.
\]
which is the exponential integrability required by Novikov's condition.\\

\noindent We give a simple example to illustrate this limitation. Let $Z$ have an exponential distribution with rate $\frac{1}{2}$, namely, $f_{Z}(z)=\frac{1}{2}e^{-z/2}$ for $z\geq0$. Then $\mathbb{E}[Z]=2<+\infty$, whereas $\mathbb{E}[e^{Z/2}]=\int_{0}^{\infty}e^{z/2}\frac{1}{2}e^{-z/2}dz=+\infty$. In fact, $\mathbb{E}[Z^{n}]=2^{n}n!<+\infty$ for every positive integer $n$. Thus, even finiteness of all positive integer moments does not, by itself, guarantee the required exponential integrability.

\subsection{A concise outline of the proof strategy}
\noindent From now on, we focus on proving that the Dol\'{e}ans–Dade exponential $M_{t}$ in (2.8) is a true martingale on every finite time interval, without relying on Novikov's or Kazamaki's condition.\\

\begin{mdframed}[
  linecolor=black,linewidth=2pt,
  roundcorner=5pt,
  backgroundcolor=white
]
\noindent \textbf{Overall Objective}\\
We want to verify that our Dol\'{e}ans-Dade exponential (2.8) of the form
\begin{equation}\nonumber
	\mathbb{M}_{t}\overset{\text{def}}{=}\text{exp}\Big\{\int_{0}^{t\wedge\tau}q(\lambda_{u})dW_{u}-\frac{1}{2}\int_{0}^{t\wedge\tau}[q(\lambda_{u})]^{2}du\Big\},~t\in[0,+\infty),
\end{equation}
is indeed a true martingale. We will confirm the correctness of this by applying \textbf{Theorem 4.1}, as originally established by \hyperlink{Mijatovic2012}{Mijatović and Urusov (2012)} based on Feller's test for explosion proposed by \hyperlink{Feller1952}{Feller (1952)}.\\
\noindent \textbf{A Step-by-step Verification}\\
\begin{enumerate}
	\item \textbf{To Verify Assumptions of Theorem 4.1}\\
		\emph{Assumption 1.} With $J=(0,+\infty)$, the CKLS model's drift $\mu(z)=a-bz$ and diffusion $\nu(z)=\sigma z^{k}$ satisfy:
			\begin{equation}\nonumber
				\text{Both }\mu(z)[\nu(z)]^{-2}\text{ and }[\nu(z)]^{-2}\text{ belong to the class of locally integrable functions } \mathcal{L}^1_{\mathrm{loc}}\big(J\big).
			\end{equation}
	    \emph{Assumption 2.} The Dol\'{e}ans-Dade exponential's kernel $q(z)$ satisfies
			\begin{equation}\nonumber
				[q(z)]^{2}[\nu(z)]^{-2}\in \mathcal{L}^1_{\mathrm{loc}}\big(J\big).
			\end{equation}
	\item \textbf{To Define an Auxiliary Process}
		\begin{enumerate}
			\item \emph{Construction of the auxiliary process}
            Define an auxiliary diffusion $\dot{\lambda}_{t}$
            on the same state space $J$ of $\lambda_{t}$ (which is $(0,+\infty)$), whose drift and diffusion coefficients incorporate both the original CKLS parameters $\mu(z)$, $\nu(z)$ and Dol\'{e}ans-Dade exponential's kernel $q(z)$:
            \begin{equation}\nonumber
                   d\dot{\lambda}_{t}=\underbrace{(\mu+q\nu)(\dot{\lambda}_{t})}_{=:~\gamma(\dot{\lambda}_{t})}dt+\nu(\dot{\lambda}_{t})d\dot{W}_{t},~\dot{\lambda}_{0}=\lambda_{0}\in J.
            \end{equation}
		\item \emph{Key Criterion established by Mijatovi\'{c} and Urusov (2012)} (\textbf{Theorem 4.1}) 
			\begin{equation}\nonumber
				\mathbb{M}_{t}\text{ is a (true) martingale}\quad\Longleftrightarrow\quad \dot{\lambda}_{t}\text{ does not exit }J.
			\end{equation}
        \end{enumerate}    
        \item \textbf{To Carry out the Test for Explosion via Theorem 4.2}
		\begin{enumerate}
		\item \emph{Feller's Tests for Inaccessibility of Endpoints} (\textbf{Theorem 4.2})\\
            In order to verify that $\dot{\lambda}_{t}$ truly does not exit $J$, we need the help of Feller's test for explosion.
			\begin{itemize}
				\item \textbf{Sufficient condition (the simple test, easier to calculate).}
				Compute the testing function for explosion (for some $c\in J$):
            			\begin{equation}\nonumber
					\psi(x)\overset{\text{def}}{=}\int_{c}^{x}\text{exp}\Big\{-2\int_{c}^{y}\frac{\gamma(z) dz}{[\nu(z)]^{2}}\Big\}dy,~x\in J,
				\end{equation}
		          and check if in our case both $\lim_{x\to +\infty}\psi(x)=+\infty$ and $\lim_{x\to 0+}\psi(x)=-\infty$ hold.
			\end{itemize}
		\item \emph{Application to CKLS Auxiliary Process.}
			\begin{itemize}
				\item \textbf{As }$x\to+\infty$: The simple test suffices to show inaccessibility via $\psi(x)\to+\infty$.
				\item \textbf{As }$x\to0+$: The simple test also suffices to show inaccessibility, since $\psi(x)\to-\infty$.
			\end{itemize}
		\end{enumerate}
\item \textbf{To Conclude}\\
Under Assumptions 1 and 2, and having shown that $\dot{\lambda}_{t}$ cannot exit its state space at either endpoint, the Dol\'{e}ans-Dade exponential $M_{t}$ is a true martingale by \textbf{Theorem 4.1} established by \hyperlink{Mijatovic2012}{Mijatovi\'{c} and Urusov (2012)} based on Feller's test for explosion proposed by \hyperlink{Feller1952}{Feller (1952)}.
\end{enumerate}
\end{mdframed}

\subsection{Some settings and the martingale theorem}
\noindent Consider the state space $J\overset{\text{def}}{=}(l,r)\subset\mathbb{R}\cup\{\pm\infty\}$, $-\infty\leq l<r \leq+\infty$ and a $J$-valued diffusion process $Z_{t}$, $t\in[0,T]$, defined on some filtered probability space $(\Omega,\mathcal{F},\{\mathcal{F}_{t}\}_{t\in[0,T]},\mathbb{P})$ admitting the following dynamics:
\begin{equation}
    dZ_{t}=\mu(Z_{t})dt+\nu(Z_{t})dW_{t},~Z_{0}=z_{0}\in J,
\end{equation}
where $\mu(Z_{t})$ and $\nu(Z_{t})$ are both $J\xrightarrow{}\mathbb{R}$ Borel functions.\\

\noindent \textbf{Assumption 1} (Engelbert-Schmidt conditions)
\begin{align}\nonumber
    \nu(z)\neq&0, \forall z\in J;\\
    \nu^{-2}\in \mathcal{L}^{1}_{\text{loc}}(J)&;~\mu\nu^{-2}\in \mathcal{L}^{1}_{\text{loc}}(J),\nonumber
\end{align}
where $\mathcal{L}^{1}_{\text{loc}}(J)$ denotes the function class of local integrability, which means that the functions inside of the class are $J\to\mathbb{R}$ integrable on compact subsets of $J$.$\hfill\blacksquare$\\

\noindent According to the original papers \hyperlink{Engelbert1984}{Engelbert and Schmidt (1984)} and \hyperlink{Engelbert1991}{Engelbert and Schmidt (1991)} (see also \hyperlink{Karatzas2012}{Karatzas and Shreve (2012)} [Chapter 5, Theorem 5.15, page 341; Theorem 5.7 pages 335-336]), (4.1) has a unique in-law-weak solution that has the chance to exit its state space. Denote the possible exit time by $\xi$ as:
\begin{equation}\nonumber
    \xi\overset{\text{def}}{=}\text{inf}\{t\geq 0\lvert Z_{t}\notin J\}.
\end{equation}
The diffusion $Z$ never exits $J$ almost surely if $\mathbb{P}(\xi=+\infty)=1$. If $\mathbb{P}(\xi<+\infty)>0$, it exits $J$ with positive probability. We further assume that over the set $\{\xi<+\infty\}$, the solution $Z$ remains at the boundary point of $J$ where it exits, after time $\xi$, i.e. the left and right boundaries, denoted by $l$ and $r$ respectively, are the so-called $\verb+"+$absorbing boundaries$\verb+"+$. This convention specifies the behavior after exit and does not alter the exit time from the open interval $J$. We therefore give the following definitions:\\

\noindent \textbf{Definitions} (Exits of the state space)\\
\noindent (1) $Z$ exits $J$ at $r$ means: $\mathbb{P}(\xi<+\infty, \lim_{t\uparrow \xi}Z_{t}=r)>0$.\\
\noindent  (2) $Z$ exits $J$ at $l$ means: $\mathbb{P}(\xi<+\infty, \lim_{t\uparrow\xi}Z_{t}=l)>0$.\\

\noindent \textbf{Assumption 2}\\
\noindent Suppose $q$ is a Borel measurable function $J\xrightarrow{}\mathbb{R}$ satisfying:
\begin{equation}
    q^{2}\nu^{-2}\in\mathcal{L}^{1}_{\text{loc}}(J)
\end{equation}

\noindent Now of our interest is whether the following stochastic process is a (true) martingale:
\begin{equation}
    \mathbb{M}_{t}=\text{exp}\Big\{\int_{0}^{t\wedge\xi} q(Z_{u})dW_{u}-\frac{1}{2}\int_{0}^{t\wedge\xi}[q(Z_{u})]^{2}du\Big\},~t\in[0,T],
\end{equation}
and we set $\mathbb{M}_{t}=0$ for $t\geq\xi$ on $\{\xi<+\infty,\int_{0}^{\xi}[q(Z_{u})]^{2}du=+\infty\}$. A stochastic process of this form is just what we have in (2.8). $\hfill\blacksquare$\\

\noindent Consider an auxiliary $J$-valued diffusion process $\dot{Z}_{t}$, $t\in[0,T]$, defined on $(\dot{\Omega},\dot{\mathcal{F}},(\dot{\mathcal{F}}_{t})_{t\in[0,T]},\dot{\mathbb{P}})$ admitting the following dynamics:
\begin{equation}
    d\dot{Z}_{t}=(\mu+q\nu)(\dot{Z}_{t})dt+\nu(\dot{Z}_{t})d\dot{W}_{t},~\dot{Z}_{0}=z_{0}\in J.
\end{equation}
Here $\dot{W}$ is a standard Brownian motion under $\dot{\mathbb{P}}$ with respect to the filtration $\dot{\mathcal{F}}_{t}$, and the auxiliary diffusion starts from the same initial value as $Z$.\\

\noindent Suppose that \textbf{Assumption 1} for $\mu$ and $\nu$ in (4.1) is satisfied. As long as \textbf{Assumption 2} for $q$ and $\nu$ in (4.2) is satisfied too, we can immediately conclude that the auxiliary SDE (4.4) meets the Engelbert-Schmidt criteria for both the drift coefficient $\mu+q\nu$ and the diffusion function $\nu$. This is because checking if $(\mu+q\nu)\nu^{-2}=\mu\nu^{-2}+q\nu^{-1}\in\mathcal{L}^{1}_{\text{loc}}(J)$ holds is equivalent to checking if $\mu\nu^{-2}\in\mathcal{L}^{1}_{\text{loc}}(J)$ and $q{\nu}^{-1}\in\mathcal{L}^{1}_{\text{loc}}(J)$ hold (the correctness of which can be directly induced by (4.2)). Therefore, it is ensured that equation (4.4) possesses a unique weak solution in law that can exit its defined state space, which is $J$ in (4.1).\\

\noindent Now we state the most important tool being used in this paper, which is a sufficient and necessary condition for judging whether a stochastic process of the form (4.3) is a (true) martingale or not under a certain probability measure:\\

\noindent \textbf{Theorem 4.1} (A sufficient and necessary condition for (4.3) to be a (true) martingale)\\
Suppose that \textbf{Assumption 1} and \textbf{Assumption 2} are satisfied and that $Z$ in (4.1) does not exit $J$ under $\mathbb{P}$. Then $\mathbb{M}_{t}$ in (4.3) is a (true) $\mathbb{P}$-martingale on every finite time interval if and only if the auxiliary diffusion $\dot{Z}$ in (4.4), started from $z_{0}$, does not exit $J$ under $\dot{\mathbb{P}}$.\\

\noindent \textit{Proof}: See \hyperlink{Mijatovic2012}{Mijatovi\'c and Urusov (2012)} [Corollary 2.2]. One may also refer to \hyperlink{Lewis2016}{Lewis (2016)} Chapter 6 \S 12, pages 340-345, or \hyperlink{Karatzas2016}{Karatzas and Ruf (2016)} [Section 3.1, Theorem 3.2, pages 1031-1034; Chapter 5 \S 3, Theorem 5.9, pages 1050-1051]. The detailed proof falls outside the scope and depth of this paper and will therefore not be provided here.\qed

\subsection{Verification that \texorpdfstring{$M_{t}$}{Mt} in our case is a (true) martingale}
\noindent Now we go back to our measure-transformed CKLS model again. Recall that in (2.8) we have $M_{t}=\text{exp}\{\int_{0}^{t}q_{s}dW_{s}-\frac{1}{2}\int_{0}^{t}(q_{s})^{2}ds\}$ being of the same form as stated in (4.3). We introduce the state space: $S=(l,r)$ with $-\infty\leq l<r \leq+\infty$, and introduce the time of $\lambda_{t}$ exiting $S$: $\tau\in\mathbb{R}\cup\{+\infty\}$. We now rewrite the Dol\'{e}ans-Dade exponential in (2.8), $M_{t}=\text{exp}\{\int_{0}^{t}q_{s}dW_{s}-\frac{1}{2}\int_{0}^{t}(q_{s})^{2}ds\}$, as:
\begin{align}\nonumber
    M_{t}=\text{exp}\Big\{\int_{0}^{t\wedge\tau}q(\lambda_{u})&dW_{u}-\frac{1}{2}\int_{0}^{t\wedge\tau}[q(\lambda_{u})]^{2}du\Big\},~t\in[0,+\infty),\\
    q(\lambda_{t})=\frac{\sigma}{2}&(\lambda_{t})^{k-1}-\frac{a}{\sigma}(\lambda_{t})^{-k},~t\in[0,+\infty).\nonumber
\end{align}
Our aim is to prove that $M_{t}$ is a (true) martingale with the help of \textbf{Theorem 4.1}.\\

\noindent Firstly, we need to verify that \textbf{Assumption 1} and \textbf{Assumption 2} are satisfied. Recall (1.1) is the original interest of us: 
\begin{equation}\nonumber
    d\lambda_{t}=(a-b\lambda_{t})dt+\sigma (\lambda_{t})^{k}dW_{t},~\lambda_{t}\big\lvert_{t=0}=\lambda_{0}>0.
\end{equation}
So in our case what we need to verify are:
\begin{align}
    \lambda_{0}\in S;&\tag{\RomanNumeralCaps{1}}\\
    \sigma y^{k}\neq0,~\forall y\in S;&\tag{\RomanNumeralCaps{2}}\\
    y\mapsto\frac{1}{(\sigma y^{k})^{2}}\in\mathcal{L}^{1}_{\text{loc}}(S),~\forall y\in S;&\tag{\RomanNumeralCaps{3}.\romannumeral1}\\
    y\mapsto\frac{a-by}{(\sigma y^{k})^{2}}\in\mathcal{L}^{1}_{\text{loc}}(S),~\forall y\in S;&\tag{\RomanNumeralCaps{3}.\romannumeral2}\\
    y\mapsto\frac{[q(y)]^{2}}{(\sigma y^{k})^{2}}\in\mathcal{L}^{1}_{\text{loc}}(S),~\forall y\in S.&\tag{\RomanNumeralCaps{4}}
\end{align}

\noindent Secondly, we must know what $S=(l,r)$ is like in our case before starting the verification. Recall that the parameter restrictions in (2.10) for the CIR transformation are imposed throughout the following verification. Thus, from now on we always assume:
\begin{align}\nonumber
    \frac{1}{2}\leq k<1 \text{ (when }k=\frac{1}{2},~2a\geq \sigma^{2}),~L>0~,a~>0,~b>0,~\sigma>0.
\end{align}
Recall in \textbf{Theorem 1.2}, we already know from the statements (1)-(4) what values $\lambda_{t}$ can take for different cases. In our current setting, either for the case $k\in(\frac{1}{2},1)$, or for the case $k=\frac{1}{2}$ with $2a\geq \sigma^{2}$, $\lambda_{t}\in(0,+\infty)$.\\

\noindent Therefore, in our case, we can say that:\\
Case \ding{172} For $k\in(\frac{1}{2},1)$, $S=(0,+\infty)$;\\
Case \ding{173} For $k=\frac{1}{2}$ with $2a\geq \sigma^{2}$, $S=(0,+\infty)$.\\

\noindent Obviously, \RomanNumeralCaps{1} and \RomanNumeralCaps{2} are satisfied for both Case \ding{172} and Case \ding{173}. We now check \RomanNumeralCaps{3}.\romannumeral1, \RomanNumeralCaps{3}.\romannumeral2~and \RomanNumeralCaps{4} for Case \ding{172}.
\begin{align}
    \int_{c}^{c+\epsilon}\frac{1}{\sigma^{2}y^{2k}}dy<+\infty,~\forall &c\in(0,+\infty),~\epsilon>0;\\
    \int_{c}^{c+\epsilon}\frac{\lvert a-by\lvert}{\sigma^{2}y^{2k}}dy<+\infty,~\forall &c\in(0,+\infty),~\epsilon>0;\\
    \int_{c}^{c+\epsilon}\Big(\frac{1}{4}y^{-2}+\frac{a^{2}}{\sigma^{4}}y^{-4k}-\frac{a}{\sigma^{2}}y^{-1-2k}\Big)dy&<+\infty,~\forall c\in(0,+\infty),~\epsilon>0.
\end{align}
For (4.5), since $k\neq \frac{1}{2}$:
\begin{align}\nonumber
    \int_{c}^{c+\epsilon}\frac{1}{\sigma^{2}y^{2k}}dy=\frac{1}{\sigma^{2}}\int_{c}^{c+\epsilon}y^{-2k}dy=\frac{1}{\sigma^{2}}\Big[\frac{1}{1-2k}y^{1-2k}\Big]^{c+\epsilon}_{c}=\frac{1}{\sigma^{2}}\Big[\frac{1}{1-2k}(c+\epsilon)^{1-2k}-\frac{1}{1-2k}c^{1-2k}\Big]<+\infty,~\forall c\in(0,+\infty),~\epsilon>0.\nonumber
\end{align}
For (4.6), since $k\neq\frac{1}{2},1$ and $a,b>0$, we have:
\begin{align}\nonumber
    &\int_{c}^{c+\epsilon}\frac{\lvert a-by\lvert}{\sigma^{2}y^{2k}}dy\leq\frac{1}{\sigma^{2}}\int_{c}^{c+\epsilon}\Big(ay^{-2k}+by^{1-2k}\Big)dy=\frac{1}{\sigma^{2}}\Big[\frac{a}{1-2k}y^{1-2k}\Big]_{c}^{c+\epsilon}+\frac{1}{\sigma^{2}}\Big[\frac{b}{2-2k}y^{2-2k}\Big]_c^{c+\epsilon}\\
    &=\frac{1}{\sigma^{2}}\Big[\frac{a}{1-2k}(c+\epsilon)^{1-2k}-\frac{a}{1-2k}c^{1-2k}+\frac{b}{2-2k}(c+\epsilon)^{2-2k}-\frac{b}{2-2k}c^{2-2k}\Big]<+\infty,~\forall c>0,~\epsilon>0.\nonumber
\end{align}
For (4.7), since $k\neq 0, \frac{1}{4}$:
\begin{align}\nonumber
    &\int_{c}^{c+\epsilon}\Big(\frac{1}{4}y^{-2}+\frac{a^{2}}{\sigma^{4}}y^{-4k}-\frac{a}{\sigma^{2}}y^{-1-2k}\Big)dy=\Big[-\frac{1}{4}y^{-1}\Big]^{c+\epsilon}_{c}+\Big[\frac{a^{2}}{\sigma^{4}}\frac{1}{1-4k}y^{1-4k}\Big]^{c+\epsilon}_{c}-\Big[\frac{a}{\sigma^{2}}\frac{1}{-2k}y^{-2k}\Big]^{c+\epsilon}_{c}\nonumber\\
    &=\frac{\epsilon}{4c(c+\epsilon)}+\frac{a^{2}\big((c+\epsilon)^{1-4k}-c^{1-4k}\big)}{\sigma^{4}(1-4k)}-\frac{a\big(c^{-2k}-(c+\epsilon)^{-2k}\big)}{2k\sigma^{2}}<+\infty,~\forall c\in(0,+\infty),~\epsilon>0.\nonumber
\end{align}

\noindent We now check \RomanNumeralCaps{3}.\romannumeral1, \RomanNumeralCaps{3}.\romannumeral2~and \RomanNumeralCaps{4} for Case \ding{173}. In this case, $k=\frac{1}{2}$ and $c\in(0,+\infty)$, $y\in(c,c+\epsilon]$. We rewrite (4.5), (4.6), (4.7), which are the assumptions to be checked, as:
\begin{align}
    \int_{c}^{c+\epsilon}\frac{1}{\sigma^{2}y}dy&<+\infty,~\forall c\in(0,+\infty);\nonumber\\
    \int_{c}^{c+\epsilon}\frac{\lvert a-by\lvert}{\sigma^{2}y}dy&<+\infty,~\forall c\in(0,+\infty);\nonumber\\
    \int_{c}^{c+\epsilon}\Big(\frac{1}{4}+\frac{a^{2}}{\sigma^{4}}-\frac{a}{\sigma^{2}}\Big)&y^{-2}dy<+\infty,~\forall c\in(0,+\infty).\nonumber
\end{align}
It is easy to verify that:
\begin{align}
    \int_{c}^{c+\epsilon}&\frac{1}{\sigma^{2}y}dy=\frac{1}{\sigma^{2}}\Big[\text{log}y\Big]^{c+\epsilon}_{c}=\frac{1}{\sigma^{2}}\text{log}\frac{c+\epsilon}{c}<+\infty,~\forall c\in(0,+\infty),~\epsilon>0;\nonumber\\
    \int_{c}^{c+\epsilon}\frac{\lvert a-by\lvert}{\sigma^{2}y}dy&\leq\frac{1}{\sigma^{2}}\Big[a\text{log}y\Big]^{c+\epsilon}_{c}+\frac{1}{\sigma^{2}}\Big[by\Big]^{c+\epsilon}_{c}=\frac{1}{\sigma^{2}}(a\text{log}\frac{c+\epsilon}{c}+b\epsilon)<+\infty,~\forall c\in(0,+\infty),~\epsilon>0;\nonumber\\
    \int_{c}^{c+\epsilon}\Big(\frac{1}{4}+\frac{a^{2}}{\sigma^{4}}-\frac{a}{\sigma^{2}}\Big)y^{-2}dy&=\Big[-\Big(\frac{1}{4}+\frac{a^{2}}{\sigma^{4}}-\frac{a}{\sigma^{2}}\Big)y^{-1}\Big]_{c}^{c+\epsilon}
    =\Big(\frac{a}{\sigma^{2}}-\frac{1}{4}-\frac{a^{2}}{\sigma^{4}}\Big)\Big(\frac{1}{c+\epsilon}-\frac{1}{c}\Big)<+\infty,~\forall c\in(0,+\infty),~\epsilon>0.\nonumber
\end{align}
We thus conclude that for the Case \ding{172} and the Case \ding{173}, \textbf{Assumption 1} and \textbf{Assumption 2} are satisfied.\\

\noindent According to \textbf{Theorem 4.1}, we now consider the auxiliary $J$-valued diffusion process of the form (4.4):
\begin{equation}
    d\dot{\lambda}_{t}=\underbrace{(\mu+q\nu)(\dot{\lambda}_{t})}_{=:\gamma(\dot{\lambda}_{t})}dt+\nu(\dot{\lambda}_{t})d\dot{W}_{t},~\dot{\lambda}_{t}\big\lvert_{t=0}=\lambda_{0},
\end{equation}
where in our case:
\begin{equation}
\left\{
\begin{aligned}
    \mu(\dot{\lambda}_{t})&=a-b\dot{\lambda}_{t}\\ 
    \nu(\dot{\lambda}_{t})&=\sigma\big(\dot{\lambda}_{t}\big)^{k}\\
    q(\dot{\lambda}_{t})&=\frac{\sigma}{2}\big(\dot{\lambda}_{t}\big)^{k-1}-\frac{a}{\sigma}\big(\dot{\lambda}_{t}\big)^{-k}
\end{aligned}
\right.
\end{equation}
Thus, given $\dot{\lambda}_{t}\big\lvert_{t=0}=\lambda_{0}\in S=(0,+\infty)$, (4.8) will admit the following expression:
\begin{align}\nonumber
    d\dot{\lambda}_{t}&=\Big[a-b\dot{\lambda}_{t}+\big(\frac{\sigma}{2}\big(\dot{\lambda}_{t}\big)^{k-1}-\frac{a}{\sigma}\big(\dot{\lambda}_{t}\big)^{-k}\big)\sigma\big(\dot{\lambda}_{t}\big)^{k} \Big]dt+\sigma\big(\dot{\lambda}_{t}\big)^{k}d\dot{W}_{t}\\
    &=\Big[\frac{\sigma^{2}}{2}\big(\dot{\lambda}_{t}\big)^{2k-1}-b\dot{\lambda}_{t}\Big]dt+\sigma\big(\dot{\lambda}_{t}\big)^{k}d\dot{W}_{t}.\nonumber
\end{align}

\noindent Now we state the second most important theorem used in this paper, which is used to test whether the solution to the $J$-auxiliary SDE of $\dot{Z}$ exits the state space $J$ of $Z$ or not. For this reason, this theorem can be regarded as a natural continuation/sequel of \textbf{Theorem 4.1}.\\

\noindent \textbf{Theorem 4.2} (Feller's test for explosions)\footnote{Actually, the theorem is intended to be applied to SDE (4.8) in our case. However, to avoid cumbersome notation, we state the assumptions, lemmas, and the theorem itself using the notation from (4.1) throughout this section, without loss of generality. The correspondence to (4.8) will be clarified when we explain how the theorem is applied in practice.}\\
Assume \textbf{Assumption 1} (Engelbert-Schmidt conditions) holds, so that there is an in-law-weak solution $Z_{t}$ to (4.1) or (4.8) (Note that in the case of (4.8), one should substitute $\mu$ appearing in the following expressions with $\gamma$) existing in $J$. Given a non-random initial condition $Z_{0}=z_{0}\in J$. For some $c\in J$, define the scale function for testing:
\begin{equation}\nonumber
    \psi(x)\overset{\text{def}}{=}\int_{c}^{x}\text{exp}\Big\{-2\int_{c}^{y}\frac{\mu(z) dz}{[\nu(z)]^{2}}\Big\}dy,~x\in J,
\end{equation}
and define the test function for Feller's test:
\begin{equation}\nonumber
    \phi(x)\overset{\text{def}}{=}\int_{c}^{x}\psi'(y)\int_{c}^{y}\frac{2dz}{\psi'(z)[\nu(z)]^{2}}dy,~x\in J,
\end{equation}
then $\mathbb{P}(\xi=+\infty)=1$ or $\mathbb{P}(\xi<+\infty)>0$ according to whether $\lim_{x\uparrow r}\phi(x)=\lim_{x\downarrow l}\phi(x)=+\infty$ or not.\\

\noindent \textit{Proof}: See the original paper by \hyperlink{Feller1952}{Feller (1952)} [Sections 20-23, pages 507-519] or \hyperlink{Karatzas2016}{Karatzas and Shreve (2012)} [Chapter 5 \S 5 C, 5.29 Theorem, pages 348-349]. The detailed proof falls outside the scope and depth of this paper and will therefore not be provided here. For a clearer explanation of Feller's boundary classification, it is strongly suggested to refer to the final section of the Appendix.\\

\noindent \textbf{Remark 4.3}\\
Note that the function $\psi(x)$ has a continuous, strictly positive-valued derivative, and $\psi''(x)$ exists almost everywhere and satisfies $\psi''(x)=-2\mu(x)[\nu(x)]^{-2}\psi'(x)$.$\hfill\blacksquare$\\

\noindent Now we compute $\psi(x)$ in our case. Computing the case for (4.9) leads to:
\begin{align}\nonumber
    \psi(x)&=\int_{c}^{x}\text{exp}\Big\{-2\int_{c}^{y}\frac{(\frac{\sigma^{2}}{2}z^{2k-1}-bz)dz}{\sigma^{2}z^{2k}}\Big\}dy=\int_{c}^{x}\text{exp}\Big\{\int_{c}^{y}\Big(-z^{-1}+\frac{2b}{\sigma^{2}}z^{1-2k}\Big)dz\Big\}dy\nonumber
\end{align}
We may let $c=1$ to alleviate the computation burden, so that 
\begin{equation}\nonumber
    \int_{1}^{y}-z^{-1}dz+\int_{1}^{y}\frac{2b}{\sigma^{2}}z^{1-2k}dz=-\text{log}y+\frac{b}{\sigma^{2}(1-k)}\Big(y^{2(1-k)}-1\Big).
\end{equation}
resulting in:
\begin{align}\nonumber
    \psi(x)=\int_{1}^{x}\text{exp}\Big\{-\text{log}y+\frac{b}{\sigma^{2}(1-k)}\Big(y^{2(1-k)}-1\Big)\Big\}dy
    =\text{exp}\Big\{\frac{-b}{\sigma^{2}(1-k)}\Big\}\int_{1}^{x}y^{-1}\text{exp}\Big\{\frac{b}{\sigma^{2}(1-k)}y^{2(1-k)}\Big\}dy.\nonumber
\end{align}
Let $M\overset{\text{def}}{=}\frac{b}{\sigma^{2}(1-k)}$. Since $b>0$, $1-k>0$, $\sigma^{2}>0$, it is easy to see that $M>0$, so:
\begin{align}\nonumber
    \psi(x)&=e^{-M}\int_{1}^{x}y^{-1}\text{exp}\{My^{2(1-k)}\}dy=e^{-M}\int_{1}^{x}y^{-1}\sum_{j=0}^{+\infty}\frac{M^{j}}{j!}y^{2j(1-k)}dy=e^{-M}\int_{1}^{x}\sum_{j=0}^{+\infty}\frac{M^{j}}{j!}y^{2j-2jk-1}dy.
\end{align}
For $x>1$, the summands are nonnegative and measurable, so Tonelli's theorem permits the interchange of summation and integration. For $0<x<1$, the same argument applies after reversing the integration limits.
\begin{align}\nonumber
    \psi(x)&=e^{-M}\int_{1}^{x}\sum_{j=0}^{+\infty}\frac{M^{j}}{j!}y^{2j-2jk-1}dy=e^{-M}\sum_{j=0}^{+\infty}\frac{M^{j}}{j!}\int_{1}^{x}y^{2j-2jk-1}dy\\
    &=e^{-M}\Big[\int_{1}^{x}\frac{dy}{y}+\sum_{j=1}^{+\infty}\frac{M^{j}}{j!}\int_{1}^{x}y^{2j-2jk-1}dy\Big]=e^{-M}\Big[\text{log}x+\sum_{j=1}^{+\infty}\frac{M^{j}}{j!}\Big(\frac{x^{2j(1-k)}}{2j(1-k)}-\frac{1}{2j(1-k)}\Big)\Big]\nonumber
\end{align}
Since $\frac{1}{2}\leq k<1$, we have $0<1-k\leq \frac{1}{2}$. For each $j\in\mathbb{N}$, we then have $0<2j(1-k)\leq j$ for $j\geq 1$, so $\frac{1}{2j(1-k)}x^{2j(1-k)}$ is a power function that has a positive-valued coefficient and a positive-valued power.\\

\noindent We now compute the limit of $\psi(x)$ when $x$ approaches $+\infty$. 
\begin{align}\nonumber
    \psi(x)=e^{-M}\Big[\text{log}x+\sum_{j=1}^{+\infty}\frac{M^{j}}{j!}\frac{x^{2j(1-k)}-1}{2j(1-k)}\Big]\geq e^{-M}\text{log}x,
\end{align}
since $\frac{x^{2j(1-k)}-1}{2j(1-k)}\geq0$ for every $j\geq 1$ when $x\geq 1$. Due to this, we have $\lim_{x\to+\infty}\psi(x)=+\infty$. This truly verifies that our $J$-valued diffusion process $\dot{\lambda}_{t}$ does not attain its upper boundary $+\infty$.\\

\noindent We now compute the limit of $\psi(x)$ as $x\to0$. Since $\frac{1}{2}\leq k<1$, we have $1-k>0$. Therefore, for $0<x<1$ and every integer $j\geq 1$, $0<x^{2j(1-k)}<1,~\frac{x^{2j(1-k)}-1}{2j(1-k)}<0$. Moreover, $M>0$, so $M^{j}/j!>0$. Thus, every term in the series appearing after $\text{log}x$ is nonpositive. Using the series representation established above, we obtain
\begin{equation}\nonumber
	\psi(x)=e^{-M}\Big[\text{log}x+\sum_{j=1}^{\infty}\frac{M^{j}}{j!}\frac{x^{2j(1-k)}-1}{2j(1-k)}\Big]\leq e^{-M}\text{log}x,\qquad 0<x<1.
\end{equation}
Since $e^{-M}>0$ and $\text{log}x\to-\infty$ as $x\downarrow0$, it follows that $\lim_{x\downarrow0}\psi(x)=-\infty$. Hence, the auxiliary diffusion $\dot{\lambda}$ cannot attain its lower boundary $0$ in finite time.\\

\noindent Combining the lower-boundary calculation with the upper-boundary calculation above gives $\lim_{x\downarrow0}\psi(x)=-\infty$, and $\lim_{x\to+\infty}\psi(x)=+\infty$. By Feller's test for explosions, these two limits imply that the auxiliary diffusion $\dot{\lambda}$ does not exit its state space $(0,+\infty)$ under $\dot{\mathbb{P}}$.\\

\noindent Recall that our parameter assumptions are $\frac{1}{2}\leq k<1$, $a,b,\sigma,L>0$, and $\lambda_{0}>0$, with the additional condition $2a\geq\sigma^{2}$ when $k=\frac{1}{2}$. Under these assumptions, \textbf{Theorem 1.2} ensures that the original CKLS process $\lambda$ remains strictly positive and does not explode in finite time under $\mathbb{P}$. Thus, the original process also does not exit $(0,+\infty)$. Since \textbf{Assumptions 1 and 2} have already been verified, \textbf{Theorem 4.1} applies and shows that the Dol\'{e}ans--Dade exponential
\begin{equation}\nonumber
	M_{t}=\text{exp}\Big\{\int_{0}^{t}q_{s}dW_{s}-\frac{1}{2}\int_{0}^{t}(q_{s})^{2}ds\Big\}, \qquad q(x)=\frac{\sigma}{2}x^{k-1}-\frac{a}{\sigma}x^{-k},
\end{equation}
is a true $\mathbb{P}$-martingale on every finite time interval with fixed $T>0$. In particular, since $M_{0}=1$, $\mathbb{E}_{\mathbb{P}}[M_{t}]=1$, $0\leq t\leq T$.\\

\noindent To verify that the resulting change of measure is equivalent, we also check the strict positivity of $M_{T}$. Fix $T>0$. By the continuity, strict positivity, and non-explosion of $\lambda$, almost surely,
\begin{equation}\nonumber
	0<\underline{\lambda}_{T}:=\min_{0\leq s\leq T}\lambda_{s}\leq\overline{\lambda}_{T}:=\max_{0\leq s\leq T}\lambda_{s}<+\infty.
\end{equation}
Since $q$ is continuous on $(0,+\infty)$, it is bounded on the compact interval $[\underline{\lambda}_{T},\overline{\lambda}_{T}]$. Consequently,
\begin{equation}\nonumber
	\int_{0}^{T}(q_{s})^{2}ds\leq T\sup_{x\in[\underline{\lambda}_{T},\overline{\lambda}_{T}]}q(x)^{2}<+\infty, \qquad \mathbb{P}\text{-a.s.}
\end{equation}
Thus, the stochastic exponential is strictly positive: $M_{T}>0,~\mathbb{P}\text{-a.s.}$. We now define a probability measure $\mathbb{Q}$ on $\mathcal{F}_{T}$ by
\begin{equation}\nonumber
	\mathbb{Q}(A):=\mathbb{E}_{\mathbb{P}}\big[M_{T}\mathbbm{1}_{A}\big],\qquad A\in\mathcal{F}_{T}, \quad\text{or, equivalently,}\quad \left.\frac{d\mathbb{Q}}{d\mathbb{P}}\right|_{\mathcal{F}_{T}}=M_{T}.
\end{equation}
The identity $\mathbb{E}_{\mathbb{P}}[M_{T}]=1$ ensures that $\mathbb{Q}$ is a probability measure. Moreover, because $M_{T}>0$ $\mathbb{P}-a.s.$, for every $A\in\mathcal{F}_{T}$,
\begin{equation}\nonumber
	\mathbb{Q}(A)=0 \quad\Longleftrightarrow\quad \mathbb{P}(A)=0.
\end{equation}
Hence, $\mathbb{Q}$ and $\mathbb{P}$ are equivalent on $\mathcal{F}_{T}$. Finally, the martingale property gives
\begin{equation}\nonumber
	M_{t}=\mathbb{E}_{\mathbb{P}}\big[M_{T}\mid\mathcal{F}_{t}\big], \qquad 0\leq t\leq T,
\end{equation}
so $M_{t}$ is indeed the density process of this change of measure throughout $[0,T]$. This establishes the claimed equivalent change of measure on every finite time interval.

\section*{Acknowledgements}
\noindent The authors gratefully acknowledge the insightful and constructive comments provided by the two anonymous reviewers, which have significantly enhanced the quality of this paper.

\section*{References}

\noindent \hangindent=1.5em \hangafter=1 \hypertarget{AitSahalia1999}{A\"it-Sahalia, Y. (1999). Transition densities for interest rate and other nonlinear diffusions. \textit{The Journal of Finance}, 54(4), 1361--1395.}

\noindent \hangindent=1.5em \hangafter=1 \hypertarget{AitSahalia2002}{A\"it-Sahalia, Y. (2002). Maximum likelihood estimation of discretely sampled diffusions: A closed-form approximation approach. \textit{Econometrica}, 70(1), 223--262.}

\noindent \hangindent=1.5em \hangafter=1 \hypertarget{Andersen2007}{Andersen, L. B. G. \& Piterbarg, V. V. (2007). Moment explosions in stochastic volatility models. \textit{Finance and Stochastics}, 11(1), 29--50.}

\noindent \hangindent=1.5em \hangafter=1 \hypertarget{Ahn1999}{Ahn, D. H. \& Gao, B. (1999). A parametric nonlinear model of term structure dynamics. \textit{The Review of Financial Studies}, 12(4), 721--762.}

\noindent \hangindent=1.5em \hangafter=1 \hypertarget{Alaya2012}{Alaya, M. B. \& Kebaier, A. (2012). Parameter estimation for the square-root diffusions: Ergodic and nonergodic cases. \textit{Stochastic Models}, 28(4), 609--634.}

\noindent \hangindent=1.5em \hangafter=1 \hypertarget{Alaya2013}{Alaya, M. B. \& Kebaier, A. (2013). Asymptotic behavior of the maximum likelihood estimator for ergodic and nonergodic square-root diffusions. \textit{Stochastic Analysis and Applications}, 31(4), 552--573.}

\noindent \hangindent=1.5em \hangafter=1 \hypertarget{Beuermann2005}{Beuermann, D., Antoniou, A. \& Bernales, A. (2005). \textit{The Dynamics of the Short‑Term Interest Rate in the UK}. Working paper. University Library of Munich. \url{https://ideas.repec.org/p/wpa/wuwpfi/0512029.html}. Accessed 17 July 2025.}

\noindent \hangindent=1.5em \hangafter=1 \hypertarget{Baldi2011}{Baldi, P. \& Caramellino, L. (2011). General Freidlin-Wentzell large deviations and positive diffusions. \textit{Statistics \& Probability Letters}, 81(8), 1218--1229.}

\noindent \hangindent=1.5em \hangafter=1 \hypertarget{BaroneAdesi1999}{Barone-Adesi, G., Allegretto, W., Dinenis, E. \& Sorwar, G. (1999). \textit{Valuation of derivatives based on CKLS interest rate models}. Working paper. Universit\`{a} della Svizzera italiana. \url{https://citeseerx.ist.psu.edu/document?repid=rep1&type=pdf&doi=e700b308909a99c6f1a88e6dcae2862457463671}. Accessed 17 July 2025.}

\noindent \hangindent=1.5em \hangafter=1 \hypertarget{Baxter1996}{Baxter, M. \& Rennie, A. (1996). \textit{Financial calculus: An introduction to derivative pricing}. Cambridge University Press.}

\noindent \hangindent=1.5em \hangafter=1 \hypertarget{Borodin2002}{Borodin, A. N. \& Salminen, P. (2002). \textit{Handbook of Brownian motion - facts and formulae}. Birkh\"{a}user Basel.}

\noindent \hangindent=1.5em \hangafter=1 \hypertarget{Bergstrom1984}{Bergstrom, A. R. (1984). Continuous time stochastic models and issues of aggregation over time. \textit{Handbook of Econometrics}, 2, 1145--1212.}

\noindent \hangindent=1.5em \hangafter=1 \hypertarget{Brennan1980}{Brennan, M. J. \& Schwartz, E. S. (1980). Analyzing convertible bonds. \textit{Journal of Financial and Quantitative Analysis}, 15(4), 907--929.}

\noindent \hangindent=1.5em \hangafter=1 \hypertarget{Brenner1996}{Brenner, R. J., Harjes, R. H. \& Kroner, K. F. (1996). Another look at models of the short-term interest rate. \textit{Journal of Financial and Quantitative Analysis}, 31(1), 85--107.}

\noindent \hangindent=1.5em \hangafter=1 \hypertarget{Broze1995}{Broze, L., Scaillet, O. \& Zakoian, J. M. (1995). Testing for continuous-time models of the short-term interest rate. \textit{Journal of Empirical Finance}, 2(3), 199--223.}

\noindent \hangindent=1.5em \hangafter=1 \hypertarget{Byers1998}{Byers, S. L. \& Nowman, K. B. (1998). Forecasting UK and US interest rates using continuous time term structure models. \textit{International Review of Financial Analysis}, 7(3), 191--206.}

\noindent \hangindent=1.5em \hangafter=1 \hypertarget{Cai2015}{Cai, Y. \& Wang, S. (2015). Central limit theorem and moderate deviation principle for CKLS model with small random perturbation. \textit{Statistics \& Probability Letters}, 98, 6--11.}

\noindent \hangindent=1.5em \hangafter=1 \hypertarget{Carr2007}{Carr, P. \& Sun, J. (2007). A new approach for option pricing under stochastic volatility. \textit{Review of Derivatives Research}, 10(2), 87--150.}

\noindent \hangindent=1.5em \hangafter=1 \hypertarget{Chan1992}{Chan, K. C., Karolyi, G. A., Longstaff, F. A. \& Sanders, A. B. (1992). An empirical comparison of alternative models of the short-term interest rate. \textit{The Journal of Finance}, 47(3), 1209--1227.}

\noindent \hangindent=1.5em \hangafter=1 \hypertarget{Choi2007}{Choi, Y. \& Wirjanto, T. S. (2007). An analytic approximation formula for pricing zero-coupon bonds. \textit{Finance Research Letters}, 4(2), 116--126.}

\noindent \hangindent=1.5em \hangafter=1 \hypertarget{Clark2011}{Clark, I. J. (2011). \textit{Foreign exchange option pricing: A practitioner's guide}. John Wiley \& Sons.}

\noindent \hangindent=1.5em \hangafter=1 \hypertarget{Cox1980}{Cox, J. C., Ingersoll, J. E., Jr. \& Ross, S. A. (1980). An analysis of variable rate loan contracts. \textit{The Journal of Finance}, 35(2), 389--403.}

\noindent \hangindent=1.5em \hangafter=1 \hypertarget{Cox1985}{Cox, J. C., Ingersoll, J. E., Jr. \& Ross, S. A. (1985). A theory of the term structure of interest rates. \textit{Econometrica}, 53(2), 385--407.}

\noindent \hangindent=1.5em \hangafter=1 \hypertarget{Cox1996}{Cox, J. C. (1996). The constant elasticity of variance option pricing model. \textit{Journal of Portfolio Management}, 23(5), 15--17.}

\noindent \hangindent=1.5em \hangafter=1 \hypertarget{Dahlquist1996}{Dahlquist, M. (1996). On alternative interest rate processes. \textit{Journal of Banking \& Finance}, 20(6), 1093--1119.}

\noindent \hangindent=1.5em \hangafter=1 \hypertarget{DeRossi2010}{De Rossi, G. (2010). Maximum likelihood estimation of the Cox-Ingersoll-Ross model using particle filters. \textit{Computational Economics}, 36(1), 1--16.}

\noindent \hangindent=1.5em \hangafter=1 \hypertarget{Dehtiar2022}{Dehtiar, O., Mishura, Y. \& Ralchenko, K. (2022). Two methods of estimation of the drift parameters of the Cox-Ingersoll-Ross process: Continuous observations. \textit{Communications in Statistics - Theory and Methods}, 51(19), 6818--6833.}

\noindent \hangindent=1.5em \hangafter=1 \hypertarget{Delbaen2002}{Delbaen, F. \& Shirakawa, H. (2002). A note on option pricing for the constant elasticity of variance model. \textit{Asia-Pacific Financial Markets}, 9, 85--99.}
 
\noindent \hangindent=1.5em \hangafter=1 \hypertarget{Dothan1978}{Dothan, L. U. (1978). On the term structure of interest rates. \textit{Journal of Financial Economics}, 6(1), 59--69.}

\noindent \hangindent=1.5em \hangafter=1 \hypertarget{Dokuchaev2017}{Dokuchaev, N. (2017). A pathwise inference method for the parameters of diffusion terms. \textit{Journal of Nonparametric Statistics}, 29(4), 731--743.}

\noindent \hangindent=1.5em \hangafter=1 \hypertarget{Engelbert1984}{Engelbert, H. J. \& Schmidt, W. (1984). On exponential local martingales connected with diffusion processes. \textit{Mathematische Nachrichten}, 119(1), 97--115.}

\noindent \hangindent=1.5em \hangafter=1 \hypertarget{Engelbert1991}{Engelbert, H. J. \& Schmidt, W. (1991). Strong Markov continuous local martingales and solutions of one-dimensional stochastic differential equations (Part III). \textit{Mathematische Nachrichten}, 151(1), 149--197.}

\noindent \hangindent=1.5em \hangafter=1 \hypertarget{Feller1952}{Feller, W. (1952). The parabolic differential equations and the associated semi-groups of transformations. \textit{Annals of Mathematics}, 55(3), 468--519.}

\noindent \hangindent=1.5em \hangafter=1 \hypertarget{Feng2012}{Feng, J., Fouque, J. P. \& Kumar, R. (2012). Small-time asymptotics for fast mean-reverting stochastic volatility models. \textit{Annals of Applied Probability}, 22(4), 1541--1575.}

\noindent \hangindent=1.5em \hangafter=1 \hypertarget{Girsanov1960}{Girsanov, I. V. (1960). On transforming a certain class of stochastic processes by absolutely continuous substitution of measures. \textit{Theory of Probability and Its Applications}, 5(3), 285--301.}

\noindent \hangindent=1.5em \hangafter=1 \hypertarget{GoingJaeschke2003}{G\"{o}ing-Jaeschke, A. \& Yor, M. (2003). A survey and some generalizations of Bessel processes. \textit{Bernoulli}, 9(2), 313--349.}

\noindent \hangindent=1.5em \hangafter=1 \hypertarget{Groisman2007}{Groisman, P. \& Rossi, J. D. (2007). Explosion time in stochastic differential equations with small diffusion. \textit{Electronic Journal of Differential Equations}, 2007(140), 1--9.}

\noindent \hangindent=1.5em \hangafter=1 \hypertarget{Gyongy1996}{Gy\"{o}ngy, I. \& Krylov, N. (1996). Existence of strong solutions for It\^{o}'s stochastic equations via approximations. \textit{Probability Theory and Related Fields}, 105(2), 143--158.}

\noindent \hangindent=1.5em \hangafter=1 \hypertarget{Holden1996}{Holden, H., {\O}ksendal, B., Ub{\o}e, J. \& Zhang, T. (1996). \textit{Stochastic partial differential equations: A modeling, white noise functional approach}. Boston, MA: Birkh\"{a}user Boston.}

\noindent \hangindent=1.5em \hangafter=1 \hypertarget{Hansen1982}{Hansen, L. P. (1982). Large sample properties of generalized method of moments estimators. \textit{Econometrica}, 50(4), 1029--1054.}

\noindent \hangindent=1.5em \hangafter=1 \hypertarget{Hu2003}{Hu, Y. \& {\O}ksendal, B. (2003). Fractional white noise calculus and applications to finance. \textit{Infinite Dimensional Analysis, Quantum Probability and Related Topics}, 6(1), 1--32.}

\noindent \hangindent=1.5em \hangafter=1 \hypertarget{Hu2015}{Hu, Y., Lan, G. \& Zhang, C. (2015). The explicit solution and precise distribution of CKLS model under Girsanov transform. \textit{International Journal of Statistics and Probability}, 4(1), 68--75.}
 
\noindent \hangindent=1.5em \hangafter=1 \hypertarget{Ma2008}{Ma, J., Wang, C. \& Bai, M. (2008). Modeling on interest rate behavior of RMB: Comparative research based on Vasicek, CIR and CKLS model. \textit{Proceedings of the 2008 IEEE International Conference on Automation and Logistics (ICAL)}. IEEE, 2941--2945.}

\noindent \hangindent=1.5em \hangafter=1 \hypertarget{Karlin1981}{Karlin, S. \& Taylor, H. E. (1981). \textit{A second course in stochastic processes}. Elsevier.}

\noindent \hangindent=1.5em \hangafter=1 \hypertarget{Karatzas2012}{Karatzas, I. \& Shreve, S. (2012). \textit{Brownian motion and stochastic calculus}. Springer Science \& Business Media.}

\noindent \hangindent=1.5em \hangafter=1 \hypertarget{Karatzas2016}{Karatzas, I. \& Ruf, J. (2016). Distribution of the time to explosion for one-dimensional diffusions. \textit{Probability Theory and Related Fields}, 164, 1027--1069.}

\noindent \hangindent=1.5em \hangafter=1 \hypertarget{Kazamaki1977}{Kazamaki, N. (1977). On a problem of Girsanov. \textit{Tohoku Mathematical Journal, Second Series}, 29(4), 597--600.}

\noindent \hangindent=1.5em \hangafter=1 \hypertarget{Khor2012}{Khor, C. Y., Pooi, A. H. \& Ng, K. H. (2012). Bond option pricing under the CKLS model. \textit{Proceedings of the 2nd Regional Conference on Applied and Engineering Mathematics (RCAEM-II)}, 179--183.}

\noindent \hangindent=1.5em \hangafter=1 \hypertarget{Khor2014}{Khor, C. Y. \& Pooi, A. H. (2014). Prediction of interest rate using CKLS model with stochastic parameters. \textit{AIP Conference Proceedings}, 1602(1), 467--472.}

\noindent \hangindent=1.5em \hangafter=1 \hypertarget{Krylov2002}{Krylov, N. (2002). \textit{A simple proof of a result of A. Novikov}. Preprint. arXiv. \url{https://doi.org/10.48550/arXiv.math/0207013}. Accessed 17 July 2025.}

\noindent \hangindent=1.5em \hangafter=1 \hypertarget{Kubilius2021}{Kubilius, K. \& Med\v{z}i{\=u}nas, A. (2021). Positive solutions of the fractional SDEs with non-Lipschitz diffusion coefficient. \textit{Mathematics}, 9(1), 18.}

\noindent \hangindent=1.5em \hangafter=1 \hypertarget{Lamberton2011}{Lamberton, D. \& Lapeyre, B. (2011). \textit{Introduction to stochastic calculus applied to finance}. Chapman and Hall/CRC.}

\noindent \hangindent=1.5em \hangafter=1 \hypertarget{Leon2013}{Le\'{o}n, J. A., Peralta Hern\'{a}ndez, L. \& Villa Morales, J. (2013). On the distribution of explosion time of stochastic differential equations. \textit{Bolet\'{i}n de la Sociedad Matem\'{a}tica Mexicana}, 19(2), 125--138.}

\noindent \hangindent=1.5em \hangafter=1 \hypertarget{Lewis2000}{Lewis, A. L. (2000). \textit{Option valuation under stochastic volatility}. Finance Press.}

\noindent \hangindent=1.5em \hangafter=1 \hypertarget{Lewis2016}{Lewis, A. L. (2016). \textit{Option valuation under stochastic volatility II}. Finance Press.}

\noindent \hangindent=1.5em \hangafter=1 \hypertarget{Lileika2020}{Lileika, G. \& Mackevi{\v{c}}ius, V. (2020). Weak approximation of CKLS and CEV processes by discrete random variables. \textit{Lithuanian Mathematical Journal}, 60, 208--224.}

\noindent \hangindent=1.5em \hangafter=1 \hypertarget{Liptser1977}{Liptser, R. S. \& Shiryaev, A. N. (1977). \textit{Statistics of random processes: General theory}. Springer.}

\noindent \hangindent=1.5em \hangafter=1 \hypertarget{Liptser2013}{Liptser, R. S. (2013). Bene\v{s} condition for a discontinuous exponential martingale. \textit{Journal of Mathematical Sciences}, 188, 717--723.}

\noindent \hangindent=1.5em \hangafter=1 \hypertarget{Lyu2025}{Lyu, Y. \& Nkurunziza, S. (2025). Inference methods in time-varying linear diffusion processes. \textit{Electronic Journal of Statistics}, 19(1), 1633--1680.}

\noindent \hangindent=1.5em \hangafter=1 \hypertarget{Mao2006}{Mao, X., Truman, A. \& Yuan, C. (2006). Euler-Maruyama approximations in mean-reverting stochastic volatility model under regime-switching. \textit{Journal of Applied Mathematics and Stochastic Analysis}, 2006(1), 080967.}

\noindent \hangindent=1.5em \hangafter=1 \hypertarget{Mao2007}{Mao, X. (2007). \textit{Stochastic differential equations and applications}. Elsevier.}

\noindent \hangindent=1.5em \hangafter=1 \hypertarget{Mao2013}{Mao, X. \& Szpruch, L. (2013). Strong convergence and stability of implicit numerical methods for stochastic differential equations with non-globally Lipschitz continuous coefficients. \textit{Journal of Computational and Applied Mathematics}, 238, 14--28.}

\noindent \hangindent=1.5em \hangafter=1 \hypertarget{Marie2014}{Marie, N. (2014). A generalized mean-reverting equation and applications. \textit{ESAIM: Probability and Statistics}, 18, 799--828.}

\noindent \hangindent=1.5em \hangafter=1 \hypertarget{Maruyama1954}{Maruyama, G. (1954). On the transition probability functions of Markov processes. \textit{Natural Science Report Ochanomizu University}, 5(1), 10--20.}

\noindent \hangindent=1.5em \hangafter=1 \hypertarget{Maruyama1955}{Maruyama, G. (1955). Continuous Markov processes and stochastic equations. \textit{Rendiconti del Circolo Matematico di Palermo}, 4, 48--90.}

\noindent \hangindent=1.5em \hangafter=1 \hypertarget{Mazzonetto2024}{Mazzonetto, S. \& Nieto, B. (2024). Parameters estimation of a threshold CKLS process from continuous and discrete observations. Preprint. HAL. \url{https://hal.science/hal-04524431v2/document}. Accessed 17 July 2025.}

\noindent \hangindent=1.5em \hangafter=1 \hypertarget{Mendy2013}{Mendy, I. (2013). Parametric estimation for sub-fractional Ornstein-Uhlenbeck process. \textit{Journal of Statistical Planning and Inference}, 143(4), 663--674.}

\noindent \hangindent=1.5em \hangafter=1 \hypertarget{Merton1974}{Merton, R. C. (1974). On the pricing of corporate debt: The risk structure of interest rates. \textit{The Journal of Finance}, 29(2), 449--470.}

\noindent \hangindent=1.5em \hangafter=1 \hypertarget{Mijatovic2012}{Mijatovi{\'c}, A. \& Urusov, M. (2012). On the martingale property of certain local martingales. \textit{Probability Theory and Related Fields}, 152(1), 1--30.}

\noindent \hangindent=1.5em \hangafter=1 \hypertarget{Mishura2022}{Mishura, Y., Ralchenko, K. \& Dehtiar, O. (2022). Parameter estimation in CKLS model by continuous observations. \textit{Statistics \& Probability Letters}, 184, 109391.}

\noindent \hangindent=1.5em \hangafter=1 \hypertarget{Naouara2016}{Naouara, N. \& Trabelsi, F. (2016). A short review on boundary behavior of linear diffusion processes. \textit{Graduate Journal of Mathematics}, 1(2), 138--149.}

\noindent \hangindent=1.5em \hangafter=1 \hypertarget{Newey1997}{Newey, W. K. \& Steigerwald, D. G. (1997). Asymptotic bias for quasi-maximum-likelihood estimators in conditional heteroskedasticity models. \textit{Econometrica}, 65(3), 587--599.}

\noindent \hangindent=1.5em \hangafter=1 \hypertarget{Novikov1972}{Novikov, A. (1972). On an identity for stochastic integrals. \textit{Theory of Probability and Its Applications}, 17(4), 717--720.}

\noindent \hangindent=1.5em \hangafter=1 \hypertarget{Nowman1997}{Nowman, K. B. (1997). Gaussian estimation of single-factor continuous time models of the term structure of interest rates. \textit{The Journal of Finance}, 52(4), 1695--1706.}
 
\noindent \hangindent=1.5em \hangafter=1 \hypertarget{Nowman1998}{Nowman, K. B. (1998). Continuous-time short term interest rate models. \textit{Applied Financial Economics}, 8(4), 401--407.}

\noindent \hangindent=1.5em \hangafter=1 \hypertarget{Nowman1999a}{Nowman, K. B. \& Sorwar, G. (1999a). Pricing UK and US securities within the CKLS model: Further results. \textit{International Review of Financial Analysis}, 8(3), 235--245.}

\noindent \hangindent=1.5em \hangafter=1 \hypertarget{Nowman1999b}{Nowman, K. B. \& Sorwar, G. (1999b). An evaluation of contingent claims using the CKLS interest rate model: an analysis of Australia, Japan, and the United Kingdom. \textit{Asia-Pacific Financial Markets}, 6, 205--219.}

\noindent \hangindent=1.5em \hangafter=1 \hypertarget{Nowman2005}{Nowman, K. B. \& Sorwar, G. (2005). Derivative prices from interest rate models: Results for Canada, Hong Kong, and United States. \textit{International Review of Financial Analysis}, 14(4), 428--438.}

\noindent \hangindent=1.5em \hangafter=1 \hypertarget{Overbeck1997}{Overbeck, L. \& Ryd\'{e}n, T. (1997). Estimation in the Cox-Ingersoll-Ross model. \textit{Econometric Theory}, 13(3), 430--461.}

\noindent \hangindent=1.5em \hangafter=1 \hypertarget{Overbeck1998}{Overbeck, L. (1998). Estimation for continuous branching processes. \textit{Scandinavian Journal of Statistics}, 25(1), 111--126.}

\noindent \hangindent=1.5em \hangafter=1 \hypertarget{Ozaki1992}{Ozaki, T. (1992). A bridge between nonlinear time series models and nonlinear stochastic dynamical systems: a local linearization approach. \textit{Statistica Sinica}, 2(1), 113--135.}

\noindent \hangindent=1.5em \hangafter=1 \hypertarget{Revuz2013}{Revuz, D. \& Yor, M. (2013). \textit{Continuous martingales and Brownian motion}. Vol. 293. Springer Science \& Business Media.}

\noindent \hangindent=1.5em \hangafter=1 \hypertarget{Shoji1998}{Shoji, I. \& Ozaki, T. (1998). Estimation for nonlinear stochastic differential equations by a local linearization method. \textit{Stochastic Analysis and Applications}, 16(4), 733--752.}

\noindent \hangindent=1.5em \hangafter=1 \hypertarget{Skorokhod2009}{Skorokhod, A. V. (2009). \textit{Asymptotic methods in the theory of stochastic differential equations}. American Mathematical Society.}

\noindent \hangindent=1.5em \hangafter=1 \hypertarget{Skouras2000}{Skouras, K. (2000). Strong consistency in nonlinear stochastic regression models. \textit{Annals of Statistics}, 28(3), 871--879.}

\noindent \hangindent=1.5em \hangafter=1 \hypertarget{Stehlikova2009}{Stehl{\'\i}kov{\'a}, B. \& \v{S}ev\v{c}ovi\v{c}, D. (2009). Approximate formulae for pricing zero-coupon bonds and their asymptotic analysis. \textit{International Journal of Numerical Analysis and Modeling}, 6(2), 274--283.}

\noindent \hangindent=1.5em \hangafter=1 \hypertarget{Stehlikova2013}{Stehl{\'\i}kov{\'a}, B. (2013). A simple analytic approximation formula for the bond price in the Chan-Karolyi-Longstaff-Sanders model. \textit{International Journal of Numerical Analysis and Modeling}, 4(3), 224--234.}

\noindent \hangindent=1.5em \hangafter=1 \hypertarget{Tangman2011}{Tangman, D. Y., Thakoor, N., Dookhitram, K. \& Bhuruth, M. (2011). Fast approximations of bond option prices under CKLS models. \textit{Finance Research Letters}, 8(4), 206--212.}

\noindent \hangindent=1.5em \hangafter=1 \hypertarget{Tsumurai2020}{Tsumurai, S. (2020). Malliavin differentiability of CEV-type Heston model. \textit{Journal of Mathematical Finance}, 10(1), 173--199.}

\noindent \hangindent=1.5em \hangafter=1 \hypertarget{VanSchuppen1974}{Van Schuppen, J. H. \& Wong, E. (1974). Transformation of local martingales under a change of law. \textit{Annals of Probability}, 2(5), 879--888.}

\noindent \hangindent=1.5em \hangafter=1 \hypertarget{Vasicek1977}{Va\v{s}\'{\i}\v{c}ek, O. (1977). An equilibrium characterization of the term structure. \textit{Journal of Financial Economics}, 5(2), 177--188.}

\noindent \hangindent=1.5em \hangafter=1 \hypertarget{Wei2020}{Wei, C. (2020). Parameter estimation for Chan-Karolyi-Longstaff-Saunders model driven by small L{\'e}vy noises from discrete observations. \textit{International Journal of Dynamical Systems and Differential Equations}, 10(4), 373--382.}

\noindent \hangindent=1.5em \hangafter=1 \hypertarget{Wu2008}{Wu, F., Mao, X. \& Chen, K. (2008). A highly sensitive mean-reverting process in finance and the Euler-Maruyama approximations. \textit{Journal of Mathematical Analysis and Applications}, 348(1), 540--554.}

\noindent \hangindent=1.5em \hangafter=1 \hypertarget{Xu2015}{Xu, R., Wu, D. \& Yi, R. (2015). Existence theorem for mean-reverting CEV process with regime switching. \textit{Proceedings of the 2015 International Conference on Mechatronics, Electronic, Industrial and Control Engineering (MEIC-15)}. Atlantis Press, 1560--1563.}

\noindent \hangindent=1.5em \hangafter=1 \hypertarget{Yamada1971}{Yamada, T. \& Watanabe, S. (1971). On the uniqueness of solutions of stochastic differential equations. \textit{Journal of Mathematics of Kyoto University}, 11(1), 155--167.}

\noindent \hangindent=1.5em \hangafter=1 \hypertarget{Yang2020}{Yang, H., Wu, F., Kloeden, P. E. \& Mao, X. (2020). The truncated Euler-Maruyama method for stochastic differential equations with H\"{o}lder diffusion coefficients. \textit{Journal of Computational and Applied Mathematics}, 366, 112379.}

\noindent \hangindent=1.5em \hangafter=1 \hypertarget{Zahle1998}{Z\"{a}hle, M. (1998). Integration with respect to fractal functions and stochastic calculus. I. \textit{Probability Theory and Related Fields}, 111, 333--374.}

\noindent \hangindent=1.5em \hangafter=1 \hypertarget{Zhao2022}{Zhao, J. \& Xu, Z. (2022). Simultaneous identification of volatility and mean-reverting parameter for European option under fractional CKLS model. \textit{Fractal and Fractional}, 6(7), 344.}

\section*{Appendix}
\noindent \textit{Proof of }\textbf{Theorem 1.2}: The proof will be presented in 3 parts separately as follows.\\
\noindent \textbf{Part \RomanNumeralCaps{1}: Ranges of the solution in (1)-(4)}:\\

\noindent \textit{Proof} \textbf{of (1)} This is another application of \textbf{Theorem 4.2} (Feller's test for explosion) and its generalization - the boundary classification criteria, which describe the boundary behaviors of some diffusion processes of prescribed types. For a clearer explanation, see the last part of the appendix. One may also refer to \hyperlink{Karlin1981}{Karlin and Taylor (1981)} [Chapter 15 \S 6, pages 226-242] or \hyperlink{Naouara2016}{Naouara and Trabelsi (2016)} [Section 2.1, Lemma 2.1, Lemma 2.3, pages 143-144] or \hyperlink{Borodin2002}{Borodin and Salminen (2002)} [Chapter 2, pages 14-15], in which the theory of boundary classification for regular (linear) diffusion processes is explained elaborately. To be specific, one may adopt either the classical Feller boundary classification criteria (as in this paper) or the classical Russian (Gikhman and Skorokhod) boundary classification criteria to complete the proof.\qed\\ 

\noindent \textit{Proof} \textbf{of (3)}: For $0<k<\frac{1}{2}$, consider the diffusion on $(0,+\infty)$ up to its first hitting time of zero. Fix $c>0$ and define its scale density by
$\rho(y)=\text{exp}\big\{-2\int_{c}^{y}\frac{a-bz}{\sigma^{2}z^{2k}}dz\big\}$. Since $2k<1$, the integral in the exponent has a finite limit as $y\downarrow0$. Hence $\rho(y)\longrightarrow\rho_{0}\in(0,+\infty)$. Consequently, $\int_{0}^{c}\rho(y)dy<+\infty$ and $\int_{0}^{c}\frac{2}{\sigma^{2}y^{2k}\rho(y)}dy<+\infty$. The scale and speed integrals are therefore finite near zero, so zero is a regular, attainable boundary. A boundary condition must be specified to select an extension of the interior diffusion after it reaches zero. For the reflected process considered here, we choose the instantaneously reflecting extension. See also \hyperlink{Andersen2007}{Andersen and Piterbarg (2007)} [Section 2, Proposition 2.1, Proposition 2.2, pages 32-33].\qed\\

\noindent \textit{Proof} \textbf{of (4)}: For the borderline exponent $k=\frac{1}{2}$ with $2a \geq \sigma^{2}$, which turns out to be the case when the CKLS process degenerates to the CIR process, it is possible to show that the solution never reaches the origin. See the proof of \textbf{Lemma 3.1}. Briefly speaking, we can fix $c\in\mathbb{R}_{+}$ and observe that $\lim_{x\to +\infty}\psi(x)=+\infty$ and $\lim_{x\to 0+}\psi(x)=-\infty$, where $\psi(x)=\int_{c}^{x}\text{exp}\{-2\int_{c}^{y}\frac{a-bz}{\sigma^{2}z}dz\}dy$. In contrast, for the case $k=\frac{1}{2}$ with $2a<\sigma^{2}$, one can compute that for all $c\in(0,+\infty)$, $\lim_{x\downarrow 0}\phi(x)=\lim_{x\downarrow 0}\int_{c}^{x}\psi'(y)\int_{c}^{y}\frac{2dz}{\psi'(z)[\nu(z)]^{2}}dy<+\infty$, indicating that the origin is an attainable boundary. In contrast, when $2a<\sigma^{2}$, $\lambda_{t}$ will reach $0$  with probability one (see also the proof of \textbf{Lemma 3.1}).\qed\\

\noindent Moreover, by invoking \textbf{Lemma 3.5}, the CIR process can be mapped, via an appropriate space-time change and It\^{o}'s lemma, to a squared Bessel process of dimension $d=\frac{4a^{*}b^{*}}{\sigma^{*2}}$. In other words, the CIR process is a scaled, drift-adjusted version of a $\mathrm{BESQ}_{d,r_{0}}$ process. This connection allows key properties of the CIR process — such as boundary behavior at zero, recurrence, and the existence of a stationary distribution — to be analyzed using classical results from the theory of Bessel processes, see \hyperlink{Revuz2013}{Revuz and Yor (2013)} [Chapter XI, (1.5) Proposition, page 442]. As a result, it is shown that within this context, for $k=\frac{1}{2}$ with $2a<\sigma^{2}$, the origin acts as a strong reflector (i.e. the origin is strongly reflecting). This means that the time spent by the process at $\lambda_{t}=0$ has a Lebesgue measure zero, and no explicit boundary condition at $\lambda_{t}=0$ is required. In other words, while the process may reach the boundary point $0$, it immediately reflects and moves into the positive interior. Therefore, a stationary distribution for $\lambda_{t}$ is expected to exist in this case, again without the need for an explicit boundary condition at the origin. For further details on the case $2a<\sigma^{2}$, see \hyperlink{Andersen2007}{Andersen and Piterbarg (2007)} [Section 2, Proposition 2.1 and Proposition 2.2, pages 32-33].\\

\noindent \textit{Proof} \textbf{of (2)}: In the following several paragraphs, we provide a detailed proof for the case $k>\frac{1}{2}$. Note that for all values of $k$, the unattainability at $+\infty$ has already been discussed in the first paragraph. Let $\tau_{\lambda}=\text{inf}\{t\geq 0;\lambda_{t}=0 \text{ or } \lambda_{t}=+\infty\}$ with $\text{inf}\{\varnothing\}=+\infty$. We want to verify that for a fixed number $c\in\mathbb{R}_{+}$, $\lim_{x\to +\infty}\psi(x)=+\infty$ and $\lim_{x\to 0+}\psi(x)=-\infty$, $\forall k>\frac{1}{2}$, where $\psi(x)=\int_{c}^{x}\text{exp}\{-2\int_{c}^{y}\frac{a-bz}{\sigma^{2}z^{2k}}dz\}dy$. If so, we can readily claim $\mathbb{P}(\tau_{\lambda}=+\infty)=1$. We can assume without loss of generality that $c=1$.\\

\noindent \textbf{Case \RomanNumeralCaps{1}: Assume firstly that} $k\neq1$, we have:
\begin{align}\nonumber
-2\int_{1}^{y}\frac{a-bz}{\sigma^{2} z^{2k}}dz&=-2\int_{1}^{y}\Big(\frac{a}{\sigma^{2} z^{2k}}-\frac{b}{\sigma^{2} z^{2k-1}}\Big)dz=\frac{-2a}{\sigma^{2}(1-2k)}[z^{1-2k}]^{y}_{1}-\frac{-2b}{\sigma^{2}(2-2k)}[z^{2-2k}]^{y}_{1}\nonumber\\
&=\frac{2a}{\sigma^{2}(2k-1)}\Big(y^{1-2k}-1\Big)-\frac{b}{\sigma^{2}(k-1)}\Big(y^{2-2k}-1\Big)\overset{\text{def}}{=}I(y).\nonumber
\end{align}

$\bm{\langle 1.1\rangle}$ \textbf{Upper boundary calculation for} $k>1$:\\
Since $\lim_{y\to+\infty}y^{1-2k}=0$ (as $1-2k<-1$), $\lim_{y\to+\infty}y^{2-2k}=0$ (as $2-2k<0$), we have:
\begin{align}\nonumber
    \lim_{y\to+\infty}I(y)=\frac{-2a}{\sigma^{2}(2k-1)}+\frac{b}{\sigma^{2}(k-1)}\overset{\text{def}}{=}C_{1}\in\mathbb{R}.
\end{align}

$\bm{\langle 1.2\rangle}$ \textbf{Upper boundary calculation for} $\frac{1}{2}<k<1$:
\begin{equation}\nonumber
    I(y)=\frac{2a}{\sigma^{2}(2k-1)}\Big(\frac{1}{y^{2k-1}}-1\Big)-\frac{b}{\sigma^{2}(k-1)}y^{2-2k}+\frac{b}{\sigma^{2}(k-1)}>\frac{2a}{\sigma^{2}(2k-1)}\Big(\frac{1}{y^{2k-1}}-1\Big)+\frac{b}{\sigma^{2}(k-1)}\overset{\text{def}}{=}i(y),~y>1.
\end{equation}
The sign $>$ requires additionally that $k<1$, which results in $-\frac{b}{\sigma^{2}(k-1)}>0$.  Since $y\in(1,+\infty)$, $y^{2-2k}>0$, we have $-\frac{b}{\sigma^{2}(k-1)}y^{2-2k}>0$ (in fact, this term possibly diverges as $y$ approaches $+\infty$, but this does not influence the final result to be proven). Also, since $2k-1>0$ when $k>\frac{1}{2}$, we have $\frac{2a}{\sigma^{2}(2k-1)}>0$. Taking the limit $y\to+\infty$ gives $\lim_{y\to+\infty}\frac{1}{y^{2k-1}}=0$ and thus $\lim_{y\to+\infty}i(y)=-\frac{2a}{\sigma^{2}(2k-1)}+\frac{b}{\sigma^{2}(k-1)}=C_{2}\in\mathbb{R}_{-}$.\\

\noindent As a result, for each $j=1,2$, there exists some $R_{j}>1$ such that, for $y\geq R_{j}$, $I(y)>C'_{j}:=\min(0,C_{j})-1$. Hence, $e^{I(y)}>e^{C'_{j}}=:C_{j}^{*}>0$. Thus, for $x\geq R_{j}$,
\begin{align}\nonumber
    \psi(x)=\psi(R_{j})+\int_{R_{j}}^{x}e^{I(y)}dy\geq \psi(R_{j})+C_{j}^{*}(x-R_{j})\xrightarrow[x\to+\infty]{}+\infty.
\end{align}

$\bm{\langle 2.1\rangle}$ \textbf{Lower boundary calculation for} $\frac{1}{2}<k<1$:\\
Note that $0<2k-1<1$ and $-1<2k-2<0$. For $0<y<1$, we have $0<y^{2-2k}<1$. Let $K_{1}:=\frac{2a}{\sigma^{2}(2k-1)}>0$ and $K_{2}:=\frac{b}{\sigma^{2}(1-k)}>0$. Then $-K_{2}<-\frac{b}{\sigma^{2}(k-1)}(y^{2-2k}-1)<0$, and therefore $I(y)>K_{1}(y^{1-2k}-1)-K_{2}$.
\begin{align}\nonumber
    \psi(x)&=\int_{1}^{x}e^{I(y)}dy=-\int_{x}^{1}e^{I(y)}dy<-\int_{x}^{1}\text{exp}\Big\{{K_{1}y^{1-2k}-K_{1}-K_{2}}\Big\}dy\\
    &\overset{y=\frac{1}{w}}{=}-e^{-K_{1}-K_{2}}\int_{1}^{\frac{1}{x}}\frac{1}{w^{2}}\text{exp}\Big\{K_{1}w^{2k-1}\Big\}dw\xrightarrow[x\to0+]{}-\infty.\nonumber
\end{align}
The last limit follows from $\frac{\text{exp}\{K_{1}w^{2k-1}\}}{w^{2}}\to+\infty$, $\text{as }w\to+\infty$, since $K_{1}>0$ and $2k-1>0$.\\

$\bm{\langle 2.2\rangle}$ \textbf{Lower boundary calculation for} $k>1$:\\
Note that $1-2k<-1$ and $2-2k<0$, so both $y^{1-2k}$ and $y^{2-2k}$ tend to $+\infty$ as $y\downarrow0$. However,$y^{2-2k}/y^{1-2k}=y\to0$. Thus,
$I(y)\sim K_{1}y^{1-2k}$, $K_{1}=\frac{2a}{\sigma^{2}(2k-1)}>0$, as $y\downarrow0$, and consequently $\lim_{y\downarrow0}I(y)=+\infty$. More precisely, there exists $\delta\in(0,1)$ such that $I(y)\geq\frac{K_{1}}{2}y^{1-2k}\geq 2\text{log}\frac{1}{y},~0<y\leq\delta$, because $y^{1-2k}/\text{log}(1/y)\to+\infty$ as $y\downarrow0$. Thus, for $0<x<\delta$,
\begin{equation}\nonumber
    \psi(x)=-\int_x^1e^{I(y)}dy\leq-\int_x^\delta e^{I(y)}dy\leq-\int_x^\delta\frac1{y^{2}}dy=-\frac1x+\frac1\delta\xrightarrow[x\downarrow0]{}-\infty.
\end{equation}

\noindent \textbf{Case \RomanNumeralCaps{2}: Now we assume} $k=1$. We have: 
\begin{align}\nonumber
    -2\int_{1}^{y}\frac{a-bz}{\sigma^{2} z^{2k}}dz=-2\int_{1}^{y}\frac{a-bz}{\sigma^{2} z^{2}}dz=\frac{2a}{\sigma^{2}}\Big[\frac{1}{z}\Big]^{y}_{1}+\frac{2b}{\sigma^{2}}\Big[\text{log}z\Big]^{y}_{1}=\frac{2a}{\sigma^{2}}\Big(\frac{1}{y}-1\Big)+\frac{2b}{\sigma^{2}}\text{log}y=:I^{*}(y).
\end{align}

\textbf{Lower boundary calculation}: When $y\to 0+$, $\text{exp}\{\frac{2a}{\sigma^{2}}(\frac{1}{y}-1)\}$ explodes to $+\infty$ at an exponential rate, while $\text{exp}\{\frac{2b}{\sigma^{2}}\text{log}y\}=y^{\frac{2b}{\sigma^{2}}}$ decays to $0$ at a polynomial rate. Together, $\frac{I^{*}(y)}{\text{log}(1/y)}\to+\infty\text{ as }y\downarrow0$. Hence, there exists $\delta\in(0,1)$ such that $I^{*}(y)\geq2\text{log}\frac{1}{y},~0<y\leq\delta$. Consequently, for $0<x<\delta$,
\begin{equation}\nonumber
    \psi(x)=-\int_{x}^{1}e^{I^{*}(y)}dy\leq-\int_{x}^{\delta}\frac1{y^{2}}dy=-\frac{1}{x}+\frac{1}{\delta}\xrightarrow[x\downarrow0]{}-\infty.
\end{equation}

\textbf{Upper boundary calculation}: When $y\to +\infty$, $\text{exp}\{\frac{2a}{\sigma^{2}}(\frac{1}{y}-1)\}$ converges to $\text{exp}\{-\frac{2a}{\sigma^{2}}\}$ while $\text{exp}\{\frac{2b}{\sigma^{2}}\text{log}y\}=y^{\frac{2b}{\sigma^{2}}}$ diverges to $+\infty$. Together,
\begin{align}\nonumber
    &\lim_{y\to+\infty}I^{*}(y)=\lim_{y\to+\infty}\Big[-2\int_{1}^{y}\frac{a-bz}{\sigma^{2} z^{2k}}dz\Big]=+\infty\\
    \xRightarrow{}&~\text{$e^{I^{*}(y)}\geq1$ for all sufficiently large $y$.}\nonumber\\
    \xRightarrow{}&\lim_{x\to +\infty}\psi(x)=\lim_{x\to +\infty}\int_{1}^{x}e^{I^{*}(y)}dy=+\infty.\nonumber
\end{align}
\textbf{Combining Case I and Case II, we arrive at the conclusion}: $\lim_{x\to +\infty}\psi(x)=+\infty$, $\lim_{x\to 0+}\psi(x)=-\infty$, for all $k>\frac{1}{2}$. By Feller's criterion, we proved \textbf{(2)}. \qed\\

\noindent \textbf{Remark*}: We would also like to highlight a particular paper in the literature here, which gives a different method for the proof: For the case $k\in(\frac{1}{2},1)$ and the case $k=\frac{1}{2}$ with $2a\geq\sigma^{2}$, \hyperlink{Xu2015}{Xu et al. (2015)} prove the global existence and strict positivity of the (regime-switching) CKLS process $\lambda_{t}$ in an innovative way by constructing a Lyapunov function $V(\lambda_{t})=\theta_{1} (\lambda_{t})^{\frac{1}{2}}+\theta_{2}(\lambda_{t})^{-2}$ tailored to control the process both near zero and at infinity for some $\theta_{1}$ and $\theta_{2}$. Applying It\^{o}'s lemma to this function yields a bound on its expected growth through the generator $LV(\lambda_{t})$, which is shown to be at most linear in $V(\lambda_{t})$. A contradiction argument is then employed: assuming the process hits the boundary with positive probability (i.e. $\mathbb{P}(\tau_{\lambda}<+\infty)>0)$ leads to a divergence in the expected Lyapunov values, which contradicts the boundedness derived through Gr{\"{o}}nwall's inequality. This contradiction implies that the stopping time associated with boundary exit is almost surely infinite (i.e. $\mathbb{P}(\tau_{\lambda}=+\infty)=1)$.\\

\noindent \textbf{Part \RomanNumeralCaps{2}: Uniqueness and Strongness/Weakness of the solution in (1)-(4)}
Note: This section does not address the first item of \textbf{Theorem 1.2}. Consequently, we only need to prove items (2)-(4), excluding the previous part concerning the ranges of the solution. Henceforth, we will denote the items concerning uniqueness and strongness/weakness of the solution as \textbf{(2$\ast$)}-\textbf{(4$\ast$)}.\\

\noindent \textit{Proof} \textbf{of (2$\ast$)-(4$\ast$)}: \textbf{Yamada-Watanabe-Engelbert Theorem} (commonly known as Yamada's condition, see the original paper by \hyperlink{Yamada1971}{Yamada and Watanabe (1971)}, or \hyperlink{Revuz2013}{Revuz and Yor (2013)} [Chapter IX, Theorem 3.5, page 390]) claims:
\begin{mdframed}[
  linecolor=black,linewidth=2pt,
  roundcorner=5pt,
  backgroundcolor=white
]
Consider the stochastic differential equation $dZ_{t}=\mu(t,Z_{t})dt+\nu(t,Z_{t})dW_{t}, t\in[0,+\infty)$ with $Z_{t}$ defined on some filtered probability space, assume that there exists a constant $\tau^{*}>0$, a constant $A$ and a function $B:[0,\tau^{*}]\xrightarrow{}[0,+\infty)$ such that $\lvert\mu(t,x)-\mu(t,y)\lvert\leq A\lvert x-y\lvert$ (Lipschitz continuous) and $\lvert\nu(t,x)-\nu(t,y)\lvert\leq B(\lvert x-y\lvert)$, $\forall t\in[0,+\infty)$ (H\"{o}lder continuous) where $B(u)$ should be non-decreasing, strictly positive-valued $\forall u\in(0,\tau^{*}]$, and its square should satisfy the Osgood condition (see \hyperlink{Leon2013}{Le{\'o}n et al. (2013}) and \hyperlink{Groisman2007}{Groisman and Rossi (2007)}): $\int_{0+}^{\tau^{*}}\frac{1}{B^{2}(u)}du=+\infty$. Then the strong uniqueness of $Z_{t}$ is ensured. 
\end{mdframed}

\noindent In our case, for $k>\frac{1}{2}$, we may let $A=b+1$ with $b>0$, then it is checked that: 
\begin{equation}\nonumber
    \text{(Lipschitz continuous) }\lvert\mu(x)-\mu(y)\lvert=\lvert(a-bx)-(a-by)\lvert=\lvert-b(x-y)\lvert\leq A\lvert x-y\lvert.
\end{equation}
Note that for any differential function $F$ on $\mathbb{R}$, let $u=y+\delta(x-y)$ and thus $du=(x-y)d\delta$ with $\delta\in[0,1]$: 
\begin{equation}\nonumber
    F(x)-F(y)=(x-y)\int_{0}^{1}F'(\delta x+(1-\delta)y)d\delta,~x>0 \text{ and } y\geq 0,
\end{equation}
we have, for $F(u)=u^{k}$:
\begin{equation}\nonumber
    \lvert x^{k}-y^{k}\lvert=\lvert x-y\lvert\int_{0}^{1} k\big(\delta x+(1-\delta)y\big)^{k-1}d\delta.
\end{equation}
When $\frac{1}{2}<k<1$, the function $u\mapsto u^{k-1}$ is decreasing on $(0,\infty)$. The case $x=y$ is immediate, so we may assume without loss of generality that $x>y\geq0$. Then, for $\delta\in(0,1]$, we have
\begin{equation}\nonumber
    \delta x+(1-\delta)y \geq \delta(x-y)
    \quad\Longrightarrow\quad
    \big(\delta x+(1-\delta)y\big)^{k-1} \leq \big(\delta \lvert x-y\lvert\big)^{k-1}.
\end{equation}
Plugging this into the previous display yields the following:
\begin{equation}\nonumber
    \lvert x^{k}-y^{k}\lvert \leq \lvert x-y\lvert\int_{0}^{1} k\big(\delta \lvert x-y\lvert \big)^{k-1} d\delta=\lvert x-y\lvert^{k}.
\end{equation}
Hence for the diffusion term,
\begin{equation}\nonumber
    \lvert\nu(x)-\nu(y)\lvert=\sigma\lvert x^{k}-y^{k}\lvert\leq \sigma \lvert x-y\lvert^{k}:=B(\lvert x-y\lvert),\quad B(u):=\sigma u^{k}.
\end{equation}
Here $B$ depends only on $u=\lvert x-y\lvert$, is non-decreasing, concave, and satisfies $B(0)=0$.\\

\noindent \textbf{Osgood condition (this is where $k>\frac{1}{2}$ is used):} On any fixed interval $[0,\tau^{*})$:
\begin{equation}\nonumber
    \int_{0+}^{\tau^{*}}\frac{du}{(B(u))^{2}}=\frac{1}{\sigma^{2}}\int_{0+}^{\tau^{*}}u^{-2k}du=\frac{1}{\sigma^{2}(1-2k)}\Big[(\tau^{*})^{1-2k}-\lim_{u\to0+}u^{1-2k}\Big]=+\infty, \quad\text{for}~\frac{1}{2}<k<1,
\end{equation}
because $-1<1-2k<0$ and $\lim_{u\to0+}u^{1-2k}=+\infty$. Together with the Lipschitz drift, Yamada-Watanabe then gives pathwise uniqueness.\\

\noindent When $k\geq 1$, the mean value theorem on any bounded interval $[0,\tau^{*}]$ (locally) gives
\begin{equation}\nonumber
    \lvert x^{k}-y^{k}\lvert\leq k (\tau^{*})^{k-1}\lvert x-y\lvert,
\end{equation}
hence
\begin{equation}\nonumber
    \lvert \nu(x)-\nu(y)\lvert \leq B_{\tau^{*}}(\lvert x-y\lvert), \quad B_{\tau^{*}}(u):=\sigma k (\tau^{*})^{k-1}u=Mu.
\end{equation}
and 
\begin{equation}\nonumber
    \int_{0+}^{\tau^{*}}\frac{du}{(B(u))^{2}}=\frac{1}{M^{2}}\int_{0+}^{\tau^{*}}u^{-2}du=\frac{-1}{M^{2}}\Big[(\tau^{*})^{-1}-\lim_{u\to0+}u^{-1}\Big]=+\infty.
\end{equation}
Since $\tau^{*}$ is arbitrary, the pathwise uniqueness holds up to the first time a trajectory leaves $[0,\tau^{*}]$; by increasing $\tau^{*}$ step by step to infinity, we extend the pathwise uniqueness to the entire time axis.\qed\\

\noindent Meanwhile, when $k=\frac{1}{2}$ (either $2a\geq\sigma^{2}$ or $2a<\sigma^{2}$), $\lambda_{t}$ is a pathwise unique strong solution over $[0,+\infty)$, because the drift function $\mu(x)=a-bx$ is Lipschitz continuous and the diffusion function $\nu(x)=\sigma x^{\frac{1}{2}}$ is H\"{o}lder continuous, leading to the existence of a constant $C>0$ such that $\lvert \nu(x)-\nu(y)\lvert \leq C (\lvert x-y\lvert)^{\frac{1}{2}}$, for all $x>0$ and $y\geq 0$. See also \hyperlink{Andersen2007}{Andersen and Piterbarg (2007)} [Section 2, Proposition 2.1, Proposition 2.2, pages 32-33].\qed\\

\noindent \textbf{Part \RomanNumeralCaps{3}: Proofs of (5)-(7)}\\
\noindent \textit{Proof} \textbf{of (5)}: The result represents a straightforward application of the general ergodic theory applicable to homogeneous diffusion processes, see \hyperlink{Skorokhod2009}{Skorokhod (2009)} [Chapter 1, §3, Theorem 16, page 46]. Having supposed that the process $\lambda_{t}$ is instantaneously reflected at the origin for the case $0<k<\frac{1}{2}$, we know that the specific boundary condition that the transition density function $p(t,\lambda_{0},x)$ should satisfy, which is the zero-flux boundary condition for the case $\lambda_{t}=x$ with initial value $\lambda_{0}$, is:
\begin{equation}\nonumber
    \lim_{x\downarrow 0}\Big\{\frac{\partial}{\partial x}\Big(\frac{\sigma^{2}x^{2k}}{2}p(t,\lambda_{0},x)\Big)-(a-bx)p(t,\lambda_{0},x)\Big\}=0.
\end{equation}
This is the zero-flux boundary condition at the origin. To determine the stationary density $p_{\infty}(x)$, we solve the stationary Fokker–Planck–Kolmogorov equation
\begin{equation}\nonumber
    0=-\frac{d}{dx}\big[(a-bx)p_{\infty}(x)\big]+\frac{\sigma^{2}}{2}\frac{d^{2}}{dx^{2}}\big[x^{2k}p_{\infty}(x)\big], \qquad x>0.
\end{equation}
Integrating once gives
\begin{equation}\nonumber
    (a-bx)p_{\infty}(x)-\frac{\sigma^{2}}{2}\frac{d}{dx}\big[x^{2k}p_{\infty}(x)\big]=J,
\end{equation}
where $J$ is a constant probability flux. Under the zero-flux condition, $J=0$. Expanding the derivative therefore gives:
\begin{equation}\nonumber
   (a-bx)p_{\infty}=\frac{\sigma^{2}}{2}\Big(2kx^{2k-1}p_{\infty}+x^{2k}\frac{dp_{\infty}}{dx}\Big).
\end{equation}
This turns out to be:
\begin{align}\nonumber
    \frac{2}{\sigma^{2}}(ax^{-2k}-bx^{1-2k})p_{\infty}-2kx^{-1}p_{\infty}&=\frac{dp_{\infty}}{dx}\\
    \frac{2}{\sigma^{2}}\Big(\frac{a}{1-2k}x^{1-2k}-\frac{b}{2-2k}x^{2-2k}\Big)-2k\text{log}x+C&=\text{log}p_{\infty}\tag{$\ast$}\nonumber\\
    p_{\infty}\propto x^{-2k}e^{\frac{2}{\sigma^{2}}(\frac{a}{1-2k}x^{1-2k}-\frac{b}{2-2k}x^{2-2k})}&.\nonumber
\end{align}
Note that ($\ast$) holds if and only if $k\neq\frac{1}{2}$ and $k\neq 1$. When $k=\frac{1}{2}$, $\int ax^{-2k}dx$ should equal $a\text{log}x$ and $\int bx^{1-2k}dx=bx$ (integral constants omitted, the same hereinafter); When $k=1$, $\int bx^{1-2k}dx$ should equal $b\text{log}x$ and $\int ax^{-2k}dx=-ax^{-1}$, respectively. In the end, by direct computations of the scale function and the speed measure of $\lambda_{t}$ (see the part \textbf{Feller's boundary classification} in the Appendix for these two concepts), $C_{k}$ is obtained by normalizing the density: $C_{k}=\big(\int_{0}^{\infty}x^{-2k}e^{\Lambda(x;k)}dx\big)^{-1}$ provided that $0<\int_{0}^{\infty}x^{-2k}e^{\Lambda(x;k)}dx<\infty$.
See also \hyperlink{Andersen2007}{Andersen and Piterbarg (2007)} [Section 2, Proposition 2.1, Proposition 2.2, pages 32-33].\qed\\

\noindent \textit{Proof} \textbf{of (6)}:\\
For this property, we don't even need to restrict the drift coefficient to be linear, i.e. $\mu(x)=a-bx$ as in the CKLS model (Recall \textbf{Remark 1.3 (1)}). This condition can be relaxed to any globally Lipschitz continuous function $\mu(x)$ satisfying $\mu(0)>0$. Throughout this proof, fix $T>0$, $k\in(\frac{1}{2},1)$, and $\lambda_{0}>0$.\\

\noindent \textbf{Firstly, we consider positive moments with $p\geq 2$}. For sufficiently large integers $n$, define the stopping time $\tau_{n}=\text{inf}\{t\geq0:\lambda_{t}\notin(1/n,n)\}$ with $\text{inf}\{\varnothing\}=+\infty$. By It\^{o}'s lemma, we have

\begin{align}\nonumber
&(\lambda_{t\wedge\tau_{n}})^{p}=(\lambda_{0})^{p}+\int_{0}^{t\wedge\tau_{n}}p(\lambda_{s})^{p-1}d\lambda_{s}+\frac{1}{2}\int_{0}^{t\wedge\tau_{n}}p(p-1)(\lambda_{s})^{p-2}(d\lambda_{s})^{2}\\
\leq&(\lambda_{0})^{p}+p\int_{0}^{t\wedge\tau_{n}}(\lambda_{s})^{p-1}\mu(\lambda_{s})ds+p\sigma\int_{0}^{t\wedge\tau_{n}}(\lambda_{s})^{p-1+k}dW_{s}+\frac{p(p-1)\sigma^{2}}{2}\int_{0}^{t\wedge\tau_{n}}(\lambda_{s})^{p-2+2k}ds.\nonumber
\end{align}

\noindent Given the Lipschitz continuity of the drift function $\mu(x)$, there exists a constant $K>0$ such that $\mu(\lambda_{s})-\mu(0)\leq K\lambda_{s}$. Applying Young's inequality, with $1/m+1/n=1$, we have the following two:
\begin{align}\nonumber
    (\lambda_{s})^{p-1}\mu(0)&\leq \frac{[(\lambda_{s})^{p-1}]^{m}}{m}+\frac{[\mu(0)]^{n}}{n}\xlongequal[n=p]{m=\frac{p}{p-1}}=\frac{(\lambda_{s})^{p}}{\frac{p}{p-1}}+\frac{[\mu(0)]^{p}}{p};\\
    (\lambda_{s})^{p-2+2k}&\leq\frac{[(\lambda_{s})^{p-2+2k}]^{m}}{m}+\frac{1}{n}\xlongequal[n=\frac{p}{2-2k}]{m=\frac{p}{p-2+2k}}
    \frac{(\lambda_{s})^{p}}{\frac{p}{p-2+2k}}+\frac{1}{\frac{p}{2-2k}}.\nonumber
\end{align}
We have, by taking the expectation (note that an It\^{o}'s integral has $0$ expectation), there exist constants $C_{1}, C_{2}$ that do not depend on $n$:
\begin{align}\nonumber
&\mathbb{E}[(\lambda_{t\wedge\tau_{n}})^{p}]\leq (\lambda_{0})^{p}+p\mathbb{E}\Big[\int_{0}^{t\wedge\tau_{n}}(\lambda_{s})^{p-1}\mu(\lambda_{s})ds\Big]+\frac{p(p-1)\sigma^{2}}{2}\mathbb{E}\Big[\int_{0}^{t\wedge\tau_{n}}(\lambda_{s})^{p-2+2k}ds\Big]\\
\leq&(\lambda_{0})^{p}+pK\mathbb{E}\Big[\int_{0}^{t\wedge\tau_{n}}(\lambda_{s})^{p}ds\Big]+p\mathbb{E}\Big[\int_{0}^{t\wedge\tau_{n}}(\lambda_{s})^{p-1}\mu(0)ds\Big]+\frac{p(p-1)\sigma^{2}}{2}\mathbb{E}\Big[\int_{0}^{t\wedge\tau_{n}}(\lambda_{s})^{p-2+2k}ds\Big]\nonumber\\
\leq&(\lambda_{0})^{p}+pK\mathbb{E}\Big[\int_{0}^{t\wedge\tau_{n}}(\lambda_{s})^{p}ds\Big]+p\frac{\mathbb{E}[\int_{0}^{t\wedge\tau_{n}}(\lambda_{s})^{p}ds]}{\frac{p}{p-1}}+p\frac{[\mu(0)]^{p}}{p}\mathbb{E}[t\wedge\tau_{n}]+\frac{p(p-1)\sigma^{2}}{2}\Bigg(\frac{\mathbb{E}[\int_{0}^{t\wedge\tau_{n}}(\lambda_{s})^{p}ds]}{\frac{p}{p-2+2k}}+\frac{\mathbb{E}[t\wedge\tau_{n}]}{\frac{p}{2-2k}}\Bigg)\nonumber\\
=&(\lambda_{0})^{p}+\Big[[\mu(0)]^{p}+\frac{p(p-1)\sigma^{2}}{2}\frac{1}{\frac{p}{2-2k}}\Big]\mathbb{E}[t\wedge\tau_{n}]+\Big[pK+\frac{p}{\frac{p}{p-1}}+\frac{p(p-1)\sigma^{2}}{2}\frac{1}{\frac{p}{p-2+2k}}\Big]\mathbb{E}\Big[\int_{0}^{t\wedge\tau_{n}}(\lambda_{s})^{p}ds\Big]\nonumber\\
=&(\lambda_{0})^{p}+\Big[[\mu(0)]^{p}+(p-1)(1-k)\sigma^{2}\Big]\mathbb{E}[t\wedge\tau_{n}]+\Big[pK+p-1+\frac{1}{2}(p-1)(p-2+2k)\sigma^{2}\Big]\mathbb{E}\Big[\int_{0}^{t\wedge\tau_{n}}(\lambda_{s})^{p}ds\Big]\nonumber\\
\leq&C_{1}+C_{2}\mathbb{E}\Big[\int_{0}^{t\wedge\tau_{n}}(\lambda_{s})^{p}ds\Big]\xlongequal{\text{Fubini-Tonelli}}C_{1}+C_{2}\int_{0}^{t}\mathbb{E}[\mathbbm{1}_{\{s\leq \tau_{n}\}}(\lambda_{s})^{p}]ds\leq C_{1}+C_{2}\int_{0}^{t}\mathbb{E}[(\lambda_{s\wedge\tau_{n}})^{p}]ds.\nonumber
\end{align}
By Gr{\"{o}}nwall's inequality, $\sup_{n}\sup_{0\leq t\leq T}\mathbb{E}[(\lambda_{t\wedge\tau_{n}})^{p}]
\leq C_{1}e^{C_{2}T}<\infty$, for $p\geq2$. Thus, for every $p\geq 2$, the stopped positive moments are bounded uniformly in $n$ and $t\in[0,T]$.\\

\noindent \textbf{Secondly, we consider the negative moments with $p>0$}, namely moments of $(\lambda_{t})^{-p}$. We use the same stopping times $\tau_{n}$. By It\^{o}'s lemma, we have
\begin{align}\nonumber
&(\lambda_{t\wedge\tau_{n}})^{-p}=(\lambda_{0})^{-p}+\int_{0}^{t\wedge\tau_{n}}(-p)(\lambda_{s})^{-(p+1)}d\lambda_{s}+\frac{1}{2}\int_{0}^{t\wedge\tau_{n}}p(p+1)(\lambda_{s})^{-(p+2)}(d\lambda_{s})^{2}\\
=&(\lambda_{0})^{-p}-p\int_{0}^{t\wedge\tau_{n}}\frac{\mu(\lambda_{s})}{(\lambda_{s})^{p+1}}ds-p\sigma\int_{0}^{t\wedge\tau_{n}}\frac{(\lambda_{s})^{k}}{(\lambda_{s})^{p+1}}dW_{s}+\frac{p(p+1)\sigma^{2}}{2}\int_{0}^{t\wedge\tau_{n}}\frac{(\lambda_{s})^{2k}}{(\lambda_{s})^{p+2}}ds.\nonumber
\end{align}

\noindent Given the Lipschitz continuity of the drift function $\mu(x)$, there exists a constant $K>0$ such that $\mu(\lambda_{s})-\mu(0)\geq -K\lambda_{s}$, that is, $-\mu(\lambda_{s})\leq K\lambda_{s}-\mu(0)$. We have, by taking the expectation (note that an It\^{o}'s integral has $0$ expectation)
\begin{align}\nonumber
\mathbb{E}[(\lambda_{t\wedge\tau_{n}})^{-p}]&=(\lambda_{0})^{-p}+p\mathbb{E}\Big[\int_{0}^{t\wedge\tau_{n}}\frac{-\mu(\lambda_{s})}{(\lambda_{s})^{p+1}}ds\Big]+\frac{1}{2}p(p+1)\sigma^{2}\mathbb{E}\Big[\int_{0}^{t\wedge\tau_{n}}\frac{1}{(\lambda_{s})^{2(1-k)+p}}ds\Big]\\
&\leq(\lambda_{0})^{-p}+p\mathbb{E}\Big[\int_{0}^{t\wedge\tau_{n}}\frac{K\lambda_{s}-\mu(0)}{(\lambda_{s})^{p+1}}ds\Big]+\mathbb{E}\Big[\int_{0}^{t\wedge\tau_{n}}\frac{p(p+1)\sigma^{2}}{2(\lambda_{s})^{2(1-k)+p}}ds\Big]\nonumber\\
&\leq(\lambda_{0})^{-p}+pK\int_{0}^{t}\mathbb{E}\Big[(\lambda_{s\wedge\tau_{n}})^{-p}\Big]ds+\mathbb{E}\Big[\int_{0}^{t\wedge\tau_{n}}\Big(\frac{p(p+1)\sigma^{2}}{2(\lambda_{s})^{2(1-k)+p}}-\frac{p\mu(0)}{(\lambda_{s})^{p+1}}\Big)ds\Big].\nonumber
\end{align}

\noindent Let $l(x)=\frac{p(p+1)\sigma^{2}}{2x^{2(1-k)+p}}-\frac{p\mu(0)}{x^{p+1}}$, so $l'(x)=\frac{\partial l(x)}{\partial x}=\frac{-p(p+1)\sigma^{2}(2-2k+p)}{2x^{3-2k+p}}+\frac{p(p+1)\mu(0)}{x^{p+2}}$, the extreme point will be the value of $x$ (say $x^{*}$) which makes $\frac{p(p+1)\sigma^{2}(2-2k+p)}{2x^{3-2k+p}}=\frac{p(p+1)\mu(0)}{x^{p+2}}$, that is, $\frac{\sigma^{2}(2-2k+p)}{2x^{1-2k}}=\mu(0)$, which means $x^{*}=\big[\frac{\sigma^{2}(2-2k+p)}{2\mu(0)}\big]^{\frac{1}{1-2k}}$. As a result,
\begin{align}\nonumber
    l(x^{*})&=\frac{p(p+1)\sigma^{2}}{2\Big[\frac{\sigma^{2}(2-2k+p)}{2\mu(0)}\Big]^{\frac{2-2k+p}{1-2k}}}-\frac{p\mu(0)}{\Big[\frac{\sigma^{2}(2-2k+p)}{2\mu(0)}\Big]^{\frac{p+1}{1-2k}}}
    =\frac{p(p+1)\sigma^{2}}{2\Big[\frac{\sigma^{2}(2-2k+p)}{2\mu(0)}\Big]^{\frac{2-2k+p}{1-2k}}}-\frac{2p\mu(0)\Big[\frac{\sigma^{2}(2-2k+p)}{2\mu(0)}\Big]}{2\Big[\frac{\sigma^{2}(2-2k+p)}{2\mu(0)}\Big]^{\frac{p+1}{1-2k}+\frac{1-2k}{1-2k}}}\\
    &=\frac{p(p+1)\sigma^{2}-p\sigma^{2}(2-2k+p)}{2\Big[\frac{\sigma^{2}(2-2k+p)}{2\mu(0)}\Big]^{\frac{2-2k+p}{1-2k}}}=\frac{p\sigma^{2}(2k-1)}{2}\Big[\frac{\sigma^{2}(2-2k+p)}{2\mu(0)}\Big]^{\frac{2-2k+p}{2k-1}}\overset{\text{def}}{=}L.\nonumber
\end{align}
Note that $l''(x)=\frac{\partial^{2} l}{\partial x^{2}}=\frac{p(p+1)\sigma^{2}(2-2k+p)(3-2k+p)}{2x^{4-2k+p}}-\frac{p(p+1)(p+2)\mu(0)}{x^{p+3}}$. Therefore, 
\begin{equation}\nonumber
\hspace*{-0.5cm}
    l''(x^{*})=\frac{p(p+1)\sigma^{2}(2-2k+p)(3-2k+p)}{2}\Big[\frac{\sigma^{2}(2-2k+p)}{2\mu(0)}\Big]^{\frac{4-2k+p}{2k-1}}-p(p+1)(p+2)\mu(0)\Big[\frac{\sigma^{2}(2-2k+p)}{2\mu(0)}\Big]^{\frac{p+3}{2k-1}}
    \overset{\text{def}}{=}C_{1}^{*}M^{s_{1}}-C_{2}^{*}M^{s_{2}},\nonumber
\end{equation}
where $C_{1}^{*}\overset{\text{def}}{=}\frac{p(p+1)\sigma^{2}(2-2k+p)(3-2k+p)}{2}$, $C_{2}^{*}\overset{\text{def}}{=}p(p+1)(p+2)\mu(0)$, $M\overset{\text{def}}{=}\frac{\sigma^{2}(2-2k+p)}{2\mu(0)}$, $s_{1}\overset{\text{def}}{=}\frac{4-2k+p}{2k-1}$ and $s_{2}\overset{\text{def}}{=}\frac{p+3}{2k-1}$. Note that $s_{1}-s_{2}=\frac{4-2k+p-p-3}{2k-1}=\frac{1-2k}{2k-1}=-1$, which leads to $l''(x^{*})=M^{s_{2}}(\frac{C_{1}^{*}}{M}-C_{2}^{*})$. We can easily obtain $\frac{C_{1}^{*}}{M}=p(p+1)(3-2k+p)\mu(0)$ and then $\frac{C_{1}^{*}}{M}-C_{2}^{*}=p(p+1)\mu(0)(1-2k)$.\\

\noindent Assuming $\frac{1}{2}<k<1$ makes $2-2k>0$, so $M>0$ and thus $M^{s_{2}}>0$; Assuming $\frac{1}{2}<k<1$ also makes $1-2k<0$, so $\frac{C_{1}^{*}}{M}-C_{2}^{*}<0$. As a result, $l''(x^{*})<0$. Moreover, $x^{*}$ is the unique critical point of $l$ on $(0,\infty)$, and
\begin{equation}\nonumber
    \lim_{x\downarrow0}l(x)=-\infty, \qquad \lim_{x\to+\infty}l(x)=0.
\end{equation}
Therefore, $x^{*}$ is the global maximizer, and $l(x)\leq l(x^{*})=:L$ for any $x>0$.
\noindent In summary, we have $\mathbb{E}[(\lambda_{t\wedge\tau_{n}})^{-p}]\leq (\lambda_{0})^{-p}+pK\int_{0}^{t}\mathbb{E}[(\lambda_{s\wedge\tau_{n}})^{-p}]ds+Lt$, and from Gr{\"{o}}nwall's inequality, we finally have $\mathbb{E}[(\lambda_{t\wedge\tau_{n}})^{-p}]\leq[(\lambda_{0})^{-p}+Lt]\text{exp}\{pKt\}$. Consequently,
\begin{equation}\nonumber
    \sup_{n}\sup_{0\leq t\leq T}\mathbb{E}[(\lambda_{t\wedge\tau_{n}})^{-p}]\leq((\lambda_{0})^{-p}+LT)e^{pKT}<\infty, \qquad p>0.
\end{equation}
Using the positive- and negative-moment bounds with $p=2$, together with path continuity, we obtain
\begin{equation}\nonumber
    \mathbb{P}(\tau_{n}\leq T)\leq\frac{1}{n^{2}}\big(\mathbb{E}[(\lambda_{T\wedge\tau_{n}})^{2}]+\mathbb{E}[(\lambda_{T\wedge\tau_{n}})^{-2}]\big)\leq\frac{C_{T}}{n^{2}}.
\end{equation}
Since $T>0$ is arbitrary, $\tau_{n}\uparrow+\infty$ $a.s.$. Fatou's lemma then gives the corresponding fixed-time moment bounds for $\lambda_{t}$.\\

\noindent It remains to establish the moment bounds for the supremum. Fix $p\geq2$, and define
\begin{equation}\nonumber
    M_{t}^{+,n}=p\sigma\int_{0}^{t\wedge\tau_{n}}(\lambda_{s})^{p-1+k}dW_{s},\qquad M_{t}^{-,n}=-p\sigma\int_{0}^{t\wedge\tau_{n}}(\lambda_{s})^{k-p-1}dW_{s}.
\end{equation}
By the BDG inequality, Jensen's inequality, and the preceding bounds for higher-order fixed-time moments,
\begin{align}\nonumber
    \mathbb{E}\Big[\sup_{0\leq t\leq T}\lvert M_{t}^{+,n}\lvert\Big]&\leq C_{p}\mathbb{E}\Big[\Big(\int_{0}^{T\wedge\tau_{n}}(\lambda_{s})^{2p-2+2k}ds\Big)^{1/2}\Big]\leq C_{p}\Big(\int_{0}^{T}\mathbb{E}[(\lambda_{s\wedge\tau_{n}})^{2p-2+2k}]ds\Big)^{1/2}\leq C_{p,T},\\
    \mathbb{E}\Big[\sup_{0\leq t\leq T}\lvert M_{t}^{-,n}\lvert\Big]&\leq C_{p}\mathbb{E}\Big[\Big(\int_{0}^{T\wedge\tau_{n}}(\lambda_{s})^{-(2p+2-2k)}ds\Big)^{1/2}\Big]\leq C_{p}\Big(\int_{0}^{T}\mathbb{E}[(\lambda_{s\wedge\tau_{n}})^{-(2p+2-2k)}]ds\Big)^{1/2}\leq C_{p,T}.\nonumber
\end{align}
Here and below, $C_{p,T}$ may change from line to line but is independent of $n$.\\
Taking the supremum in the stopped It\^{o} formulas and using the preceding drift estimates yields
\begin{align}\nonumber
    \mathbb{E}\Big[\sup_{0\leq t\leq T}(\lambda_{t\wedge\tau_{n}})^{p}\Big]&\leq C_{1}+C_{2}\int_{0}^{T}\mathbb{E}[(\lambda_{s\wedge\tau_{n}})^{p}]ds+\mathbb{E}\Big[\sup_{0\leq t\leq T}\lvert M_{t}^{+,n}\lvert\Big]\leq C_{p,T},\\
    \mathbb{E}\Big[\sup_{0\leq t\leq T}(\lambda_{t\wedge\tau_{n}})^{-p}\Big]&\leq(\lambda_{0})^{-p}+LT+pK\int_{0}^{T}\mathbb{E}[(\lambda_{s\wedge\tau_{n}})^{-p}]ds+\mathbb{E}\Big[\sup_{0\leq t\leq T}\lvert M_{t}^{-,n}\lvert\Big]\leq C_{p,T}.\nonumber
\end{align}
Letting $n\to\infty$ and applying Fatou's lemma gives
\begin{equation}\nonumber
    \mathbb{E}\Big[\sup_{0\leq t\leq T}\lambda_{t}^{p}\Big]+\mathbb{E}\Big[\sup_{0\leq t\leq T}\lambda_{t}^{-p}\Big]<\infty, \qquad p\geq 2.
\end{equation}
For $0<p<2$, the same conclusion follows from
\begin{equation}\nonumber
    x^{p}\leq 1+x^{2}, \qquad x^{-p}\leq 1+x^{-2}, \qquad x>0.
\end{equation}
The case $p=0$ is immediate.\qed\\

\noindent \textit{Proof} \textbf{of (7)}: It is not difficult to see that $p_{\infty}(x)$ tends to zero faster than any negative power of $x$ as $x$ approaches $+\infty$, and thus for arbitrary $q$, the integrand, i.e. $x^{q}$ times $p_{\infty}(x)$, will always tend to zero no matter what value $q$ takes (Note that using L'H\^{o}pital's rule will give the same result). Further, if we take an arbitrarily large $q'>1$ to check the limit behavior of:
\begin{equation}\nonumber
    \frac{x^{q} p_{\infty}(x)}{x^{-q'}}=Ax^{q+q'-2k}\text{exp}\bigg\{\frac{2}{\sigma^{2}}\bigg(\frac{ax^{1-2k}}{1-2k}-\frac{bx^{2-2k}}{2-2k}\bigg)\bigg\}\xrightarrow{}0.
\end{equation}
Clearly, $\forall\epsilon>0$, $\int_{\epsilon}^{\infty}x^{q} p_{\infty}(x)dx<+\infty$, since $\int_{\epsilon}^{\infty}x^{-q'}dx<+\infty$ always holds. We also need to check integrability near zero. Since $\frac{1}{2}<k<1$, we have $2k-1>0$ and $2-2k>0$. Consequently, as $x\downarrow0$,
\begin{equation}\nonumber
    x^{q}p_{\infty}(x)\sim A x^{q-2k}\text{exp}\Big\{-\frac{2a}{\sigma^{2}(2k-1)}x^{-(2k-1)}\Big\}\to0.
\end{equation}
for every $q\in\mathbb{R}$, because $a>0$ and the exponential factor dominates any power singularity at zero. Together with continuity on $(0,+\infty)$, this implies $\int_{0}^{\epsilon}x^{q}p_{\infty}(x)dx<+\infty$, $\forall\epsilon>0$. Combining this with the estimate at infinity, we obtain
\begin{equation}\nonumber
    \int_{0}^{\infty}x^{q}p_{\infty}(x)dx<+\infty, \qquad \forall q\in\mathbb{R}.
\end{equation}
Finally, the ergodic theorem implies that for $\mathbb{R}\ni q\neq 0$, $\frac{1}{T}\int_{0}^{T}(\lambda_{t})^{q}dt\xrightarrow[a.s.]{T\xrightarrow{}+\infty}\int_{0}^{\infty}x^{q} p_{\infty}(x)dx$. See also \hyperlink{Andersen2007}{Andersen and Piterbarg (2007)} [Section 2, Proposition 2.1, Proposition 2.2, pages 32-33].\qed\\

\noindent \textit{Proof of }\textbf{Lemma 3.1}:\\
This is another application of \textbf{Theorem 4.2} (Feller's test for explosion) and its generalization - the boundary classification criteria, which describes the boundary behaviors of some diffusion processes of prescribed types. In this case, $\frac{\mu(z)}{[\nu(z)]^{2}}=\frac{a^{*}b^{*}}{\sigma^{*2}z}-\frac{a^{*}}{\sigma^{*2}}$, $\int_{c}^{y}\frac{\mu(z)}{[\nu(z)]^{2}}dz=\frac{a^{*}b^{*}}{\sigma^{*2}}\text{log}\frac{y}{c}-\frac{a^{*}}{\sigma^{*2}}(y-c)$ (integral constants omitted, the same hereinafter), so the scale function for testing is $\psi(x)=\int_{c}^{x}\text{exp}\big\{-\frac{2a^{*}b^{*}}{\sigma^{*2}}\text{log}\frac{y}{c}+\frac{2a^{*}}{\sigma^{*2}}(y-c)\big\}dy=\int_{c}^{x}(\frac{y}{c})^{-\frac{2a^{*}b^{*}}{\sigma^{*2}}}\text{exp}\big\{\frac{2a^{*}}{\sigma^{*2}}(y-c)\big\}dy$. Whether the CIR process touches zero with probability one can be shown by calculating the value of $\lim_{x\to0+}\psi(x)$. As $y\downarrow0$, the exponential factor $\text{exp}\big\{\frac{2a^{*}}{\sigma^{*2}}(y-c)\big\}$ converges to $\text{exp}\{-2a^{*}c/\sigma^{*2}\}>0$. Consequently, the two improper integrals
\begin{equation}\nonumber
    \int_{0}^{c}\Big(\frac{y}{c}\Big)^{-\frac{2a^{*}b^{*}}{\sigma^{*2}}}\text{exp}\Big\{\frac{2a^{*}}{\sigma^{*2}}(y-c)\Big\}dy \quad\text{and}\quad \int_{0}^{c}y^{-\frac{2a^{*}b^{*}}{\sigma^{*2}}}dy
\end{equation}
are either both finite or both infinite. Obviously, when $-\frac{2a^{*}b^{*}}{\sigma^{*2}}\leq-1$, namely when Feller's condition $2a^{*}b^{*}\geq\sigma^{*2}$ is satisfied, the scale function for testing $\psi(x)$ explodes to $-\infty$. One can find a more detailed proof, among others, in \hyperlink{Clark2011}{Clark (2011)} [Chapter 6 \S 3.1, pages 98-104, method B] or \hyperlink{Lamberton2011}{Lamberton and Lapeyre (2011)} [Chapter 6 \S 2 Proposition 6.2.4, page 130]. It is also worth comparing this result to the proofs detailed for \textbf{Theorem 1.2}.\\

\noindent We now outline an alternative way of proof following \hyperlink{Clark2011}{Clark (2011)} [Chapter 6 \S 3.1, pages 98-104, method A]. This approach relies on the fact that a CIR process can be transformed into (and studied via) a Bessel process. Let $\textbf{R}_{t}$ satisfy the $d$-dimensional standard Bessel SDE:
\begin{equation}\nonumber
	d\textbf{R}_{t}=\frac{d-1}{2\textbf{R}_{t}}dt+dW_{t}. 
\end{equation}
It has already been an established result that, provided $d\geq 2$, the path of $\mathbf{R}_{t}$ will never hit the origin $\forall t>0$. For integer $d$, one may refer to Proposition 3.22 in \hyperlink{Karatzas2012}{Karatzas and Shreve (2012)} [Chapter 3 \S 3, pages 161-162] for the proof; for the more general case when $d$ is real-valued, one may refer to \hyperlink{GoingJaeschke2003}{G{\"{o}}ing-Jaeschke and Yor (2003)} [Section 2.1, pages 319-321] or \hyperlink{Revuz2013}{Revuz and Yor (2013)} [Chapter XI, \S 1 (1.5) Proposition, pages 442-443], where the result is obtained following an analysis of the scale function and speed measure (see the part \textbf{Feller's boundary classification} in the Appendix for these two concepts). In the analysis, it is shown that the Bessel process:
\begin{equation}\nonumber
	d\textbf{R}_{t}=\frac{1-2\nu}{2\textbf{R}_{t}}dt+dW_{t},
\end{equation}
never reaches the origin if and only if $\nu\leq 0$, which corresponds to $d\geq 2$.\\

\noindent As will be detailed in the proof \textbf{Lemma 3.5} (a) right after, the CIR process $r_{t}$ is related to a squared Bessel process by an exact time–space transformation. To identify the corresponding dimension, we consider the following square-root diffusion with constant drift:
\begin{equation}\nonumber
 	\text{(\RomanNumeralCaps{1})}~dR_{t}=a^{*}b^{*}dt+\sigma^{*}(R_{t})^{\frac{1}{2}}dW_{t},
\end{equation}
The approach taken there is to put $R_{t}\overset{\text{def}}{=}\frac{\sigma^{*2}}{4}(\textbf{R}_{t})^{2}$ and see what the required dimension (which will reveal whether $0$ is an attainable point) is. Applying It\^{o}'s lemma for $f(x)=\frac{1}{4}\sigma^{*2}x^{2}$ gives
\begin{equation}\nonumber
	dR_{t}=\frac{\sigma^{*2}d}{4}dt+\frac{\sigma^{*2}}{2}\textbf{R}_{t}dW_{t}.
\end{equation}
Since $R_{t}=\frac{\sigma^{*2}}{4}(\textbf{R}_{t})^{2}$, we have $\textbf{R}_{t}=\frac{2}{\sigma^{*}}(R_{t})^{\frac{1}{2}}$, so the dynamics of $R_{t}$ become
\begin{equation}\nonumber
 	\text{(\RomanNumeralCaps{2})}~dR_{t}=\frac{\sigma^{*2}d}{4}dt+\sigma^{*}(R_{t})^{\frac{1}{2}}dW_{t}.
\end{equation}
Comparing (\RomanNumeralCaps{1}) with (\RomanNumeralCaps{2}), we see that the diffusion terms coincide and that if $d$ is chosen such that $\frac{1}{4}\sigma^{*2}d=a^{*}b^{*}$, then the drift terms also coincide. Hence,
\begin{equation}\nonumber
 	d=\frac{4a^{*}b^{*}}{\sigma^{*2}}.
\end{equation}
Requiring $d\geq2$ (to ensure $R_{t}$ never reaches $0$) is therefore equivalent to requiring
\begin{equation}\nonumber
 	(\text{Feller's condition})~2a^{*}b^{*}\geq{\sigma^{*2}},
\end{equation}
which guarantees that $R$, started from $R_{0}>0$, never reaches zero. To transfer this conclusion to the CIR process, \textbf{Lemma 3.5} (a) to be explained later gives that
\begin{equation}\nonumber
    r_{t}=e^{-a^{*}t}\operatorname{BESQ}_{(d,r_{0})}\big(g(t)\big), \qquad g(t)=\frac{\sigma^{*2}}{4a^{*}}\big(e^{a^{*}t}-1\big),
\end{equation}
where $r_{0}>0$ and $d=4a^{*}b^{*}/\sigma^{*2}$. Since $a^{*}>0$, the function $g$ is a strictly increasing bijection from $[0,+\infty)$ onto itself. Moreover, $e^{-a^{*}t}>0$ for every finite $t$. Consequently, the CIR process reaches zero in finite time if and only if the corresponding squared Bessel process does. Thus, $r_{t}=0$ is unattainable when Feller's condition is satisfied. On the other hand, if $2a^{*}b^{*}<(\sigma^{*})^{2}$, then $0<d<2$, and the corresponding squared Bessel process, started from $r_{0}>0$, hits zero in finite time almost surely. Since $g$ is a strictly increasing bijection from $[0,+\infty)$ onto itself and $e^{-a^{*}t}>0$ for every finite $t$, the CIR process also hits zero in finite time almost surely. Thus, $\mathbb{P}(\tau_{0}^{r}<+\infty)=1$, which proves \textbf{Lemma 3.1 (2)}.\qed\\

\noindent \textit{Proof of }\textbf{Lemma 3.3}:\\
The CIR model (2.7) is equivalent to:
\begin{equation}\nonumber
    dr_{t}+a^{*}r_{t}dt=a^{*}b^{*}dt+\sigma^{*}(r_{t})^{\frac{1}{2}}dW_{t},
\end{equation}
thus multiplying both sides by $e^{a^{*}t}$ results in:
\begin{align}\nonumber
    e^{a^{*}t}dr_{t}+a^{*}e^{a^{*}t}r_{t}dt&=a^{*}b^{*}e^{a^{*}t}dt+\sigma^{*}e^{a^{*}t}(r_{t})^{\frac{1}{2}}dW_{t}\\
    d(e^{a^{*}t}r_{t})&=a^{*}b^{*}e^{a^{*}t}dt+\sigma^{*}e^{a^{*}t}(r_{t})^{\frac{1}{2}}dW_{t}\nonumber\\
    e^{a^{*}t}r_{t}-r_{0}&=a^{*}b^{*}\int_{0}^{t}e^{a^{*}s}ds+\sigma^{*}\int_{0}^{t}e^{a^{*}s}\big(r_{s}\big)^{\frac{1}{2}}dW_{s}\nonumber\\
    \text{Thus,}~r_{t}=e^{-a^{*}t}r_{0}+b^{*}(1&-e^{-a^{*}t})+\sigma^{*}e^{-a^{*}t}\int_{0}^{t}e^{a^{*}s}\big(r_{s}\big)^{\frac{1}{2}}dW_{s},\nonumber
\end{align}
which is the exact solution (the CIR process) to the CIR model (2.7).\qed\\

\noindent \textit{Proof of} \textbf{Lemma 3.4}:\\
\textbf{Firstly}, we show how Bessel process can be constructed from a series of independent OU processes, and how Bessel process is related to CIR model, which serves as a complement to \textbf{Lemma 3.5}.\\

\noindent We first consider the case $d=\frac{4a^{*}b^{*}}{\sigma^{*2}}\in\{1,2,\ldots\}$, and choose deterministic initial values satisfying $\sum_{i=1}^{d}(Z_{0}^{i})^{2}=\frac{4a^{*}}{\sigma^{*2}}r_{0}$. Suppose that $Z^{1}_{t},\cdots,Z^{d}_{t}$ are $d$ independent OU processes:
\begin{equation}\nonumber
    dZ^{i}_{t}=-\frac{1}{2}a^{*}Z^{i}_{t}dt+(a^{*})^{\frac{1}{2}}dB^{i}_{t},
\end{equation}
where $B^{i}_{t}$ are independent standard Wiener processes. Consider the squared radius $R_{t}\overset{\text{def}}{=}\sum_{i=1}^{d}(Z^{i}_{t})^{2}$ in $\mathbb{R}^{d}$ of the vector process $Z^{i}_{t}$. Note that $d(Z^{i}_{t})=-\frac{1}{2}a^{*}Z^{i}_{t}dt+(a^{*})^{\frac{1}{2}}dB^{i}_{t}$, which leads to $(dZ^{i}_{t})^{2}=a^{*}dt$. By It\^{o}'s lemma:
\begin{equation}\nonumber
    d(Z^{i}_{t})^{2}=2Z^{i}_{t}\Big(-\frac{1}{2}a^{*}Z^{i}_{t}dt+(a^{*})^{\frac{1}{2}}dB^{i}_{t}\Big)+a^{*}dt=-a^{*}(Z^{i}_{t})^{2}dt+2(a^{*})^{\frac{1}{2}}Z^{i}_{t}dB^{i}_{t}+a^{*}dt.
\end{equation}
Consequently:
\begin{align}\nonumber
    dR_{t}&=\sum_{i=1}^{d}\big[(2Z^{i}_{t}dZ^{i}_{t})+(dZ^{i}_{t})^{2}\big]=-a^{*}\sum_{i=1}^{d}(Z^{i}_{t})^{2}dt+2\sum_{i=1}^{d}Z^{i}_{t}(a^{*})^{\frac{1}{2}}dB^{i}_{t}+da^{*}dt=a^{*}(d-R_{t})dt+(4a^{*}R_{t})^{\frac{1}{2}}dW_{t}.\nonumber
\end{align}
where $W_{t}$ is another one-dimensional Wiener process. The process $R_{t}$ is a CIR-type process. To obtain a standard squared Bessel process, define $\tau=e^{a^{*}t}-1$ and $\widehat{R}_{\tau}=(1+\tau)R_{\text{log}(1+\tau)/a^{*}}$. By the product rule,
\begin{equation}\nonumber
    d(e^{a^{*}t}R_{t})=a^{*}de^{a^{*}t}dt+2\sqrt{a^{*}e^{a^{*}t}}\sqrt{e^{a^{*}t}R_{t}}dW_{t}.
\end{equation}
Define $\widehat{W}_{\tau}=\int_{0}^{\text{log}(1+\tau)/a^{*}}\sqrt{a^{*}e^{a^{*}s}}dW_{s}$. Its quadratic variation is $\tau$, so it is a Brownian motion with respect to the time-changed filtration. Consequently,
\begin{equation}\nonumber
    d\widehat{R}_{\tau}=dd\tau+2\sqrt{\widehat{R}_{\tau}}d\widehat{W}_{\tau}.
\end{equation}
Thus, $\widehat{R}$ is a standard squared Bessel process of dimension $d$, starting from $\widehat{R}_{0}=R_{0}=\frac{4a^{*}}{\sigma^{*2}}r_{0}$, often denoted as $\mathrm{BESQ}_{(d,\widehat{R}_{0})}$ (as it is in \textbf{Lemma 3.5}). An equivalent representation of this process is to obtain the SDE that the root of $\widehat{R}_{\tau}$ satisfies. Define the standard Bessel process $\mathbf{R}_{\tau}$ as the root of $\widehat{R}_{\tau}$: $\mathbf{R}_{\tau}\overset{\text{def}}{=}(\widehat{R}_{\tau})^{\frac{1}{2}}$. If $f(x)=x^{\frac{1}{2}}$, then $\frac{\partial}{\partial x}f(x)=\frac{1}{2}x^{-\frac{1}{2}}$ and $\frac{\partial^{2}}{\partial x^{2}}f(x)=-\frac{1}{4}x^{-\frac{3}{2}}$. Let $\tau_{0}=\inf\{\tau\geq0:\widehat{R}_{\tau}=0\}$. For $\tau<\tau_{0}$, we have by It\^{o}'s lemma:
\begin{equation}\nonumber
    d\mathbf{R}_{\tau}=\frac{1}{2(\widehat{R}_{\tau})^{\frac{1}{2}}}d\widehat{R}_{\tau}-\frac{1}{2}\frac{1}{4(\widehat{R}_{\tau})^{\frac{3}{2}}}(d\widehat{R}_{\tau})^{2}=\frac{1}{2(\widehat{R}_{\tau})^{\frac{1}{2}}}(dd\tau+2(\widehat{R}_{\tau})^{\frac{1}{2}}d\widehat{W}_{\tau})-\frac{1}{8(\widehat{R}_{\tau})^{\frac{3}{2}}}4\widehat{R}_{\tau}d\tau=\frac{d-1}{2\mathbf{R}_{\tau}}d\tau+d\widehat{W}_{\tau}.\nonumber
\end{equation}
The It\^{o} calculation is understood before the first hitting time of zero. Note that when $d=1$, $\mathbf{R}_{\tau}$ is a reflected Brownian motion. Its global SDE includes a local time term at zero, which does not contribute before the first hitting time of zero.\\

\noindent If we take $\zeta=\frac{2a^{*}}{\sigma^{*2}}$ and $d=\frac{4a^{*}b^{*}}{\sigma^{*2}}$, we have $r_{t}=\frac{R_{t}}{2\zeta}$ and:
\begin{equation}\nonumber
    dr_{t}=a^{*}(b^{*}-r_{t})dt+\sigma^{*}(r_{t})^{\frac{1}{2}}dW_{t}.
\end{equation}
This sum-of-squares representation requires $d$ to be a positive integer and gives a geometric interpretation of the CIR model.\\

\noindent \textbf{Secondly}, we show how the concrete expression of the transition probability density function of the Feller square-root process $r_{t}$ is derived based on the relationship between the CIR process and the Bessel process.\\

\noindent Recall the definition: Let $U_{i}$, $i=1\cdots,d$ be $d$ independent and identically distributed standard normal random variables, and let $\vartheta_{i}$, $i=1,\cdots,d$ be $d$ real numbers with any value. Let $R=\sum_{i=1}^{d}(U_{i}+\vartheta_{i})^{2}$ and $\tilde{\theta}=\sum_{i=1}^{d}\vartheta_{i}^{2}$. Then $R$ has a non-central chi-squared with $d$ degrees of freedom and non-centrality parameter $\tilde{\theta}$. Since $Z^{i}_{t}$ above are all normally distributed with variance $1-e^{-a^{*}t}$, $Z_{t}^{i}\sim\mathcal{N}\big(e^{-a^{*}t/2}Z_{0}^{i},1-e^{-a^{*}t}\big)$), we see that $R\overset{\text{def}}{=}\frac{R_{t}}{1-e^{-a^{*}t}}$ has a non-central chi-squared distribution. Finally we have: for $d=\frac{4a^{*}b^{*}}{\sigma^{*2}}$ (alternatively we can define $\kappa\overset{\text{def}}{=}\frac{2a^{*}b^{*}}{\sigma^{*2}}-1$ so $d=2(\kappa+1)$) and $\omega\overset{\text{def}}{=}\frac{2a^{*}}{\sigma^{*2}(1-e^{-a^{*}t})}$. Then $R\overset{\text{def}}{=}2\omega r_{t}$ will have a non-central chi-squared distribution with $d$ degrees of freedom and the non-centrality parameter $\tilde{\theta}=\frac{e^{-a^{*}t}R_{0}}{1-e^{-a^{*}t}}=2\omega e^{-a^{*}t}r_{0}=2\theta$ and $\theta=\omega e^{-a^{*}t} r_{0}$. See, e.g. \hyperlink{Liptser1997}{Liptser and Shiryaev (1977)} or \hyperlink{Mao2007}{Mao (2007)} for more details.\\

\noindent To establish the transition law for arbitrary real $d>0$, we use the affine transform of the CIR process. Let $s\leq t$ and $\mathbb{P}(r_{t}\in dx\lvert r_{s}=y)=f(s,y;t,x)dx$. The transition density function $f(s,y;t,x)$ satisfies the Kolmogorov backward equation:
\begin{align}\nonumber
    \frac{\partial f}{\partial s}+a^{*}(&b^{*}-y)\frac{\partial f}{\partial y}+\frac{\sigma^{*2}}{2}y\frac{\partial^{2}f}{\partial y^{2}}=0\\
    f(s,y;t,x)dx~\Rightarrow~\delta_{y}(dx)&~\text{(Dirac probability measure at $y$), as }s\uparrow t.\nonumber
\end{align}
We define $g(\tau,u,r_{s})=\mathbb{E}[e^{iur_{t}}\lvert r_{s}]$, where $\tau=t-s$. We know that the CIR process is an affine process, thus:
\begin{equation}\nonumber
    g(\tau,u,r_{s})=\text{exp}\{A(\tau,u)+B(\tau,u)r_{s}\},
\end{equation}
where $A(0,u)=0$, $B(0,u)=iu$. Substituting for $g$ in the Kolmogorov backward equation gives:
\begin{equation}\nonumber
    \frac{\partial g}{\partial\tau}=a^{*}(b^{*}-r_{s})\frac{\partial g}{\partial r_{s}}+\frac{\sigma^{*2}}{2}r_{s}\frac{\partial^{2} g}{\partial r^{2}_{s}}.
\end{equation}
Note that $\frac{\partial g}{\partial \tau}=(\frac{\partial A}{\partial\tau}+r_{s}\frac{\partial B}{\partial \tau})g$, $\frac{\partial g}{\partial r_{s}}=Bg$, $\frac{\partial^{2} g}{\partial r^{2}_{s}}=B^{2}g$. As a result:
\begin{equation}\nonumber
    \frac{\sigma^{*2}}{2} r_{s}B^{2}+a^{*}(b^{*}-r_{s})B-\frac{\partial B}{\partial \tau}r_{s}-\frac{\partial A}{\partial\tau}=0.
\end{equation}
Comparing the constant and linear terms in $r_{s}$ yields
\begin{equation}\nonumber
    \frac{\partial A}{\partial\tau}=a^{*}b^{*}B, \qquad \frac{\partial B}{\partial\tau}+a^{*}B=\frac{\sigma^{*2}}{2}B^{2}.
\end{equation}
which is the Riccati equation. Solving these two differential equations results in the expressions of $A$ and $B$, which then form the expression of $g$ as:
\begin{equation}\nonumber
    g(\tau,u,r_{s})=\Big(1-\frac{iu}{\omega^{*}}\Big)^{-\kappa-1}\text{exp}\Big\{\frac{iue^{-a^{*}\tau}}{1-\frac{iu}{\omega^{*}}}r_{s}\Big\},
\end{equation}
where $\omega^{*}=\frac{2a^{*}}{(1-e^{-a^{*}\tau})\sigma^{*2}}$, $\kappa=\frac{2a^{*}b^{*}}{\sigma^{*2}}-1$. By application of Inverse Fourier Transform (IFT), we obtain the analytical expression of $f(s,y;t,x)$:
\begin{equation}\nonumber
    f(s,y;t,x)=\omega^{*}e^{-\theta^{*}-\gamma^{*}}\Big(\frac{\gamma^{*}}{\theta^{*}}\Big)^{\frac{\kappa}{2}}I_{\kappa}\Big(2(\gamma^{*}\theta^{*})^{\frac{1}{2}}\Big),
\end{equation}
where $\gamma^{*}=\omega^{*}x$ and $\theta^{*}=\omega^{*}e^{-a^{*}\tau}y$. The case $y=0$ is obtained by taking the limit $y\downarrow0$. $I_{\kappa}(\cdot)$ is a modified Bessel function of the first kind of order $\kappa$: $I_{\kappa}(x)=(\frac{x}{2})^{\kappa}\sum_{n=0}^{+\infty}\frac{(x/2)^{2n}}{n!\Gamma(\kappa+n+1)}$. Finally, letting $s=0$ gives the result shown in \textbf{Lemma 3.4}.\qed\\

\noindent \textit{Proof of} \textbf{Lemma 3.5}:\\
\textit{Proof} \textbf{of (a)}\\
Denote $\mathrm{BESQ}_{(d,R_{0})}$ by $R_{t}$. Let $f(t)\overset{\text{def}}{=}e^{-a^{*}t}$, $g(t)\overset{\text{def}}{=}\frac{\sigma^{*2}}{4a^{*}}(e^{a^{*}t}-1)$, for $t\geq0$. Since $\frac{\partial}{\partial t}g(t)=\tfrac{\sigma^{*2}}{4}e^{a^{*}t}>0$, the Dambis-Dubins-Schwarz theorem ensures that the time-changed process
\begin{equation}\nonumber
    \widetilde{W}_{t}\overset{\text{def}}{=}\int_{0}^{t}\Big[\frac{\partial}{\partial u}g(u)\Big]^{-\frac{1}{2}}dW_{g(u)},~t\ge0,
\end{equation}
is a standard Wiener process. Applying It\^{o}s formula to $R_{g(t)}$ yields\footnote{Note that in the drift term $d\frac{\partial}{\partial t}g(t)dt$, the first $d$ denotes the dimension of the squared Bessel process, not the differential symbol.}
\begin{equation}\nonumber
    dR_{g(t)}=d \frac{\partial}{\partial t}g(t)dt+2\big[R_{g(t)}\big]^{\frac{1}{2}}\Big[\frac{\partial}{\partial t}g(t)\Big]^{\frac{1}{2}}d\widetilde{W}_{t}.
\end{equation}
As $f$ depends only on $t$, we obtain
\begin{equation}\nonumber
    dr_{t}=\frac{\partial}{\partial t}f(t)R_{g(t)}dt+f(t)dR_{g(t)}=-a^{*}e^{-a^{*}t}R_{g(t)}dt+e^{-a^{*}t}\Big[d \frac{\partial}{\partial t}g(t)dt+2\big[R_{g(t)}\big]^{\frac{1}{2}}\Big[\frac{\partial}{\partial t}g(t)\Big]^{\frac{1}{2}}d\widetilde{W}_{t}\Big].
\end{equation}
Note $R_{g(t)}=e^{a^{*}t}r_{t}$ and $\frac{\partial}{\partial t} g(t)=\frac{\sigma^{*2}}{4}e^{a^{*}t}$, with $d=\frac{4a^{*}b^{*}}{\sigma^{*2}}$, we split
\begin{align}\nonumber
\text{drift:}\quad & -a^{*}r_{t}+e^{-a^{*}t}d \frac{\partial}{\partial t}g(t)=-a^{*}r_{t}+d\frac{\sigma^{*2}}{4}=-a^{*}r_{t}+a^{*}b^{*},\\
\text{diffusion:}\quad & e^{-a^{*}t}2\big[R_{g(t)}\big]^{\frac{1}{2}}\Big[\frac{\partial}{\partial t}g(t)\Big]^{\frac{1}{2}}=2e^{-a^{*}t}(e^{a^{*}t}r_{t})^{\frac{1}{2}}(\frac{\sigma^{*2}}{4}e^{a^{*}t})^{\frac{1}{2}}=\sigma^{*}(r_{t})^{\frac{1}{2}}.\nonumber
\end{align}
Combining the last two displays gives
\begin{equation}\nonumber
    dr_{t}=a^{*}(b^{*}-r_{t})dt+\sigma^{*}(r_{t})^{\frac{1}{2}}d\widetilde{W}_{t},
\end{equation}
which is exactly the CIR equation with $r_{0}=f(0)R_{0}=R_{0}$.\\

\noindent \textit{Proof} \textbf{of (b)}\\
Let $h(x)=x^{-\delta}$, then $\frac{\partial}{\partial x}h(x)=-\delta x^{-\delta-1}$, $\frac{\partial^{2}}{\partial x^{2}}h(x)=\delta(\delta+1)x^{-\delta-2}$. Let $Z_{t}=h(\eta_{t})$, we have by It\^{o}'s lemma
\begin{align}\nonumber
    dZ_{t}&=\frac{\partial}{\partial x}h(x)\lvert _{x=\eta_{t}}d\eta_{t}+\frac{1}{2}\frac{\partial^{2}}{\partial x^{2}}h(x)\lvert _{x=\eta_{t}}(d\eta_{t})^{2}=-\delta(\eta_{t})^{-\delta-1}(\mu \eta_{t}dt+\gamma(\eta_{t})^{K}dW_{t})+\frac{1}{2}\delta(\delta+1)(\eta_{t})^{-\delta-2}(\gamma^{2}(\eta_{t})^{2K})dt\\
    &=\Big(-\delta\mu(\eta_{t})^{-\delta}+\frac{1}{2}\delta(\delta+1)\gamma^{2}(\eta_{t})^{2K-\delta-2}\Big)dt-\delta\gamma(\eta_{t})^{K-\delta-1}dW_{t}.\nonumber
\end{align}
As $2K-\delta-2=0$, $K-\delta-1=K-2(K-1)-1=1-K=-\frac{\delta}{2}$, we know that $(\eta_{t})^{2K-\delta-2}=1$, $(\eta_{t})^{K-\delta-1}=(Z_{t})^{\frac{1}{2}}$, so
\begin{equation}\nonumber
    dZ_{t}=\Big(\frac{1}{2}\delta(\delta+1)\gamma^{2}-\delta\mu Z_{t}\Big)dt-\delta\gamma(Z_{t})^{\frac{1}{2}}dW_{t},
\end{equation}
which is just the desired result.\\

\noindent \textit{Proof} \textbf{of (c)}\\
This follows by applying the same space and time change as in (a) to the equation in (b), before the first hitting time of zero. Indeed, let $g(t)=\frac{\Gamma^{2}}{4}\int_{0}^{t}e^{\mathbf{B}s}ds$ and $R_{g(t)}=e^{\mathbf{B}t}\eta_{t}^{-\delta}$. The deterministic time change gives squared Bessel dynamics for $R$, up to its first hitting time of zero, with dimension $d=\frac{4\mathbf{A}}{\Gamma^{2}}=\frac{2K-1}{K-1}$. Solving for $\eta_{t}$ and evaluating $g(t)$ gives the stated representation.\qed\\

\noindent \textit{Proof of} \textbf{Lemma 3.6}:\\
Given the analytical expression of the solution, we also have the first and the second moments for $r_{t}$, which are:
\begin{align}\nonumber
    \mathbb{E}[r_{t}]&=r_{0}e^{-a^{*}t}+b^{*}(1-e^{-a^{*}t})+\sigma^{*}e^{-a^{*}t}\mathbb{E}\Big[\int_{0}^{t}e^{a^{*}s}\big(r_{s}\big)^{\frac{1}{2}}dW_{s}\Big]=r_{0}e^{-a^{*}t}+b^{*}(1-e^{-a^{*}t}),
\end{align}
because $\int_{0}^{t}e^{a^{*}s}\big(r_{s}\big)^{\frac{1}{2}}dW_{s}]$ is an It\^{o}'s integral.
\begin{align}\nonumber
    &\text{Var}(r_{t})=\mathbb{E}[r^{2}_{t}]-(\mathbb{E}[r_{t}])^{2}=2\Big(e^{-a^{*}t} r_{0}+b^{*}(1-e^{-a^{*}t})\Big)\sigma^{*}e^{-a^{*}t}\mathbb{E}\Big[\int_{0}^{t}e^{a^{*}s}\big(r_{s}\big)^{\frac{1}{2}}dW_{s}\Big]+\sigma^{*2}e^{-2a^{*}t}\mathbb{E}\Big[\Big(\int_{0}^{t}e^{a^{*}s}\big(r_{s}\big)^{\frac{1}{2}}dW_{s}\Big)^{2}\Big]
    \nonumber\\
    =&\sigma^{*2}e^{-2a^{*}t}\mathbb{E}\Big[\int_{0}^{t}e^{2a^{*}s}r_{s}ds\Big]
    =\sigma^{*2}e^{-2a^{*}t}\int_{0}^{t}e^{2a^{*}s}\mathbb{E}[r_{s}]ds
    =\sigma^{*2}e^{-2a^{*}t}\int_{0}^{t}e^{2a^{*}s}\Big[r_{0}e^{-a^{*}s}+b^{*}(1-e^{-a^{*}s})\Big]ds\nonumber\\
    =&\sigma^{*2}e^{-2a^{*}t}\int_{0}^{t}\Big[r_{0}e^{a^{*}s}+b^{*}(e^{2a^{*}s}-e^{a^{*}s})\Big]ds=\sigma^{*2}e^{-2a^{*}t}\Big[\frac{r_{0}}{a^{*}}(e^{a^{*}t}-1)+\frac{b^{*}}{2a^{*}}(e^{2a^{*}t}-1)-\frac{b^{*}}{a^{*}}(e^{a^{*}t}-1)\Big]\nonumber\\
    =&\frac{r_{0}\sigma^{*2}}{a^{*}}\Big(e^{-a^{*}t}-e^{-2a^{*}t}\Big)+\frac{b^{*}\sigma^{*2}}{2a^{*}}\Big(1-e^{-2a^{*}t}-2e^{-a^{*}t}+2e^{-2a^{*}t}\Big)=\frac{r_{0}\sigma^{*2}}{a^{*}}\Big(e^{-a^{*}t}-e^{-2a^{*}t}\Big)+\frac{b^{*}\sigma^{*2}}{2a^{*}}\Big(1-2e^{-a^{*}t}+e^{-2a^{*}t}\Big)\nonumber\\
    =&\frac{r_{0}\sigma^{*2}}{a^{*}}\Big(e^{-a^{*}t}-e^{-2a^{*}t}\Big)+\frac{b^{*}\sigma^{*2}}{2a^{*}}\Big(1-e^{-a^{*}t}\Big)^{2}.
    \nonumber
\end{align}
For $0\leq t\leq t'$, It\^{o}'s isometry gives
\begin{align}\nonumber
    &\text{Cov}(r_{t},r_{t'})=\mathbb{E}\Big[\big(r_{t}-\mathbb{E}[r_{t}]\big)\big(r_{t'}-\mathbb{E}[r_{t'}]\big)\Big]=\mathbb{E}\Big[\sigma^{*}e^{-a^{*}t}\int_{0}^{t}e^{a^{*}u}\big(r_{u}\big)^{\frac{1}{2}}dW_{u}\sigma^{*}e^{-a^{*}t'}\int_{0}^{t'}e^{a^{*}v}\big(r_{v}\big)^{\frac{1}{2}}dW_{v}\Big]\nonumber\\
    =&\sigma^{*2}e^{-a^{*}(t+t')}\int_{0}^{t}e^{2a^{*}u}\mathbb{E}[r_{u}]du=\sigma^{*2}e^{-a^{*}(t+t')}\int_{0}^{t}e^{2a^{*}u}\Big(r_{0}e^{-a^{*}u}+b^{*}(1-e^{-a^{*}u})\Big)du\nonumber\\
    =&\sigma^{*2}e^{-a^{*}(t+t')}\Big(\int_{0}^{t}(r_{0}-b^{*})e^{a^{*}u}du+\int_{0}^{t}b^{*}e^{2a^{*}u}du\Big)=\sigma^{*2}e^{-a^{*}(t+t')}\Big(\frac{r_{0}-b^{*}}{a^{*}}(e^{a^{*}t}-1)+\frac{b^{*}}{2a^{*}}(e^{2a^{*}t}-1)\Big)\nonumber\\
    =&\sigma^{*2}e^{-a^{*}(t+t')}\Big(\frac{r_{0}-b^{*}}{a^{*}}e^{a^{*}t}-\frac{r_{0}-b^{*}}{a^{*}}+\frac{b^{*}}{2a^{*}}e^{2a^{*}t}-\frac{b^{*}}{2a^{*}}\Big)=\frac{\sigma^{*2}}{a^{*}}e^{-a^{*}t'}\Big(r_{0}-r_{0}e^{-a^{*}t}-b^{*}+b^{*}e^{-a^{*}t}+\frac{b^{*}}{2}e^{a^{*}t}-\frac{b^{*}}{2}e^{-a^{*}t}\Big)\nonumber\\
    =&\frac{\sigma^{*2}}{a^{*}}e^{-a^{*}t'}\Big(r_{0}-r_{0}e^{-a^{*}t}-b^{*}+\frac{b^{*}}{2}e^{-a^{*}t}+\frac{b^{*}}{2}e^{a^{*}t}\Big)=\frac{r_{0}\sigma^{*2}}{a^{*}}\Big(e^{-a^{*}t'}-e^{-a^{*}(t+t')}\Big)+\frac{b^{*}\sigma^{*2}}{2a^{*}}\Big(e^{a^{*}(t-t')}+e^{-a^{*}(t+t')}-2e^{-a^{*}t'}\Big)\nonumber
\end{align}
More generally, fix $t>0$ and let $\iota=\frac{2a^{*}b^{*}}{\sigma^{*2}}=\kappa+1$, $c_{t}=\frac{\sigma^{*2}}{2a^{*}}\big(1-e^{-a^{*}t}\big)$, $y_{t}=r_{0}e^{-a^{*}t}$. By the non-central chi-squared transition law in \textbf{Lemma 3.4}, the random variable $2r_{t}/c_{t}$ has a non-central chi-squared distribution with $2\iota$ degrees of freedom and non-centrality parameter $2y_{t}/c_{t}$. Consequently, its moment generating function gives
\begin{equation}\nonumber
    \mathbb{E}[e^{zr_{t}}]=(1-c_{t}z)^{-\iota}\text{exp}\Big\{\frac{y_{t}z}{1-c_{t}z}\Big\}, \qquad z<c_{t}^{-1}.
\end{equation}
In particular, all nonnegative integer moments of $r_{t}$ are finite. For $\lvert z\lvert<c_{t}^{-1}$, expanding the exponential and then using the generalized binomial series yields
\begin{equation}\nonumber
    \mathbb{E}[e^{zr_{t}}]=\sum_{j=0}^{\infty}\frac{y_{t}^{j}z^{j}}{j!}(1-c_{t}z)^{-(\iota+j)}=\sum_{j=0}^{\infty}\sum_{m=0}^{\infty}\frac{y_{t}^{j}}{j!}\frac{\Gamma(\iota+j+m)}{\Gamma(\iota+j)m!}c_{t}^{m}z^{j+m}.
\end{equation}
These series are absolutely convergent in a neighborhood of zero. Comparing the coefficient of $z^{n}$ with that in
\begin{equation}\nonumber
    \mathbb{E}[e^{zr_{t}}]=\sum_{n=0}^{\infty}\frac{\mathbb{E}[r_{t}^{n}]}{n!}z^{n},\qquad \text{we obtain}\qquad \frac{\mathbb{E}[r_{t}^{n}]}{n!}=\sum_{j=0}^{n}\frac{y_{t}^{j}}{j!}\frac{\Gamma(\iota+n)}{\Gamma(\iota+j)(n-j)!}c_{t}^{n-j}.
\end{equation}
Multiplying both sides by $n!$ proves the stated formula. The case $t=0$ follows by continuity.\qed\\

\noindent \textit{Proof of} \textbf{Lemma 3.7}:\\
In fact, we have already derived and proved the analytical expression of the asymptotic stationary distribution density of CKLS process in (1.1). Here we do it again with a particular focus on the case $k=\frac{1}{2}$. Denote the stationary probability density of the CIR process by $p_{\infty}(x)$, $x>0$. It satisfies the stationary Fokker–Planck–Kolmogorov equation:
\begin{equation}\nonumber
    \frac{d}{dx}\big[\Pi p_{\infty}\big]=\frac{1}{2}\frac{d^{2}}{dx^{2}}\big[\Sigma^{2}p_{\infty}\big],
\end{equation}
where $\Pi$ stands for the drift structure $\Sigma$ stands for the diffusion structure (term), and in our case $\Pi=a^{*}(b^{*}-x)$, $\Sigma=\sigma^{*}x^{\frac{1}{2}}$, respectively, which turns out to be:
\begin{equation}\nonumber
    \frac{d}{dx}\big[a^{*}(b^{*}-x)p_{\infty}\big]=\frac{\sigma^{*2}}{2}\frac{d^{2}}{d x^{2}}[xp_{\infty}].
\end{equation}
For the standard CIR process, the stationary probability current is zero. Integrating the preceding equation once gives:
\begin{equation}\nonumber
   a^{*}(b^{*}-x)p_{\infty}=\frac{\sigma^{*2}}{2}\Big(p_{\infty}+x\frac{dp_{\infty}}{dx}\Big),
\end{equation}
which turns out to be:
\begin{align}\nonumber
    \frac{2a^{*}b^{*}}{\sigma^{*2}}-\frac{2a^{*}}{\sigma^{*2}}x&=1+\frac{x}{p_{\infty}}\frac{dp_{\infty}}{dx}\\
    \frac{\kappa}{x}-\frac{\kappa+1}{b^{*}}&=\frac{d}{dx}\text{log}p_{\infty}\nonumber\\
    \kappa\text{log}x-\frac{\kappa+1}{b^{*}}x+C&=\text{log}p_{\infty}\nonumber\\
    \text{Thus, }p_{\infty}\propto x^{\kappa}e&^{-\frac{\kappa+1}{b^{*}}x}.\nonumber
\end{align}
Obviously, ranging over $[0,+\infty)$, $p_{\infty}$ is the asymptotic stationary probability density function of the gamma type (parameters $\kappa+1$ and $\frac{\kappa+1}{b^{*}}$), which means:
\begin{equation}\nonumber
    p_{\infty}=f(x\lvert a^{*},b^{*},\sigma^{*})=\frac{(\frac{\kappa+1}{b^{*}})^{\kappa+1}}{\Gamma(\kappa+1)}x^{\kappa}\text{exp}\Big\{-\frac{\kappa+1}{b^{*}}x\Big\},~x\in(0,+\infty).
\end{equation}
This gives the invariant Gamma density. Convergence in distribution to this invariant law follows by taking $t\to+\infty$ in the characteristic function obtained in the proof of \textbf{Lemma 3.4}.\qed\\

\noindent \textbf{Proofs concerning Novikov's and Kazamaki's conditions}\\
Fix a finite time horizon $T>0$. Let $\theta$ be progressively measurable and satisfy $\int_{0}^{T}(\theta_{s})^{2}ds<\infty$ almost surely. For convenience, write
\begin{equation}\nonumber
    N_{t}=\int_{0}^{t}\theta_{s}dW_{s},\qquad A_{t}=\langle N\rangle_{t}=\int_{0}^{t}(\theta_{s})^{2}ds,\qquad Z_{t}=\mathcal{E}(N)_{t}=\text{exp}\Big\{N_{t}-\frac{1}{2}A_{t}\Big\}.
\end{equation}
Since $Z$ is a nonnegative local martingale starting from $1$, it is a supermartingale and hence $\mathbb{E}[Z_{t}]\leq 1$. Suppose that Novikov's condition holds: $C_{T}:=\mathbb{E}\big[\text{exp}\{\frac{1}{2}A_{T}\}\big]<\infty$. To show that $\mathbb{E}[Z_{T}]=1$, fix $c\in(0,1)$ and define
\begin{equation}\nonumber
    Z_{t}^{(c)}=\mathcal{E}(cN)_{t}=\text{exp}\Big\{cN_{t}-\frac{c^{2}}{2}A_{t}\Big\}, \qquad \tau_{n}=\inf\{t\in[0,T]:\lvert N_{t}\lvert\geq n\}\wedge T,
\end{equation}
where $\inf\varnothing=\infty$. The stopped process $Z^{(c)}_{\cdot\wedge\tau_{n}}$ is bounded above by $e^{cn}$, and is therefore a true martingale. In particular, $\mathbb{E}[Z_{\tau_{n}}^{(c)}]=1$. Using the identity $Z_{\tau_{n}}^{(c)}
=Z_{\tau_{n}}^{c}\text{exp}\big\{\frac{c(1-c)}2A_{\tau_{n}}\big\}$, H\"{o}lder's inequality gives, for every $B\in\mathcal{F}_{T}$,
\begin{equation}\nonumber
    \mathbb{E}\big[Z_{\tau_{n}}^{(c)}\mathbbm{1}_{B}\big]\leq\big(\mathbb{E}[Z_{\tau_{n}}]\big)^{c}\Big(\mathbbm{E}\big[e^{cA_{\tau_{n}}/2}\mathbf1_{B}\big]\Big)^{1-c}\leq\Big(\mathbb{E}\big[e^{A_{T}/2}\mathbbm{1}_{B}\big]\Big)^{1-c}.
\end{equation}
Since $e^{A_{T}/2}$ is integrable, this bound establishes uniform integrability of $\{Z_{\tau_{n}}^{(c)}\}_{n\geq 1}$. By continuity, $Z_{\tau_{n}}^{(c)}\to Z_{T}^{(c)}$ almost surely, and consequently $\mathbb{E}[Z_{T}^{(c)}]=1$. Applying H\"{o}lder's inequality once more, now at time $T$, yields
\begin{equation}\nonumber
    1=\mathbb{E}\big[Z_{T}^{(c)}\big]\leq\big(\mathbb{E}[Z_{T}]\big)^{c}\Big(\mathbb{E}[e^{cA_{T}/2}]\Big)^{1-c}\leq\big(\mathbb{E}[Z_{T}]\big)^{c}(C_{T})^{1-c}.
\end{equation}
Letting $c\uparrow1$, we obtain $\mathbb{E}[Z_{T}]\geq 1$. Together with the supermartingale inequality, this gives $\mathbb{E}[Z_{T}]=1$. Thus $Z$ is a true martingale on $[0,T]$, and $Z_{t}=\mathbb{E}[Z_{T}\mid\mathcal{F}_{t}]$, $0\leq t\leq T$. In particular, it is uniformly integrable on this finite interval.\\

\noindent Novikov's criterion is sufficient but by no means necessary; it also fails to distinguish the sign of the stochastic integral $N_{t}$. For instance, one can have $\mathbb{E}\big[\mathcal{E}(\theta)_{T}\big]=1$ while $\mathbb{E}\big[\mathcal{E}(-\theta)_{T}\big]<1$. For the following example, take $T=1$. Let $W_{t}$, $t\geq 0$ be a standard one-dimensional Wiener process and set $T_{1}\overset{\text{def}}{=}\inf\{t>0: W_{t}=1\}$. Define the time-change $\tau_{t}$, $t\geq 0$ as $\tau_{t}\overset{\text{def}}{=}\frac{t}{1-t}\wedge T_{1}$ if $t<1$; $\tau_{t}\overset{\text{def}}{=}T_{1}$ if $t\geq 1$. Then $\widehat{N}_{t}:=W_{\tau_{t}}$ is a continuous local martingale with respect to the time-changed filtration $\mathcal{G}_{t}:=\mathcal{F}_{\tau_{t}}$.\\

\noindent Since $\widehat{N}_{t}\leq 1$, we have $0<e^{\widehat{N}_{t}/2}\le e^{1/2}$, $0\leq t\leq 1$. By It\^{o}'s formula, $e^{\widehat{N}/2}$ is a local submartingale. Localization and bounded convergence therefore show that it is a submartingale with respect to $\mathcal{G}_{t}$. Hence, Kazamaki's condition holds on $[0,1]$. On the other hand, $\langle\widehat{N}\rangle_{1}=T_{1}$. The first hitting time $T_{1}$ has density $f_{T_{1}}(u)=\frac{1}{\sqrt{2\pi u^{3}}}\text{exp}\{-\frac{1}{2u}\big\}$, $u>0$. Consequently,
\begin{equation}\nonumber
    \mathbb{E}\big[e^{\langle\widehat{N}\rangle_{1}/2}\big]=\mathbb{E}[e^{T_{1}/2}]=\frac{1}{\sqrt{2\pi}}\int_{0}^{\infty}u^{-3/2}\text{exp}\big(\frac {u}{2}-\frac{1}{2u}\big)du=+\infty.
\end{equation}
Thus, Novikov's condition fails on $[0,1]$. Moreover, Kazamaki's criterion gives $\mathbb{E}[\mathcal{E}(\widehat{N})_{1}]=1$, whereas $\widehat{N}_{1}=1$ implies $\mathbb{E}[\mathcal{E}(-\widehat{N})_{1}]=\mathbb{E}[e^{-1-T_{1}/2}]\leq e^{-1}<1$. This also illustrates the difference between the martingale properties of the two stochastic exponentials, although Novikov's condition is unchanged by replacing $\widehat{N}$ with $-\widehat{N}$.\\

\noindent Returning to an arbitrary finite time horizon $T>0$, Kazamaki's condition requires that $(e^{N_{t}/2})_{0\leq t\leq T}$ be a submartingale. To see that Novikov's condition implies this condition, suppose that $\mathbb{E}[e^{A_{T}/2}]<\infty$. Since this assumption implies $\mathbb{E}[A_{T}]<\infty$, the stochastic integral $N$ is a square-integrable martingale on $[0,T]$. Moreover, for every $t\in[0,T]$, $e^{N_{t}/2}=Z_{t}^{\frac{1}{2}}e^{A_{t}/4}.$ Applying the Cauchy–Schwarz inequality and using the supermartingale property of $Z$, we obtain
\begin{equation}\nonumber
    \mathbb{E}[e^{N_{t}/2}]\leq\big(\mathbb{E}[Z_{t}]\big)^{\frac{1}{2}}\big(\mathbb{E}[e^{A_{t}/2}]\big)^{\frac{1}{2}}\leq\big(\mathbb{E}[e^{A_{T}/2}]\big)^{\frac{1}{2}}<\infty.
\end{equation}
Consequently, conditional Jensen's inequality gives, for $0\leq s\leq t\leq T$,
\begin{equation}\nonumber
    \mathbb{E}[e^{N_{t}/2}\mid\mathcal{F}_{s}]\geq\text{exp}\Big\{\frac{1}{2}\mathbb{E}[N_{t}\mid\mathcal{F}_{s}]\Big\}=e^{N_{s}/2}.
\end{equation}
Thus $(e^{N_{t}/2})_{0\leq t\leq T}$ is a submartingale, which verifies Kazamaki's condition. Hence, Novikov's condition implies Kazamaki's condition, whereas the converse does not hold in general.\\

\noindent We now prove that Kazamaki's condition implies the martingale property of the Dol\'{e}ans–Dade exponential. Suppose that $Y_{t}:=e^{N_{t}/2}$, $0\leq t\leq T$, is a submartingale, and set $K_{T}=\mathbb{E}[Y_{T}]<\infty$. Fix $c\in(0,1)$, and define $p=\frac{1}{c(2-c)}>1$ and $\alpha=\frac{c}{2-c}\in(0,1)$. With $Z^{(c)}=\mathcal{E}(cN)$ as above, a direct calculation gives $\big(Z_{t}^{(c)}\big)^{p}=Z_{t}^{\alpha}Y_{t}^{1-\alpha}$. Let $\tau_{n}=\inf\{t\in[0,T]:\lvert N_{t}\lvert\geq n\}\wedge T$. As before, $Z^{(c)}_{\cdot\wedge\tau_{n}}$ is a bounded martingale, so that $\mathbb{E}[Z_{\tau_{n}}^{(c)}]=1$. By H\"{o}lder's inequality and optional sampling,
\begin{equation}\nonumber
    \mathbb{E}\big[\big(Z_{\tau_{n}}^{(c)}\big)^{p}\big]\leq\big(\mathbb{E}[Z_{\tau_{n}}]\big)^{\alpha}\big(\mathbb{E}[Y_{\tau_{n}}]\big)^{1-\alpha}\leq K_{T}^{1-\alpha}.
\end{equation}
Thus $\{Z_{\tau_{n}}^{(c)}\}_{n\geq 1}$ is bounded in $\mathcal{L}^{p}$, with $p>1$, and is therefore uniformly integrable. Letting $n\to\infty$, we obtain $\mathbb{E}[Z_{T}^{(c)}]=1$. Applying Jensen's and H\"{o}lder's inequalities at time $T$, we further obtain
\begin{equation}\nonumber
    1=\big(\mathbb{E}\big[Z_{T}^{(c)}]\big)^{p}\leq\mathbb{E}\big[\big(Z_{T}^{(c)}\big)^{p}\big]\leq\big(\mathbb{E}[Z_{T}]\big)^{\alpha}K_{T}^{1-\alpha}.
\end{equation}
As $c\uparrow1$, we have $\alpha\uparrow1$, and hence $\mathbb{E}[Z_{T}]\geq1$. Since $Z$ is a nonnegative local martingale starting from $1$, the reverse inequality also holds. Therefore, $\mathbb{E}[Z_{T}]=1$. It follows that $Z=\mathcal{E}(N)$ is a true martingale, and is uniformly integrable on $[0,T]$. This completes the proof of Kazamaki's criterion.\qed\\

\noindent \textbf{Feller's boundary classification}\\
We give a concise overview of \hyperlink{Feller1952}{Feller (1952)} boundary classification for one-dimensional SDEs, drawing on the clearer exposition found in \hyperlink{Karatzas2012}{Karatzas and Shreve (2012)}.\\

\noindent Consider the SDE
\begin{equation}\nonumber
 	dZ_{t}=\mu(Z_{t})dt+\nu(Z_{t})dW_{t},
\end{equation}
with a fixed interval $(l,r)$. To classify its boundaries, we introduce the scale density $\omega(x)$ and the speed density $\pi(x)$:
\begin{equation}\nonumber
  \omega(z)\overset{\text{def}}{=}\text{exp}\Big\{-2\int^{z}\frac{\mu(u)}{[\nu(u)]^{2}}du\Big\},
  \qquad
  \pi(z)\overset{\text{def}}{=}\frac{2}{\omega(z)[\nu(z)]^{2}}.
\end{equation}
With these, define four auxiliary set-functions:
\begin{align}\nonumber
  \psi[x,y]&\overset{\text{def}}{=} \int_{x}^{y}\omega(z)dz,\qquad
  \psi(l,y]\overset{\text{def}}{=} \lim_{x\to l+}\psi[x,y],\qquad
  \psi[x,r)\overset{\text{def}}{=} \lim_{y\to {r-}}\psi[x,y];
\\
  \xi[x,y]&\overset{\text{def}}{=} \int_{x}^{y}\pi(z)dz,\qquad
  \xi(l,y]\overset{\text{def}}{=} \lim_{x\to l+}\xi[x,y],\qquad
  \xi[x,r)\overset{\text{def}}{=} \lim_{y\to r-}\xi[x,y];\nonumber
\\
  \phi(l)&\overset{\text{def}}{=} \int_{l}^{x}\psi(l,y]\pi(y)dy\overset{*}{=}\int_{l}^{x}\xi[z,x]\omega(z)dz,\qquad
  \phi(r)\overset{\text{def}}{=}
  \int_{x}^{r}\psi[z,r)\pi(z)dz\overset{*}{=}\int_{x}^{r}\xi[x,y]\omega(y)dy;\nonumber
\\
  \Phi(l)&\overset{\text{def}}{=}
  \int_{l}^{x}\psi[z,x]\pi(z)dz\overset{*}{=}\int_{l}^{x}\xi(l,y]\omega(y)dy,\qquad
  \Phi(r)\overset{\text{def}}{=}
  \int_{x}^{r}\psi[x,y)\pi(y)dy\overset{*}{=}\int_{x}^{r}\xi[z,r]\omega(z)dz.\nonumber
\end{align}
The equalities marked by $*$ follow from Tonelli's theorem, since the integrands are nonnegative; both sides are allowed to take the value $+\infty$. The boundary classification depends on the behavior of the above functions. For an endpoint $e$, one distinguishes four cases:
\begin{align}\nonumber
  1. &\text{ regular,  if } \phi(e) \text{ and } \Phi(e) \text{ are both finite};\\
  2. &\text{ exit, if } \phi(e) \text{ is finite and } \Phi(e) \text{ is infinite};\nonumber\\
  3. &\text{ entrance, if } \phi(e)=\infty \text{ is infinite and } \Phi(e) \text{ is finite};\nonumber\\
  4. &\text{ natural, if } \phi(e)=\infty \text{ and } \Phi(e) \text{ are both infinite}.\nonumber
\end{align}
For a regular boundary, the interior coefficients alone do not determine the behavior at the boundary; an appropriate boundary condition must be specified.\\

\noindent An exit boundary can be reached from inside the domain with positive probability, but the process cannot enter the interior when started at that boundary. Conversely, an entrance boundary cannot be hit from the interior, yet one may start the process at the entrance point. A natural boundary cannot be reached in finite time from the interior, nor can the process be started there.\\

\noindent For a regular boundary, possible reflecting extensions include:
\begin{align}\nonumber
  1.~&\text{instantaneous reflection, under which the process spends zero Lebesgue time at the boundary};\\
  2.~&\text{sticky reflection, under which the process can spend positive Lebesgue time at the boundary}.\nonumber
\end{align}
For a sticky reflecting extension, the speed measure has an atom at the boundary; this is not determined by the endpoint value of the interior speed density $\pi$. Killing at the boundary is specified separately through the boundary condition and cannot be inferred from the finiteness of $\pi(e)$.

\end{document}